\documentclass[11pt]{article}

\usepackage[margin=1in]{geometry}

\usepackage[normalem]{ulem}
\usepackage[dvipsnames]{xcolor}
\usepackage[T1]{fontenc}
\usepackage{lmodern}
\usepackage[utf8]{inputenc} 

\usepackage{algorithm}
\usepackage{algpseudocode}
\usepackage{titlesec}

\titleformat*{\section}{\bfseries}

\titleformat*{\subsection}{\bfseries}

\titleformat*{\subsubsection}{\bfseries} 
\usepackage{dsfont}
\usepackage{caption}
\usepackage{subcaption}

\usepackage[dvipsnames]{xcolor}
\usepackage{graphicx} 
\usepackage{pgfplots}
\pgfplotsset{compat=1.18}
\providecolor{mediationtuningred}{RGB}{255,0,0}
\providecolor{mediationtuningdarkblue}{RGB}{0,0,139}
\providecommand{\mediationtuningaxisfont}{\footnotesize}
\providecommand{\mediationtuningtickfont}{\scriptsize}
\pgfplotsset{mediation tuning/.style={
  width=\linewidth,
  height=0.75\linewidth,
  scale only axis=false,
  axis lines=box,
  axis line style={black,line width=0.45pt},
  tick style={black,line width=0.35pt},
  tick align=outside,
  major tick length=1.8pt,
  grid=none,
  title={},
  label style={font=\mediationtuningaxisfont},
  tick label style={font=\mediationtuningtickfont},
  clip=true
}}

\usepgfplotslibrary{groupplots}
\usetikzlibrary{calc}
\providecommand{\mainqqaxisfont}{\normalfont\small}
\providecommand{\mainqqtickfont}{\normalfont\footnotesize}
\providecommand{\mainqqpanelfont}{\normalfont\small}
\definecolor{mainqqblue}{RGB}{0,0,139}
\definecolor{mainqqred}{RGB}{255,0,0}
\pgfplotsset{main qq/.style={
  compat=1.18, scale only axis, axis lines=box,
  axis line style={black,line width=0.45pt},
  tick style={black,line width=0.35pt}, tick align=outside,
  major tick length=1.8pt, tick pos=left,
  tick label style={font=\mainqqtickfont,inner sep=1pt},
  scaled ticks=false, /pgf/number format/1000 sep={},
  enlargelimits=false, grid=none, clip=true,
  xlabel={}, ylabel={}, title={}
}}

\usepgfplotslibrary{groupplots}
\usetikzlibrary{calc}
\providecommand{\univqqaxisfont}{\normalfont\small}
\providecommand{\univqqtickfont}{\normalfont\footnotesize}
\providecommand{\univqqpanelfont}{\normalfont\small}
\definecolor{univqqblue}{RGB}{0,0,139}
\definecolor{univqqred}{RGB}{255,0,0}
\pgfplotsset{univ qq/.style={
  compat=1.18, scale only axis, axis lines=box,
  axis line style={black,line width=0.45pt},
  tick style={black,line width=0.35pt}, tick align=outside,
  major tick length=1.8pt, tick pos=left,
  tick label style={font=\univqqtickfont,inner sep=1pt},
  scaled ticks=false, /pgf/number format/1000 sep={},
  enlargelimits=false, grid=none, clip=true,
  xlabel={}, ylabel={}, title={}
}}

\usepackage{amsthm,amsmath,amsfonts,amssymb}
\usepackage{mathrsfs} 
\usepackage{verbatim} 
\usepackage{natbib}
\usepackage{bbm} 

\usepackage[plainpages=false,pdfpagelabels,colorlinks=true,
linkcolor=red,urlcolor=red,citecolor=blue]{hyperref}
\usepackage[shortlabels]{enumitem}
\newcommand{\bea}{\bar{E}_a}

\numberwithin{equation}{section}
\newtheorem{theorem}{Theorem}[section]
\newtheorem{corollary}{Corollary}[section]
\newtheorem{lemma}{Lemma}[section]
\newtheorem{assumption}{Assumption}
\newtheorem{proposition}{Proposition}[section]
\newtheorem{definition}{Definition}

\newtheorem{remark}{Remark}[section]

\def\ddefloop#1{\ifx\ddefloop#1\else\ddef{#1}\expandafter\ddefloop\fi}
\def\ddef#1{\expandafter\def\csname c#1\endcsname{\ensuremath{\mathcal{#1}}}}
\ddefloop ABCDEFGHIJKLMNOPQRSTUVWXYZ\ddefloop
\def\ddef#1{\expandafter\def\csname s#1\endcsname{\ensuremath{\mathsf{#1}}}}
\ddefloop ABCDEFGHIJKLMNOPQRSTUVWXYZ\ddefloop
\def\ddef#1{\expandafter\def\csname b#1\endcsname{\ensuremath{\mathbf{#1}}}}
\ddefloop ABCDEFGHIJKLMNOPQRSTUVWXYZ\ddefloop
\def\ddef#1{\expandafter\def\csname b#1\endcsname{\ensuremath{\mathbf{#1}}}}
\ddefloop abcdefghijklmnopqrstuvwxyz\ddefloop

\def\ba{\beta_A}
\def\hba{\hat{\beta}_A}

\def\tm{\theta_M}
\def\htm{\hat{\theta}_M}
\def\btm{\breve{\theta}_M}
\def\hej{\hat{\eta}_j}
\def\ej{\eta_j}
\def\tV{\tilde{V}}
\def\tW{\tilde{W}}
\def\ctV{\tilde{\mathcal{V}}}
\def\tT{\tilde{T}}
\def\tS{\tilde{S}}
\def\baj{\beta_{A,j}}
\def\tmj{\theta_{M,j}}

\def\tr{\mathrm{Tr}}
\def\E{\mathbb{E}}
\def\stY{Y_*}
\def\stA{A_*}
\def\stM{M_*}
\def\stC{C_*}
\def\stZ{Z_*}
\def\stW{W_*}
\def\stV{V_*}
\def\n{\mathfrak n}
\def\p{\mathfrak p}
\def\m{\mathfrak m}
\def\fin{f_{1,n}(\kappa_1)}
\def\fingen{f_{1,\mathfrak{n}}(\kappa)}
\def\cin{\mathfrak{c}_{1,n}}

\theoremstyle{plain}

\theoremstyle{remark}

\newcommand{\doo}{{\rm do}}

    \DeclareMathOperator{\diag}{Diag}

\allowdisplaybreaks

\hypersetup{
    colorlinks=true,
    linkcolor=Violet, 
    citecolor=OliveGreen, 
    urlcolor=brown
}

\begin{document}

\title{\textbf{Inference in High-Dimensional Linear Mediation Models}}

\author{%
Samriddha Lahiry\\
Department of Statistics and Data Science\\
National University of Singapore\\
\href{mailto:slahiry@nus.edu.sg}
     {\texttt{slahiry@nus.edu.sg}}
\and
Rajarshi Mukherjee\\
Department of Biostatistics\\ Harvard University\\
\href{mailto:ram521@mail.harvard.edu}
     {\texttt{ram521@mail.harvard.edu}}
\and
Pragya Sur\\
Department of Statistics\\ Harvard University\\
\href{mailto:pragya@fas.harvard.edu}
     {\texttt{pragya@fas.harvard.edu}}
}

\date{}

\maketitle
\begin{abstract}
Mediation analysis plays an important role in psychology, epidemiology, and genomics. Although high-dimensional mediation has been studied extensively, existing inferential methods typically rely on sparsity assumptions that may be difficult to justify when many individually weak mediators collectively produce a substantial effect. We study inference for both the overall indirect effect and the effect transmitted through an individual mediator in a proportional high-dimensional regime, where the number of mediators and the sample size grow at comparable rates, without imposing sparsity on the mediator effects. Standard debiasing procedures for high-dimensional inference commonly rely on design-covariance information or additional structural assumptions on precision-matrix directions. In mediation models, however, the full outcome-design covariance depends on unknown mediator-model coefficients, creating an additional challenge when such information or structure is unavailable.

To address this challenge, we introduce an approximate debiasing procedure based on partial ridge regression and show that the resulting estimator is asymptotically bias-corrected and Gaussian. We also develop a cross-fitted estimator with a smaller asymptotic variance than its sample-split counterpart. In addition, we construct consistent estimators of the limiting variances of both the sample-split and cross-fitted estimators. Together, these results provide feasible inferential procedures for global and pathwise mediation effects.

\end{abstract}

%\setcounter{tocdepth}{2}
%\renewcommand{\contentsname}{}
%\tableofcontents

\section{Introduction}

Mediation analysis quantifies how much of an exposure's effect on an outcome travels through intermediate variables. 
It has long played a central role in psychology, the social
sciences, and epidemiology
\cite{mackinnon2008introduction,hayes2013introduction,
baron1986moderator,imai2010general,vanderweele2016mediation}. In modern multi-omics and imaging studies, however, the number of candidate mediators is often comparable to or larger than the sample size, and their effects may be individually weak yet collectively substantial \cite{clarkboucher2023methods,zhao2023multimodalmediation,li2025weakmediation}. Over the past two decades, statistical methodology for mediation analysis has advanced rapidly to accommodate many mediating and baseline variables 
\cite{huang2014joint,huang2016hypothesis,zhang2016estimating,
sparse_med,dai2022hima2}. This literature has extended statistical inference for low-dimensional functionals in high-dimensional generalized linear models to high-dimensional mediation models, under sparsity and signal strength assumptions. That said, such sparsity assumptions are often hard to justify in the applications above. In this paper, we develop a detailed  mediation analysis in the complementary regime, without sparsity assumptions on the linear structural equation models underlying our mediation framework. To this end, we work in the proportional asymptotics regime, which has gained traction in high-dimensional statistics and machine learning theory over the past decade \cite{donoho2009message,elkaroui_pnas_1,SurPNAS,jiang2025new,hastie2022surprises,celentano_wainwright}.

 Formally, we observe independent and identically distributed quadruples
$(Y_i, A_i, M_i, C_i)_{i=1}^m$, and 
seek to quantify the effect of an exposure $A_i \in \mathbb{R}$ on an outcome
$Y_i \in \mathbb{R}$, potentially mediated by a vector of intermediate variables
$M_i \in \mathbb{R}^p$ (the \emph{mediators}) and measured in the presence of baseline
covariates $C_i$ (the \emph{confounders}). Let $(\stY,\stA,\stM,\stC)$ denote a generic observation distributed as $(Y_i,A_i,M_i,C_i)$. We assume a linear relationship among
$(\stY, \stA, \stM, \stC)$; that is, a Linear Structural Equation Model (LSEM) of the following form:
\begin{align}\stY&=\theta_M^T\stM+\theta_A\stA+\theta_C^T \stC+\delta,\label{y_mod_vec}\\
\stM&=\beta_A \stA+\beta_C \stC+\epsilon,\label{M_mod_vec}
\end{align}
where $\delta\in\mathbb{R}$ and $\epsilon\in\mathbb{R}^p$ are unobserved error terms. Under standard cross-world
counterfactual assumptions \citep{imai2010general,imai2010identification}, $\theta_A+\beta_A^\top\theta_M$, $
\theta_A$, and $\beta_A^\top\theta_M$ correspond, respectively, to the total, direct, and indirect effects of
exposure $A$ on outcome $Y$. Our main goal is to quantify the effect of the exposure on the outcome through
the mediators, without imposing sparsity on $\tm$.

\subsection{Measures of mediation effects}

Substituting \eqref{M_mod_vec} into \eqref{y_mod_vec} shows that the indirect
effect of $A$ on $Y$ through the mediators is $\beta_A^\top\theta_M$. We
formalize this target as follows. 

\begin{definition}[Global mediation effect]
The global mediation effect is defined as
$\gamma := \beta_A^\top\theta_M.
$
When there is no ambiguity, we refer to $\gamma$ simply as the mediation
effect.
\end{definition}

We also consider the contribution associated with a particular mediator.
Without additional assumptions governing the causal relationships among the
mediators, path-specific causal effects are generally not identifiable. Within the LSEM, however, the product
$\beta_{A,j}\theta_{M,j}$ provides a natural mediator-specific target.

\begin{definition}[Pathwise mediation effect]\label{def:path_med}
For a mediator $j$, the pathwise mediation effect is
$\eta_j := \beta_{A,j}\theta_{M,j},$
where $\beta_{A,j}$ and $\theta_{M,j}$ are the $j$th coordinates of
$\beta_A$ and $\theta_M$, respectively.
\end{definition}

In the linear multiple-mediation literature, this product
is also called the \emph{specific indirect effect};
see \cite{preacher2008asymptotic}.
Our primary objective is to establish the limiting distribution of an
estimator of $\gamma$---obtained by approximate debiasing based on partial ridge regression---and develop analogous inference procedures for
$\eta_j$. In the high-dimensional, non-sparse regime, estimating these quantities is particularly challenging because of the dependence structure in the design induced by \eqref{y_mod_vec} and \eqref{M_mod_vec}.

\subsection{High-dimensional mediation and non-sparse models}
Traditional methods for estimating mediation effects target
settings with a modest number of mediators, for which a variety of established
estimation and testing procedures are available
\cite{hayes2013introduction,vanderweele2014mediation,
vanderweele2015explanation}. Modern applications, however, increasingly
involve a large number of potential mediators, often comparable to or
exceeding the sample size. Examples include genomic and epidemiological
studies where DNA-methylation measurements serve as mediators
\cite{Houtepen2016NatCommun,Guo2022JASA,Tobi2018SciAdv,Huang2019AOAS};
neuroscience studies where activity across many brain regions mediates
the relationship between a stimulus and a behavioral response
\cite{chen2015highdim}; and financial applications where sector membership
affects stock returns through multiple financial characteristics
\cite{guo2023statistical}.

A growing literature has developed methods for mediation analysis where the mediators and confounders are high-dimensional. One class of approaches reduces the mediator
dimension using principal components or other low-dimensional projections
\cite{chen2015highdim,huang2016hypothesis}. A second class screens candidate
mediators using marginal exposure-mediator or mediator-outcome associations
\cite{hanson2016computational,zhang2016estimating}. Penalized methods instead
estimate a relatively small collection of active mediation pathways by
imposing sparsity on the relevant regression coefficients
\cite{zhang2016estimating,dai2022hima2}. These procedures address important goals
such as dimension reduction and mediator selection, but their performance and
interpretation depend on the assumed low-dimensional or sparse representation
of the mediation signal.

Another strategy for inference in  mediation models is motivated by debiased or
desparsified estimators of low-dimensional functionals of high-dimensional
regression parameters. 
 In high-dimensional mediation analysis,
\cite{sparse_med} applies a debiasing correction directly to a quantity
associated with the mediation effect. Subsequently, the authors in \cite{lin2025testing} employ a variance-enhanced
projection direction, adapted from \cite{cai2021optimal}, to relax an assumption on covariance structures of the mediators and
exposures used in \cite{sparse_med}, while still assuming sparsity  on $\tm$.
A related but distinct approach is
\cite{guo2023statistical}, which uses partially penalized least squares with
a general penalty. However, it assumes sparsity on $\tm$.

Sparsity-based methods are natural when the indirect effect  concentrates on a small number of mediators. In high-throughput
biological and imaging applications, however, this assumption is not always
reasonable. Polygenic models of complex traits allow
aggregate effects to be distributed across many weak molecular
contributions \cite{boyle2017omnigenic}, while empirical studies on high-dimensional mediation
have considered non-sparse epigenetic signals
\cite{clarkboucher2023methods} and 
dense, structured mediation patterns across correlated imaging measurements
\cite{zhao2023multimodalmediation}. Recent work has also emphasized that
individually weak mediation effects may be poorly estimated or omitted even
when their aggregate contribution is scientifically meaningful
\cite{li2025weakmediation}.

In this paper we  consider a proportional asymptotics regime where
$p/m\rightarrow\kappa\in(0,\infty).$
This regime is appropriate when the number of measured mediators remains a
non-negligible fraction of the available sample size, as occurs in many omics and imaging studies.
A broad collection of tools has been developed in the past decade to obtain precise asymptotic
characterizations of high-dimensional estimators in this regime. Prominent
approaches include approximate message passing
\cite{donoho2009message,bayati_mont,SurPNAS,lnk_flra_survey,AMP_survey}; leave-one-out arguments
\cite{elkaroui_pnas_1,elkaroui_pnas_2,elkaroui_ptrf,SurPNAS,chen2021spectral}; and the convex Gaussian min-max theorem
\cite{stojnic_1,stojnic_2,Tharmpoulidis_precise}. Random-matrix deterministic equivalents and resolvent methods provide
a complementary framework for analyzing regularized estimators under
proportional asymptotics
\cite{Liao_RMT,hastie2022surprises,misiakiewicz2024non,mcgrath2026optimal}. 

Our work lies
primarily within this literature, specifically using such techniques to construct debiased estimators. Debiasing
methods typically require sparsity in the regression coefficients
or suitable precision-matrix directions
\cite{van2014asymptotically,javanmard2014confidence} and therefore do
not extend directly to non-sparse settings. However, sparsity can sometimes
be replaced by alternative geometric or structural assumptions on the
design \cite{zhu2018linear,BZ_cvx_debias,li2023spectrum}, and the results in this paper belong to the same paradigm. More broadly,
our results contribute to the emerging literature on inference for
low-dimensional functionals under proportional asymptotics
\cite{SurPNAS,jiang2025new,chen2024method, Saenz_Sur_2025}.

\subsection{Challenges in the mediation model}\label{subsec:challenges}
We utilize debiasing strategies to estimate the mediation effects. Before describing our estimators, we highlight the main challenges that arise in applying these strategies to the mediation model.
Extending debiasing methods to non-sparse mediation models presents an important difficulty. Existing procedures often require knowledge of the covariance matrix of the outcome-regression design or suitable structural assumptions on target-specific precision directions; see, for example, \cite{cai_guo,javan_mont3,BZ_dof}. In our model, the outcome-regression design consists of $(M_i^\top,A_i,C_i^\top)^\top$, whose covariance depends on the unknown mediator-model coefficients, including $\ba$. Consequently, the relevant population covariance or precision information is not available a priori. Moreover, we do not impose sparsity assumptions on $\tm$ or on the relevant precision-matrix directions that could facilitate conventional debiasing procedures. Establishing valid inference under these conditions, without prior knowledge of the full design covariance, is therefore a central technical challenge.

A second challenge is that the global mediation effect, $\ba^\top\tm$, is the inner product of two high-dimensional estimands. Even if asymptotically bias-corrected estimators $\hba$ and $\htm$ are constructed separately, establishing asymptotic normality of $\hba^\top\htm$ requires controlling their dependence and the additional product fluctuations. We address this by splitting the sample into two disjoint
folds and estimating $\beta_A$ and $\theta_M$ on
independent folds, which facilitates the asymptotic analysis
of their inner product.   To improve efficiency, we then cross-fit: the roles of the two folds are reversed to construct a second sample-split estimator, and the two estimators are averaged to form the final cross-fitted estimator. This procedure uses each fold once for estimating $\beta_A$ and once for estimating $\theta_M$, thereby exploiting the full sample more effectively and reducing variance relative to a single sample split. However, the two cross-fitted estimators are potentially dependent since they are constructed from the same observations with the fold roles reversed. Accounting for this  dependence requires substantial mathematical analysis.

Finally, the limiting variance involves several unknown signal and covariance quantities, including $\lVert\ba\rVert^2$, $\lVert\tm\rVert^2$, the mediation effect $\gamma$ itself, and functionals of $\Sigma_B$, the covariance matrix of $(A_i,C_i^\top)^\top$. These quantities must be consistently estimated in the same non-sparse regime.
Addressing these challenges requires a debiasing argument tailored to our mediation framework, rather than a direct application of existing high-dimensional debiasing methods.

\subsection{\textbf{Our Contributions}}
We now present our main contributions below: 

\vspace{0.2in}

\begin{enumerate}[label=\arabic*.]

\item We construct an asymptotically bias-corrected estimator
of the global mediation effect $\gamma$ and establish its
asymptotic normality at the $\sqrt{m}$ scale under
proportional high-dimensional asymptotics, without imposing
sparsity assumptions on the mediation coefficients (Theorem \ref{thm_plugin}). 
We overcome the covariance difficulty (see Section \ref{subsec:challenges}) through an
\emph{approximate debiasing} scheme based on
\emph{partial ridge} regression, which exploits the
structure of the linear mediation model to correct
the estimation bias without requiring knowledge of
the full design covariance.
Thus, 
beyond establishing inference for the global
mediation effect in the non-sparse regime, 
our results
demonstrate how the structural information in a mediation
model can circumvent a covariance-related bottleneck in
high-dimensional debiasing.

\item We propose a cross-fitted estimator that reduces the
variance of the plug-in estimator obtained via sample splitting (Theorem \ref{thm_cf}).
 Its asymptotic variance is exactly one half
of that of the single-split plug-in estimator under the
common $\sqrt{m}$ normalization. While the spirit of cross-fitting is similar to \cite{jiang2025new} (for an AIPW estimator), our analysis and variance characterization are fundamentally different. We further show that under the \emph{small-aspect-ratio limit} the asymptotic variance of the cross-fitted estimator matches the variance of the efficient estimator in the classical case.

\item We develop consistent estimators for the limiting variances of the global mediation effect estimator $\hat{\gamma}$ and of its cross-fitted version (Theorem \ref{fluc_est_var}). Furthermore, we quantify precisely how the variability of our estimator depends on the problem parameters. 
\item 
In addition to the global mediation effect $\gamma$, we also develop a $\sqrt{m}$ consistent estimator for the pathwise mediation effect (Theorem \ref{thm_plugin_path}), and provide a consistent estimate of its limiting variance (Theorem \ref{fluc_est_var_path}).

\item Technically, our proofs rely on the theory of \emph{deterministic equivalents} \cite{Liao_RMT}, specifically  fluctuations of bilinear forms of the resolvent \cite{yin2023limiting}. The partial ridge construction leads to bilinear and quadratic forms in random matrices; analyzing them using deterministic equivalents yields an asymptotically bias-corrected estimator of $\gamma$. Although deterministic-equivalent techniques have been used before to analyze ridge-type procedures in the proportional asymptotics regime \cite{hastie2022surprises}, the dependence structure induced by the LSEM makes our analysis substantially more involved.

\item To assess the finite-sample behavior of our
proposed estimators, we present a collection of
numerical experiments in Section~\ref{sec:num_exp}.

\end{enumerate}

\subsection{Outline of the paper}

We discuss the model assumptions and state the main results in Section \ref{sec:setup} and Section \ref{sec:main_results} respectively. This includes our estimator for the global mediation effect using partial ridge regression (Theorem \ref{thm_plugin}) as well as its cross-fitted version (Theorem \ref{thm_cf}). In Theorem \ref{fluc_est_var} we show that the limiting variances of the previous two theorems can be estimated and tests for nullity of mediation effects can be designed. The corresponding results for the pathwise mediation effect are given in Theorems \ref{thm_plugin_path}, \ref{thm_cf_path}, and \ref{fluc_est_var_path} respectively. The numerical experiments are given in Section \ref{sec:num_exp}, while the proof outline is given in Section \ref{sec:proof_outline}. For the sake of brevity we only provide an abbreviated version of the limiting variances in Section \ref{sec:main_results}, and defer the detailed formulas to Appendix \ref{sec:var_details}. Finally, Appendices \ref{sec:proofs}, \ref{sec:tech_lemma}, and \ref{sec:aux_lemma} contain the proofs of the main theorems, proofs of the technical lemmas, and proofs of auxiliary lemmas respectively.

\subsection{Notations}

We use $O(.)$ and $o(.)$ to denote asymptotic boundedness and convergence to $0$,
respectively; i.e., we say $a_n=O(\gamma_n)$ if $\limsup|a_n/\gamma_n|\leq C$ and $a_n=o(\gamma_n)$ if $\limsup|a_n/\gamma_n|=0$. Similarly $O_p(.)$ and $o_p(.)$ are reserved for the stochastic versions of the same quantities. We will use $\xrightarrow{P}$ to denote convergence in probability and $\xrightarrow{d}$ to denote convergence in distribution. The $\ell_2$ norms of vectors are denoted by $\|.\|$  unless otherwise mentioned. For matrices the $\|A\|$ will always denote the largest singular value of $A$ (operator norm) and hence the same notation $\|.\|$ will be reserved for the  $\ell_2$ norm of a vector and the operator norm of the matrix, the difference being understood from the context. Throughout the notation $C$ will denote a generic constant unless otherwise specified. We will denote the $m\times n$ matrix with $0$ entries by $0_{m,n}$ and the $m\times m$ identity matrix by $I_m$. To avoid confusion with the sample matrices $C,M$ and the vectors $Y, A$, used throughout the high-dimensional analysis, we use $(\stY,\stA,\stM,\stC)$ to denote a generic observation distributed as $(Y_i,A_i,M_i,C_i)$.

\section{Preliminaries\label{sec:setup}}

\subsection{Setup and assumptions}
We observe $m$ i.i.d.~samples $\{Y_i,M_i,C_i,A_i,1\leq i\leq m\}$ with the outcome and exposure $Y_i, A_i \in \mathbb{R}$, mediators $M_i \in \mathbb{R}^p$, and confounders $C_i \in \mathbb{R}^d$.
We operate in a high-dimensional setting where the mediator dimension $p$, and the covariate dimension $d$ are allowed to grow with the sample size $m$. To make this rigorous, we consider a sequence
of problem instances with observations 
$$\{Y_{i,m},M_{i,m},A_{i,m},C_{i,m},1\leq i\leq m\}_{m\geq 1},$$  
and a sequence of errors and parameters given by
$$\{\delta_{i,m},\epsilon_{i,m},1\leq i\leq m\}_{m\geq 1}, \quad \{\theta_M(m),\theta_C(m),\theta_A(m), \beta_C(m), \beta_A (m)\}_{m\geq 1}$$
respectively. We work under the setting of a linear structural equation model
\begin{align}
Y_{i,m}&=\theta_M^T(m)M_{i,m}+\theta_A(m)A_{i,m}+\theta_C(m)^T C_{i,m}+\delta_{i,m}\\
M_{i,m}&=\beta_A(m) A_{i,m}+\beta_C(m)C_{i,m}+\epsilon_{i,m}
\end{align}
with $Y_{i,m}, A_{i,m}, \delta_{i,m}\in \mathbb{R}$, $C_{i,m} \in \mathbb{R}^{d(m)}$, and $M_{i,m}, \epsilon_{i,m}\in \mathbb{R}^{p(m)}$. 
Further for the parameters we have $\beta_A(m), \theta_M(m) \in \mathbb{R}^{p(m)}$, $\beta_C(m) \in \mathbb{R}^{p(m) \times d(m)}, \theta_A(m) \in \mathbb{R}, \theta_C (m) \in \mathbb{R}^{d(m)}$.

We consider a triangular array indexed by the sample size $m$. The
covariate vector $(A_{i,m},C_{i,m}^\top)^\top$ has dimension $d_m+1$,
where
$
    d_m+1<\frac{m}{2}.
$
We allow $p_m$, $d_m$, and $m$ to diverge jointly such that
\[
    \frac{p_m}{m}\rightarrow \mathfrak{s}_1\in(0,\infty),
    \qquad
    \frac{d_m}{m}\rightarrow \mathfrak{s}_2\in(0,1/2).
\]
All random variables and parameters may depend on $m$, but we suppress
this dependence whenever it is clear from the context. In particular,
we write
\[
    \gamma=\ba^\top\tm
    \qquad\text{and}\qquad
    \eta_j=\beta_{A,j}\theta_{M,j},
    \quad j=1,\ldots,p,
\]
for the global and path-specific mediation effects, respectively.

For $i=1,\ldots,m$, let $M_i,\epsilon_i\in\mathbb{R}^p$,
$C_i\in\mathbb{R}^d$, and $Y_i,A_i,\delta_i\in\mathbb{R}$. We stack
the vector-valued observations row-wise as
\[
    M=
    \begin{bmatrix}
        M_1^\top\\
        \vdots\\
        M_m^\top
    \end{bmatrix}
    \in\mathbb{R}^{m\times p},
    \qquad
    C=
    \begin{bmatrix}
        C_1^\top\\
        \vdots\\
        C_m^\top
    \end{bmatrix}
    \in\mathbb{R}^{m\times d},
    \qquad
    E=
    \begin{bmatrix}
        \epsilon_1^\top\\
        \vdots\\
        \epsilon_m^\top
    \end{bmatrix}
    \in\mathbb{R}^{m\times p},
\]
and collect the scalar-valued observations into
\[
    Y=(Y_1,\ldots,Y_m)^\top,\qquad
    A=(A_1,\ldots,A_m)^\top,\qquad
    D=(\delta_1,\ldots,\delta_m)^\top
    \in\mathbb{R}^m.
\]
Finally, define the combined design matrix for the mediation  and its corresponding
coefficient matrix by
\[
    B=[\,A\;\;C\,]\in\mathbb{R}^{m\times(d+1)}
    \qquad\text{and}\qquad
    \varrho=
    \begin{bmatrix}
        \beta_A,
        \beta_C
    \end{bmatrix}
    \in\mathbb{R}^{p\times (d+1)},
\]
so that the mediator model can be written compactly as
\[
    M=B\varrho^\top+E.
\]

Write $m=2n$ and define
\[
c_2:=\lim_{n\to\infty}\frac{d+1}{n}=2\mathfrak s_2\in(0,1).
\]
We work with the following assumptions on the data generating mechanism.

\begin{assumption}{\label{ass:rows}}
    The matrix $B=(A,C)$ has i.i.d. rows $B_i$ with $\E[B_i]=0$,
    and $\Sigma_B:=\E[B_1B_1^T]\in \mathbb{R}^{(d+1) \times (d+1)}$ an unknown positive definite matrix. The matrix
$\Sigma_B$ satisfies the uniform (over $d(m)$)
spectral bounds
\[
0<\underline\sigma \le \lambda_{\min}(\Sigma_B)\le \lambda_{\max}(\Sigma_B)\le \overline\sigma<\infty.
\]
Further let $B_a$ be the matrix formed by the first $m/2$ rows of $B$. Then we assume that the following identity holds:

\begin{equation}e_1^\top(\frac{1}{n}B_a^TB_a)^{-1}e_1=\frac{e_1^\top\Sigma_B^{-1}e_1}{1-c_2}+o_p(1),  \label{inv_samp_conv}\end{equation}
\begin{equation*}e_1^\top\Sigma_B^{-1}e_1\rightarrow \tau \end{equation*}
for some $\tau\in (0,\infty)$, and
\begin{equation}P(\rm{rank}(B_a)=d+1)\rightarrow 1.\label{rank_assumption}\end{equation}

\end{assumption}

Since the identity given in \eqref{inv_samp_conv} is widely used throughout the paper, for convenience we define the following notation:
\[
f^{\Sigma_B}_5(c_2):=\frac{\tau}{1-c_2}.
\]

\begin{remark}\label{rem:ass_1}
Assumption~\ref{ass:rows} ensures that the OLS estimator of
$\beta_A$ is well behaved in the proportional asymptotics regime.
Standard calculations show (see Section \ref{OLS_construction}) that the OLS estimator constructed from the second
fold satisfies
\[
\hat{\beta}_{A,j}-\beta_{A,j}\mid B_b
\sim
N\left(0,e_1^\top(B_b^\top B_b)^{-1}e_1\right).
\]
Thus coordinatewise consistency follows whenever
$e_1^\top(B_b^\top B_b)^{-1}e_1=o_p(1)$.
Assumption~\ref{ass:rows} imposes the stronger condition
\[
ne_1^\top(B_b^\top B_b)^{-1}e_1
\xrightarrow{P}
\frac{\tau}{1-c_2},
\]
which guarantees coordinatewise consistency at the
$n^{-1/2}$ rate and, more importantly, identifies the
deterministic asymptotic variance required for our inference
procedures.\end{remark}

The two conditions
\eqref{inv_samp_conv} and \eqref{rank_assumption}
hold for the following classes of designs, as established
in the corresponding appendices.

 \begin{itemize}
     \item The rows of the matrix $B$ are i.i.d. of the form $\Sigma^{1/2}_B Z_i$, where the coordinates of $Z_i$ i.e. $Z_{ij}$-s are i.i.d from a common distribution, independent of the dimension, with mean zero, variance $1$, and finite fourth moment.
      \item Let $Y\sim N(0,I_d)$ and let
$\Sigma=\Sigma_d\in\mathbb R^{d\times d}$ be
deterministic such that, uniformly in $d$,
\begin{equation}
0<\underline{\sigma}_C
\le\lambda_{\min}(\Sigma)
\le\lambda_{\max}(\Sigma)
\le\overline{\sigma}_C<\infty.
\label{C_spectral_bounds}
\end{equation}
Let $h:\mathbb R^d\to\mathbb R$ be Lipschitz with
a uniformly bounded Lipschitz constant, and let
$
\sigma(t):=(1+e^{-t})^{-1}
$
denote the logistic function. Define
\[
A_*=
\sigma\!\bigl(h(\Sigma^{1/2}Y)\bigr)
-
\E\!\left[
\sigma\!\bigl(h(\Sigma^{1/2}Y)\bigr)
\right],
\qquad
C_*=\Sigma^{1/2}Y.
\]
The rows of $B$ are independent copies of
$X=(A_*,C_*^\top)^\top$.
We show in Appendix~\ref{sec:lip_conc}
\[
0<\underline{\sigma}_B
\le
\lambda_{\min}(\Sigma_{B})
\le
\lambda_{\max}(\Sigma_{B})
\le
\overline{\sigma}_B<\infty
\]
holds uniformly in $d$ if and only if
\begin{equation}
\inf_d
\mathrm{Var}\!\left[
\sigma\!\left(h_d(\Sigma_d^{1/2}Y_d)\right)
\right]
>0.
\label{eq:logistic_non_degeneracy_statement}
\end{equation}
Under \eqref{C_spectral_bounds} and
\eqref{eq:logistic_non_degeneracy_statement},
the covariance matrices $\Sigma_B=\E[XX^\top]$
are uniformly well conditioned.
Appendix~\ref{sec:lip_conc} establishes the
concentration and lower spectral edge properties
needed to verify \eqref{inv_samp_conv}
and \eqref{rank_assumption}.
The additional population limit
$e_1^\top\Sigma_B^{-1}e_1\to\tau$
remains a separate requirement of
Assumption~\ref{ass:rows}.

 \end{itemize}

\begin{assumption}{\label{ass:errors}}
    The errors $\epsilon_i$ in the mediation regression model are i.i.d. Gaussian random vectors with distribution $N(0,I_p)$; equivalently, $E$ is an $m \times p$ matrix with i.i.d. standard
Gaussian entries. The error in the outcome regression model \eqref{y_mod_vec} is distributed as centered Gaussian  $D\sim N(0, \sigma_D^2I_m)$, for a known variance  $\sigma_D^2$. Both $E$ and $D$ are independent of each other and of $B$.
\end{assumption}

\begin{remark}\label{whitening_remark}
 For the global mediation effect, this assumption extends to
$E=GV^{1/2}$, where $G$ is an $m\times p$ matrix with i.i.d.\
standard Gaussian entries and $V$ is a known, symmetric positive
definite matrix. Indeed for known $V$ adopting the following whitening transformation: 
$$\tilde{M}=MV^{-1/2},\quad \tilde{\varrho}=V^{-1/2}\varrho, \quad \tilde{E}=EV^{-1/2},\quad \tilde{\theta}_M=V^{1/2}\theta_M$$ we have $$\tilde{M}=B\tilde{\varrho}^T+\tilde{E}, \quad Y=A\theta_A+C\theta_C+\tilde{M}\tilde{\theta}_M+D.$$ Since $\tilde{\beta}_A^T\tilde{\theta}_M=\ba^T\tm$ the results remain unchanged. Henceforth, we work with the whitened model ($V=I$) and omit tildes from all transformed parameters. Also, for simplicity, we will assume $\sigma^2_D=1$.
\end{remark}

\begin{assumption}\label{ass:signal}
We assume that the following limits exist:
\[
\|\ba\|^2\rightarrow\mathfrak{b}_A,
\qquad
\|\tm\|^2\rightarrow\mathfrak{t}_M,
\qquad
\ba^\top\tm\rightarrow\mathfrak{g},
\]
where $\mathfrak{b}_A,\mathfrak{t}_M\in(0,\infty)$
and $\mathfrak{g}\in\mathbb{R}$.
\end{assumption}

The quantities $\mathfrak{b}_A$, $\mathfrak{t}_M$, and $\mathfrak{g}$ admit a natural interpretation in
terms of the variability of the corresponding
contributions to the structural equations.
First, since the mediator errors have identity covariance,
\[
\operatorname{Var}(\tm^\top M_i\mid A_i,C_i)
=
\|\tm\|^2.
\]
Thus $\mathfrak{t}_M$ represents the limiting strength of the mediator contribution to the outcome.
Similarly, the contribution of the exposure to the
mediator vector is $\ba A_i$. Conditioning on $C_i$,
we obtain
\[
\tr\left\{
\operatorname{Cov}(\ba A_i\mid C_i)
\right\}
=
\operatorname{Var}(A_i\mid C_i)\|\ba\|^2.
\]
Hence $\mathfrak{b}_A$ controls the magnitude of the
exposure-related signal in the mediator model,
relative to the conditional variability of the exposure.
Finally, $\mathfrak{g}$ is the limiting global mediation
effect $\ba^\top\tm$. Thus Assumption~\ref{ass:signal}
controls the aggregate strengths of the two structural
relationships \eqref{y_mod_vec} and \eqref{M_mod_vec} and their resulting mediation effect.

The next assumption is related to the tuning parameter of the partial ridge estimator which will be defined in \eqref{btm_eqn_1}.  Instead of choosing $\lambda_1$ exactly proportional to $m$ we choose an equivalent scaling, which eases our analysis.

\begin{assumption}\label{ass:tune}
Fix a constant $\kappa_1\in [a,b]$, with $0<a<b<\infty$. Then the tuning parameter $\lambda_1$ is chosen as follows:
$\lambda_1:=\kappa_1(m/2-d-1).$
\end{assumption}

While Assumptions \ref{ass:rows}-\ref{ass:tune} need to hold for all the results, the next assumption is needed only for the estimation of pathwise mediation effect $\ej$.
\begin{assumption}{\label{ass:path}}
The parameters
$\beta_{A,j}$ and $\theta_{M,j}$ remain constant as $n$ varies. Further we assume $\beta^2_{A,j}+\theta^2_{M,j}>0$.
\end{assumption}

Assumption~\ref{ass:path} is imposed only for inference on the
pathwise mediation effect $\eta_j$. Since
$\eta_j=\beta_{A,j}\theta_{M,j}$, its first-order
fluctuations depend on the gradient
\[
\nabla(\beta_{A,j}\theta_{M,j})
=
(\theta_{M,j},\beta_{A,j}).
\]
The condition
$\beta_{A,j}^2+\theta_{M,j}^2>0$
ensures that this gradient is nonzero and that the
estimator admits a nondegenerate Gaussian limit at
the $\sqrt{m}$ scale. In particular, the assumption
allows the regular null cases in which exactly one
of $\beta_{A,j}$ and $\theta_{M,j}$ vanishes, and
does not require $\eta_j\neq 0$. When both
coefficients are zero, the first-order limit is
degenerate, and a non-Gaussian second-order limit
arises at the $m$ scale; see
Appendix~\ref{sec:pathwise_nonregular}.
No analogous coordinatewise nondegeneracy condition
is needed for inference on the global mediation effect.

Assumptions~\ref{ass:signal} and~\ref{ass:path}
are compatible. For example, the $j$th coordinates
of $\beta_A$ and $\theta_M$ may remain fixed while
the remaining coordinates are of order
$O(p^{-1/2})$, provided their aggregate squared
norms and inner product admit the limits specified
in Assumption~\ref{ass:signal}. The requirement
that $\beta_{A,j}$ and $\theta_{M,j}$ remain
constant can also be relaxed to convergence
towards finite limits $b_j$ and $t_j$ satisfying
$b_j^2+t_j^2>0$, with the limiting variance
evaluated at these limits. A nonzero limiting
pathwise effect additionally requires $b_jt_j\neq0$.

We are now ready to construct the estimators for the mediation effect.

\section{Main Results}\label{sec:main_results}
This section presents our main results. We begin by describing our inferential procedure for the global mediation effect and its associated central limit theorems.

\subsection{Inference for the global mediation effect}

To construct an estimator for the global mediation effect, we split the data into two parts $\mathfrak{D}_a$ and $\mathfrak{D}_b$ with each subsample being of size $n:=m/2$. We will denote the quantities $Y,M,C,A,D,E$ obtained from the parts $\mathfrak{D}_a$ and $\mathfrak{D}_b$ by $Y_a,M_a,C_a,A_a,D_a,E_a$ and \\$Y_b,M_b,C_b,A_b,D_b,E_b$ respectively.

All main mediation central limit theorems are stated using the
normalization $\sqrt{m}$, where $m=2n$ is the total sample size.
The single-split plug-in estimator uses both folds, estimating
$\theta_M$ on one fold and $\beta_A$ on the other. Cross-fitting
averages the two estimators obtained by reversing these roles. We retain $n$ for the
fold dimensions and the corresponding matrix calculations.

From the first part $\mathfrak{D}_a$ we estimate $\theta_M$ using a partial ridge estimator and from the second part $\mathfrak{D}_b$ we estimate $\beta_A$ using an OLS estimator. We combine the two estimators to obtain  our estimator for the global mediation effect $\gamma$. The procedure is outlined in Algorithm \ref{alg:mediation}, in which we use the following convention.
Write $B_s=[A_s\;C_s]$ for $s\in\{a,b\}$. All formulas involving ordinary inverses are understood on the event
\[
\mathcal R_n=
\{\operatorname{rank}(B_a)=\operatorname{rank}(B_b)=d+1\}.
\]
On its complement, define the corresponding estimators to be zero.
Since the two folds have the same distribution,
Assumption~\ref{ass:rows} implies $\mathbb P(\mathcal R_n)\to1$;
this convention does not affect the stated asymptotic results.
Exact identities involving these inverses are understood on $\mathcal R_n$.

\begin{algorithm}
\caption{Plugin estimator for the global mediation effect}
\label{alg:mediation}
\begin{algorithmic}[1]
\vspace{0.1in}
\Statex Define $\mathcal{E}_2 := (0_{p,d+1}, I_p)$ and $I_{22}:=\left(\begin{array}{cc}0_{d+1,d+1} & 0_{d+1,p}\\0_{p,d+1} & I_p\end{array}\right)$
\vspace{0.1in}
\Statex \textbf{Partial ridge estimation}
\vspace{0.1in}
\State $X_a \gets (A_a, C_a, M_a)$
\State $\btm \gets \mathcal{E}_2 (X_a^\top X_a + \lambda_1 I_{22})^{-1} X_a^\top Y_a$
\vspace{0.1in}
\Statex \textbf{Bias correction}
\vspace{0.1in}
\State $\htm \gets \cin \breve{\theta}_M$ (for a deterministic constant $\cin$, see \eqref{htm_defn} for details)
\vspace{0.1in}
\Statex \textbf{OLS estimation}
\vspace{0.1in}
\State $B_b \gets (A_b, C_b)$
\State $\hba^T \gets e_1^\top (B_b^\top B_b)^{-1} B_b^\top M_b$
\vspace{0.1in}
\Statex \textbf{Mediation effect}
\vspace{0.1in}
\State $\hat{\gamma} \gets \hba^\top \htm$
\vspace{0.1in}
\State \Return $\hat{\gamma}$
\end{algorithmic}
\end{algorithm}

Next we discuss the estimators $\htm$ and $\hba$ in detail.

 \subsubsection{Estimator of $\tm$}
 
 We consider a partial ridge penalty for the estimation of the parameter $\tm$. Let $B_a=(A_a,C_a)$ and $X_a=(B_a,M_a)$. Denote the partial ridge estimator by
\begin{align}\btm:&=\mathcal{E}_2(X_a^TX_a+\lambda_1 I_{22})^{-1}X_a^TY_a= (M_a^TP^{\perp}_{B_a}M_a+\lambda_1 I)^{-1}M_a^TP^{\perp}_{B_a}Y_a \label{btm_eqn_1}\\
\text{ where }\mathcal{E}_2&=(0,I_p), \quad P^{\perp}_{B_a}=I- B_a(B_a^TB_a)^{-1}B_a^T.\nonumber\end{align}
  Recall that the rows of $E$ have a  known covariance structure, with $V$ as the corresponding covariance matrix. In the following discussion, we assume without loss of generality that $V=I$ (see Assumption \ref{ass:errors}). Writing $\bea=P^{\perp}_{B_a}E_a$, we obtain the following simplification.

 \begin{lemma}{\label{partial_ridge_lemma}}
  The partial ridge estimator can be rewritten as
\begin{equation}\breve{\theta}_M=\tm-\lambda_1H_{\lambda_1,a}^{-1}\tm+ H_{\lambda_1,a}^{-1}E_a^TP^{\perp}_{B_a}D_a\label{ridge_decomp_simplified},\end{equation}
where $H_{\lambda_1,a}=\bea^T\bea+\lambda_1I$.

 \end{lemma}

It follows from \eqref{ridge_decomp_simplified} that the partial ridge
estimator $\btm$ has a non-negligible shrinkage bias.
Recall the scaled tuning parameter $\kappa_1=\lambda_1/(n-d-1)$, the fold size, $n=m/2$, and the aspect ratio limit $c_2$ defined in the
preliminaries. For the residualized mediator regression we
also record the equivalent aspect-ratio convention
\begin{align}
\frac{p}{n-d-1}
&\rightarrow c_1\in(0,\infty),
\label{p_conv}\\
\frac{d+1}{n}
&\rightarrow c_2\in(0,1).
\label{d_conv}
\end{align}
Next define the scaled resolvent as \[H_{\kappa_1,a}^{-1}:=(\bea^T\bea/(n-d-1)+\kappa_1I)^{-1}\]
and note that it can be approximated as follows
\[(\bea^T\bea/(n-d-1)+\kappa_1I)^{-1}\approx \fin I\]
where the constant $\fin$ and the precise meaning of the approximation are given in Section \ref{asymp_debias}.
The bias-correction factor is then given by
\begin{equation}
\cin
:=
\left\{1-\kappa_1\fin\right\}^{-1}, \quad
\htm=\cin\btm.
\label{htm_defn}
\end{equation}

We will show that the correction obtained by premultiplying $\btm$ by the constant $\cin$ is crucial in constructing our asymptotically bias-corrected estimator for $\gamma$. While the detailed argument is deferred to the appendix we provide the intuition behind this correction in Section \ref{asymp_debias} using the theory of deterministic equivalent. 
 
Finally, we note that, although the definition of $\fin$ is involved, its limit as $n \rightarrow \infty$ is given by the following function:
\begin{equation}
f_1(\kappa_1,c_1)
:=\frac{-(1-c_1+\kappa_1)+\sqrt{(1-c_1+\kappa_1)^2+4c_1\kappa_1}}{2c_1\kappa_1}. \label{fixed_pt_1}
\end{equation}
 
The function $f_1(\kappa_1,c_1)$ appears in the fixed point equation of the Mar\v{c}enko-Pastur theorem (cf. Theorem \ref{mp_detequ}) and will feature in the construction of our estimators for the limiting variances in Theorems \ref{fluc_est_var} and \ref{fluc_est_var_path}.

\subsubsection{Estimator for $\ba$}\label{OLS_construction}

Next we will consider an OLS estimator for the parameter $\ba$, estimated from the part $\mathfrak{D}_b$. Recall that $\varrho=[\,\beta_A\;,\;\beta_C\,]$ and define
 the OLS estimator for $\varrho$ as
\begin{equation*}\hat{\varrho}^T:=(B_b^TB_b)^{-1}B_b^TM_b=\varrho^T+(B_b^TB_b)^{-1}B_b^TE_b.
\end{equation*}

Since $\ba^T=e_1^T\varrho^T$ we define the OLS estimator for $\ba$ as follows:
\begin{equation}\hba^T:=e_1^T\hat{\varrho}^T=e_1^T(\varrho^T+(B_b^TB_b)^{-1}B_b^TE_b)=\ba^T+e^T_1(B_b^TB_b)^{-1}B_b^TE_b.\label{hba_defn}
\end{equation}

\begin{remark}
The OLS expression for $\hba$ is well defined whenever $B_b$ has
full column rank, an event whose probability tends to one under
Assumption~\ref{ass:rows}. The restriction $d+1<n$ ensures positive
residual degrees of freedom for the partial-ridge construction.
Importantly, the OLS estimator of $\ba$ does not require knowledge
of the population covariance matrix $\Sigma_B$ of the exposure
and covariates.

\end{remark}

Now we are ready to construct an estimator for the global mediation effect.

\subsubsection{\textbf{Plugin estimator for $\gamma$ and its behavior}}

Consider the estimators $\hba$ and $\htm$ of $\ba$ and $\tm$. First we construct the plug-in estimator $\hat{\gamma}:=\hba^T\htm$ to estimate $\gamma$. To establish a CLT for the estimator $\hat{\gamma}$, we need to characterize the limiting variance of the same estimator.  
 Recall the definitions of $c_1$ and  $c_2$:
$$c_1=\lim_{n\rightarrow \infty}\frac{p}{n-d-1},\quad c_2=\lim_{n\rightarrow \infty}\frac{d+1}{n}$$ from \eqref{p_conv} and \eqref{d_conv}. Also recall the definition  $\kappa_1:=\frac{\lambda_1}{n-d-1}$ from Assumption \ref{ass:tune} 
and the limit of the signal strengths and overlaps
$$\lim_{n\rightarrow \infty}\|\ba\|^2=\mathfrak{b}_A, \lim_{n\rightarrow \infty}\|\tm\|^2= \mathfrak{t}_M, \lim_{n\rightarrow \infty}\ba^T\tm=\mathfrak{g}$$
from Assumption \ref{ass:signal}.

Let  
\begin{align}
&\sigma^2(\kappa_1,c_1,c_2,f^{\Sigma_B}_5(c_2),\mathfrak{b}_A,\mathfrak{t}_M,\mathfrak{g})\nonumber\\
=&\varphi_1(\kappa_1,c_1,c_2)f^{\Sigma_B}_5(c_2)+\varphi_2(\kappa_1,c_1,c_2)\mathfrak{b}_A+\varphi_3(\kappa_1,c_1,c_2)f^{\Sigma_B}_5(c_2)\mathfrak{t}_M\nonumber\\
&+\varphi_4(\kappa_1,c_1,c_2)(\mathfrak{b}_A\mathfrak{t}_M+\mathfrak{g}^2),\label{lim_var_defn}
\end{align}
where $\varphi_j(.)$ for $j\in\{1,\ldots,4\}$ are explicit functions of the tuple $(\kappa_1,c_1,c_2)$. The exact descriptions of the functions are involved and we defer them to the appendix (see Appendix \ref{sec:var_details}). For the sake of brevity we will denote $\sigma^2(\kappa_1,c_1,c_2,f^{\Sigma_B}_5(c_2),\mathfrak{b}_A,\mathfrak{t}_M,\mathfrak{g})$ by $\sigma^2$ henceforth. We also show in Appendix \ref{sec:var_details} that $\sigma^2>0$. The quantity $\sigma^2$ is a common variance factor: under the
$\sqrt{m}$ normalization, the limiting variances of the plug-in
and cross-fitted estimators are $2\sigma^2$ and $\sigma^2$,
respectively. The next theorem gives the plug-in limit.
\begin{theorem}{\label{thm_plugin}}
    Consider the estimator $\hat{\gamma}:=\hba^T\htm$ for the parameter $\gamma=\ba^T\tm$ where the estimators $\hba$ and $\htm$ are given by \eqref{hba_defn} and \eqref{htm_defn} respectively. Then under Assumptions \ref{ass:rows}-\ref{ass:tune} we have the following convergence.
 $$\sqrt{m}(\hat{\gamma}-\gamma)\xrightarrow{d}N(0,2\sigma^2),$$ where $\sigma^2$ is given by \eqref{lim_var_defn}.
\end{theorem}
 Despite the limiting variance taking a rather complicated form, we provide a consistent estimate for the variance in Theorem \ref{fluc_est_var}. In practice, cross-fitting is a strategy used widely in the causal inference literature. We next show that cross-fitting reduces the variability of our estimator by a factor of two.

 \subsubsection{\textbf{Cross-fitted estimator for $\gamma$ and its behavior}}{\label{crossfit_defn}}
Next, we study the cross-fitted estimator. Cross-fitting
has been widely used for inference on low-dimensional
parameters in the presence of estimated nuisance functions;
see, for example, \cite{ChernozhukovEtAl2018DML,NeweyRobins2018}.
More recently, under proportional asymptotics,
\cite{jiang2025new} established a central limit theorem
for a cross-fitted augmented inverse probability weighted
(AIPW) estimator of the average treatment effect.
  To construct the cross-fitted estimator, we consider two permutations: the identity permutation will be denoted as $\vartheta_0=(a,b)$ and the transposition as $\vartheta_1=(b,a)$. Recall that $\hat{\gamma}=\hba^T\htm$ where $\hba$ was estimated from the part $\mathfrak{D}_b$ of the sample and $\htm$   from the part $\mathfrak{D}_a$; we denote this estimator by $\hat{\gamma}_{\vartheta_0}$. Similarly, we construct an estimator for $\gamma$ with $\hba$ estimated from the part $\mathfrak{D}_a$  and $\htm$ from $\mathfrak{D}_b$; we call this estimator $\hat{\gamma}_{\vartheta_1}$. The cross-fit estimator is defined as $\hat{\gamma}_{cf}=\frac{1}{2}\sum_{\vartheta}\hat{\gamma}_{\vartheta}$ where the sum is over the permutations $\{\vartheta_0,\vartheta_1\}$.

Analogous to the plug-in estimator, the cross-fitted estimator
$\hat{\gamma}_{cf}$ is $\sqrt{m}$-consistent and asymptotically
Gaussian.

\begin{theorem}{\label{thm_cf}}
    Consider the cross-fitted estimator $\hat{\gamma}_{cf}$.  Then under Assumptions \ref{ass:rows}-\ref{ass:tune} we have the following convergence:
    $$\sqrt{m}(\hat{\gamma}_{cf}-\gamma)\xrightarrow{d}N(0,\sigma^2)$$ where $\sigma^2$ is given by \eqref{lim_var_defn}. 
\end{theorem}

Under the common $\sqrt{m}$ normalization, cross-fitting
reduces the limiting variance from $2\sigma^2$ to
$\sigma^2$ relative to the single-split plug-in estimator.
Our result is in the same spirit as \cite{jiang2025new};
however, both the estimand and the analysis of the
corresponding estimators are fundamentally different,
leading to distinct technical challenges.

Before constructing valid estimators for the limiting variance $\sigma^2$, it is useful to examine the behavior of the limiting variance. In particular, we compute a sequential boundary limit of the
explicit asymptotic variance function. Throughout this computation, the
limiting signal parameters
$\mathfrak b_A$, $\mathfrak t_M$, $\mathfrak g$ and the population
precision parameter $\tau$ are held fixed, while we successively let
$c_2\to0$, $c_1\to0$, and $\kappa_1\to0$. 
 In Lemma \ref{low_dim_var_lemma}, we show that the following sequential limit holds:
\begin{equation}\label{eq:seq_limit}
\lim_{\kappa_1\to 0}\ \lim_{c_1\to 0}\ \lim_{c_2\to 0}\ \sigma^2
=
\mathfrak b_A+F_5(0)\,\mathfrak t_M,
\end{equation}
where
\[
F_5(0)=\lim_{n\to \infty}e_1^\top\Sigma_B^{-1}e_1.
\]
We refer to this sequential boundary limit as the
\emph{small-aspect-ratio limit}. We note that in the  low-dimensional regime, where the dimension of the mediator and the confounders remain fixed, as the sample size $m$ grows, we have $\mathfrak{b}_A=\|\ba\|^2$, $\mathfrak{t}_M=\|\tm\|^2$, and $F_5(0)=e_1^\top\Sigma_B^{-1}e_1$. And consequently  we have
\[
\mathfrak b_A+F_5(0)\,\mathfrak t_M=\|\ba\|^2 + (\Sigma_B^{-1})_{11}\|\tm\|^2.
\]

To further interpret this expression, consider the Gaussian linear structural equation model in \eqref{y_mod_vec}-\eqref{M_mod_vec}, with parameter
\[
\eta = (\theta_M,\theta_A,\theta_C,\beta_A,\beta_C),
\]
and mediation functional
\[
\gamma = \gamma(\eta):=\beta_A^\top \theta_M.
\]
As shown in Lemma \ref{fisher_info}, the efficient asymptotic variance for estimating $\gamma$ from a single observation is
\[
V_\gamma
=
\|\beta_A\|^2
+
\frac{1}{\sigma_{A\mid C}^2}\,\|\theta_M\|^2,
\]
where
\[
\sigma_{A\mid C}^2
=
\E[\stA^2]-\E[\stA\stC^\top]\E[\stC\stC^\top]^{-1}\E[\stC\stA].
\]
Noting that $(\sigma_{A\mid C}^2)^{-1}=(\Sigma^{-1}_{B})_{11}$, the efficient asymptotic variance based on $m$ samples is
\begin{equation}\label{eff_var}
V_{\gamma,m}
=
\frac{1}{m}\left(\|\ba\|^2+(\Sigma_B^{-1})_{11}\|\tm\|^2\right).
\end{equation}

Theorem~\ref{thm_cf} gives
\[
\sqrt{m}(\hat{\gamma}_{cf}-\gamma)\xrightarrow{d}N(0,\sigma^2).
\]
Combining \eqref{eq:seq_limit} and \eqref{eff_var},
we conclude that the sequential small-aspect-ratio
limit of the $\sqrt{m}$-normalized asymptotic variance
of the cross-fitted estimator coincides with the
classical efficient variance factor $mV_{\gamma,m}$.

\subsubsection{Estimation of the limiting variance}

We note that the unknown quantities appearing in $\sigma^2$ are the limiting values of $\|\ba\|^2,\|\tm\|^2, \gamma^2$ i.e. limiting values of the signal strengths along with  $f^{\Sigma_B}_5(c_2)$. Estimating these quantities and plugging them in the variance formula we obtain $\hat{\sigma}^2$.

Recall that the limit of $\fin$ is given by $f_1(\kappa_1,c_1)$  whose closed form is provided in \eqref{fixed_pt_1}.
Define two auxiliary functions which will be used in constructing the estimator $\hat{\sigma}^2$:
\begin{align*}
     f_2(\kappa,c)
    & := f_1(\kappa,c)+\kappa\,\partial_{\kappa} f_1(\kappa,c),\nonumber\\
    f_3(\kappa_1,c_1)&:=(1-2\kappa_1f_1(\kappa_1,c_1)-\kappa^2_1\partial_{\kappa}f_1(\kappa_1,c_1))^{-1}.
\end{align*}
We have the following theorem. 

\begin{theorem}{\label{fluc_est_var}}

Consider the following estimators of $\mathfrak{t}_M,f^{\Sigma_B}_5(c_2), \mathfrak{b}_A$:

\begin{align*}
\hat{\mathfrak{t}}_M&=f_3(\kappa_1,c_1)(\|\btm\|^2-c_1f_2(\kappa_1,c_1))\\
\hat{\mathfrak{r}}&=e^T_1(B_b^TB_b/n)^{-1}e_1,\qquad
\hat{\mathfrak{b}}_A=\|\hba\|^2-c_1(1-c_2)\hat{\mathfrak{r}}
\end{align*}
where $\hba$ and $\btm$ are defined in \eqref{hba_defn} and \eqref{btm_eqn_1} respectively.  Then under Assumptions \ref{ass:rows}-\ref{ass:tune} and for
\begin{align}
\hat{\sigma}^2:=&
\varphi_1(\kappa_1,c_1,c_2)\hat{\mathfrak{r}}+\varphi_2(\kappa_1,c_1,c_2)\hat{\mathfrak{b}}_A+\varphi_3(\kappa_1,c_1,c_2)\hat{\mathfrak{r}}\hat{\mathfrak{t}}_M\nonumber\\
&+\varphi_4(\kappa_1,c_1,c_2)(\hat{\mathfrak{b}}_A\hat{\mathfrak{t}}_M+\hat{\gamma}^2),\label{lim_var_est}
\end{align}
the following convergence holds: 
    \begin{equation}\hat{\sigma}^2\xrightarrow{P}\sigma^2.\label{sigma_consis}\end{equation}
\end{theorem}

Consequently, $2\hat\sigma^2$ and $\hat\sigma^2$
consistently estimate the $\sqrt m$-normalized limiting
variances in Theorems~\ref{thm_plugin} and~\ref{thm_cf},
respectively. In implementation, the limiting aspect ratios
$c_1$ and $c_2$ in the variance expressions are replaced by
\[
c_{1,n}:=\frac{p}{n-d-1},\qquad
c_{2,n}:=\frac{d+1}{n}.
\]
The same convention is used for the pathwise variance
estimator in Theorem~\ref{fluc_est_var_path}. By continuity
of the coefficient functions at the stipulated interior
limits, these substitutions preserve consistency.

While $\hat{\mathfrak{r}}$ provides a natural estimator of
$f^{\Sigma_B}_5(c_2)$ (see Assumption~\ref{ass:rows}),
the constructions of $\hat{\mathfrak{t}}_M$ and
$\hat{\mathfrak{b}}_A$ are less immediate.
From \eqref{ridge_decomp_simplified}, we observe that
$\btm$ is a biased estimator of $\tm$, and its squared
norm $\|\btm\|^2$ is likewise a biased estimator of
$\|\tm\|^2$, requiring an appropriate correction.
On the other hand, although the OLS estimator $\hba$
is conditionally unbiased given $B_b$ on the
full-column-rank event, its squared norm contains a
nonvanishing contribution from estimation error in the
high-dimensional regime.

We correct these biases by directly analyzing
$\|\btm\|^2$ and $\|\hba\|^2$, yielding the estimators
$\hat{\mathfrak{t}}_M$ and $\hat{\mathfrak{b}}_A$,
respectively.

\subsection{Inference for the pathwise mediation effect through a potential mediator $M^j$}\label{subsection:path_specific}
Next we define the pathwise mediation effect i.e. the effect of $\stA$ on $\stY$ through the mediator $\stM^j$.
From \eqref{y_mod_vec} and \eqref{M_mod_vec}, we see $\stA$ affects $\stM^j$ through the coefficient $\beta_{A,j}$ while $\stM^{j}$ affects $\stY$ through the coefficient $\theta_{M,j}$, so the pathwise mediation effect is $\beta_{A,j}\theta_{M,j}$.
Following \cite{pearl2003causality,wright1921correlation,chakrabortty2018inference}, we use the following definition of the pathwise mediation effect via Pearl's do-operation \cite{pearl2003causality}, which simplifies under our LSEM.

\begin{definition}[Causal definition of pathwise mediation effect]{\label{pathwise_defn}}
The individual mediation effect, $\eta_j$, with respect to a potential mediator
$\stM^j$ $(j=1,\ldots,p)$ is defined as
\[
\eta_j \;:=\;
\frac{\partial}{\partial a}\,
\mathbb{E}\!\left[\stY \,\mid \mathrm{do}(\stA=a)\right]
\;-\;
\frac{\partial}{\partial a}\,
\mathbb{E}\!\left[\stY \,\mid \mathrm{do}(\stM^j=m_j,\; \stA=a)\right] 
\]
\end{definition}
where the $\mathrm{do}$ operation fixes the corresponding variables to constants under the formalism of \cite{pearl2003causality}. Indeed, path-wise mediation effects are difficult to identify from observational data \citep{avin2005identifiability}. However, an LSEM provides enough information to help identify pathwise causal mediation effects. Specifically, we state the following straightforward lemma, which shows that the above definition of the pathwise mediation effect coincides with Definition \ref{def:path_med} for the LSEM model we consider.
\begin{lemma}{\label{pathwise_med_lemma}}
    Consider the mediation model given by the equations \eqref{y_mod_vec} and \eqref{M_mod_vec}. Then we have
    $\ej=\beta_{A,j}\theta_{M,j}$
    where $\ej$ is given by Definition \ref{pathwise_defn}.  
\end{lemma}

\subsubsection{Plugin estimator for $\ej$}

We will estimate $\ej$ with 
$\hat{\eta}_j=\hat{\beta}_{A,j}\hat{\theta}_{M,j}$ where we define these in the following algorithm:

\begin{algorithm}[H]
\caption{Plugin estimator for the pathwise mediation effect}
\label{alg:mediation_path}
\begin{algorithmic}[1]
\Statex \textbf{Asymptotically debiased partial ridge estimation}
\vspace{0.1in}
\State $\htm \gets \cin\mathcal{E}_2 (X_a^\top X_a + \lambda_1 I_{22})^{-1} X_a^\top Y_a$
\State $\hat{\theta}_{M,j} \gets e^T_j \htm$
\vspace{0.1in}
\Statex \textbf{OLS estimation}
\vspace{0.1in}
\State $\hba^T \gets e_1^\top (B_b^\top B_b)^{-1} B_b^\top M_b$
\State $\hat{\beta}_{A,j} \gets \hba^Te_j$
\vspace{0.1in}
\Statex \textbf{Pathwise Mediation effect}
\vspace{0.1in}
\State $\hej \gets \hat{\beta}_{A,j}\hat{\theta}_{M,j}$
\vspace{0.1in}
\State \Return $\hej$
\end{algorithmic}
\end{algorithm}

 In this case also, the limiting variance can be precisely characterized (as a function of the signal strength) as follows:
\begin{align}
&\sigma_j^2(\kappa_1,c_1,c_2,f^{\Sigma_B}_5(c_2),\beta_{A,j},\theta_{M,j},\mathfrak{t}_M)\nonumber\\
=&f^{\Sigma_B}_5(c_2)\theta^2_{M,j}+\psi_1(\kappa_1,c_1,c_2)\beta^2_{A,j}+\psi_2(\kappa_1,c_1,c_2)\beta^2_{A,j}(\theta^2_{M,j}+\mathfrak{t}_M)\label{lim_var_sigma_j}
\end{align}
where $\psi_k(.)$ for $k\in\{1,2\}$ are explicit functions of the tuple $(\kappa_1,c_1,c_2)$. Once again, the exact descriptions of the functions $\psi_k(.)$ are involved and we defer them to Appendix \ref{sec:var_details}.
For the sake of brevity we will denote $\sigma_j^2(\kappa_1,c_1,c_2,f^{\Sigma_B}_5(c_2),\beta_{A,j},\theta_{M,j},\mathfrak{t}_M)$ by $\sigma_j^2$ henceforth. Under the $\sqrt{m}$ normalization, the limiting variances of the
plug-in and cross-fitted pathwise estimators are $2\sigma_j^2$
and $\sigma_j^2$, respectively. We have the following theorem.

\begin{theorem}{\label{thm_plugin_path}}
    Consider the estimator $\hej:=\hat{\beta}_{A,j}\hat{\theta}_{M,j}$.  Then under Assumptions \ref{ass:rows}-\ref{ass:path} we have
 $$\sqrt{m}(\hej-\ej)\xrightarrow{d}N(0,2\sigma_{j}^2).$$
\end{theorem}

Compared to the limiting variance $\sigma^2$ for the global mediation effect, we observe an asymmetry in the expression of $\sigma_j^2$ with respect to the signal vectors $\ba$ and $\tm$. In particular, while $\sigma^2$ in \eqref{lim_var_defn} depends on the overall signal strengths $\|\ba\|^2$ and $\|\tm\|^2$, the coordinate-wise variance $\sigma_j^2$ from \eqref{lim_var_sigma_j} exhibits a different structure. Specifically, $\sigma_j^2$ depends on both the global signal strength $\|\tm\|^2$ and the individual coefficient $\theta_{M,j}$ for $\tm$, whereas it depends only on the coordinate $\beta_{A,j}$ and not on the full norm $\|\beta_A\|$.

This asymmetry stems from the distinct estimation procedures for $\tm$ and $\ba$. The parameter $\tm$ appears in the high-dimensional outcome regression model \eqref{y_mod_vec} and is estimated via a partial ridge estimator, while $\ba$ is estimated from the mediator model \eqref{M_mod_vec}, via underparametrized OLS regression.

More concretely, from \eqref{ridge_decomp_simplified}, we have
\[
\hat{\theta}_{M,j}
=
\cin\theta_{M,j}
-
\cin e_j^\top \lambda_1 H_{\lambda_1,a}^{-1}\tm
+
\cin e_j^\top H_{\lambda_1,a}^{-1} E_a^\top P^{\perp}_{B_a} D_a.
\]
The bias term
$
\cin \lambda_1 e_j^\top H_{\lambda_1,a}^{-1}\tm
$
depends on the entire vector $\tm$. Consequently, the total signal strength $\|\tm\|^2$ enters the limiting variance. 
In contrast, from \eqref{hba_defn}, we have
\[
\hat{\beta}_{A,j}
=
\beta_{A,j}
+
e_1^\top (B_b^\top B_b)^{-1} B_b^\top E_b e_j.
\]
Thus the dependence of $\sigma_j^2$ on $\beta_A$ is through
$\beta_{A,j}$, rather than the full norm $\|\beta_A\|$.
This reflects the unpenalized, underparameterized mediator
regression, whose covariate dimension nevertheless grows with
the sample size.

\subsubsection{Cross-fitted estimator for $\ej$}\label{sec:cross_path}

Recall that the plug-in estimator $\hej$ uses
$\hat\beta_{A,j}$ estimated from $\mathfrak D_b$ and
$\hat\theta_{M,j}$ estimated from $\mathfrak D_a$. Denote the identity as $\vartheta_0=(a,b)$ and the transposition as $\vartheta_1=(b,a)$.  Then, similar to the global mediation effect, the cross-fitted estimator for the pathwise mediation effect is defined as follows: $\hat{\eta}_{j,cf}=\frac{1}{2}\sum_{\vartheta}\hat{\eta}_{j,\vartheta}$ where the sum is over the permutations $\{\vartheta_0,\vartheta_1\}.$ The following theorem describes the limiting behaviour of the cross-fitted estimator for the pathwise mediation effect.

\begin{theorem}{\label{thm_cf_path}}
    Consider the cross-fitted estimator $\hat{\eta}_{j,cf}$.  Then under Assumptions \ref{ass:rows}-\ref{ass:path} we have the following:
    \[
\sqrt{m}(\hat{\eta}_{j,cf}-\ej)
\xrightarrow{d}N(0,\sigma_j^2),
\]
where $\sigma_j^2$ is given by \eqref{lim_var_sigma_j}.
\end{theorem}

As in the case of the global mediation effect, we observe that under the common $\sqrt{m}$ normalization, cross-fitting reduces
the pathwise limiting variance from $2\sigma_j^2$ to $\sigma_j^2$. 

\subsubsection{Estimation of the limiting variance}

Like the global mediation case, the limiting variance $\sigma_j^2$ can be estimated as well. 
From \eqref{lim_var_sigma_j} we note that the unknown quantities appearing in $\sigma_j^2$ are the limiting values of $\beta_{A,j},\theta_{M,j}, \mathfrak{t}_M$, and   $f^{\Sigma_B}_5(c_2)$. We construct the estimators of  these quantities from the estimators $\hba$ and $\htm$  of $\ba$ and $\tm$  defined in \eqref{hba_defn} and \eqref{htm_defn} respectively. This yields the following theorem:

\begin{theorem}{\label{fluc_est_var_path}}

Let $\hba$ and $\htm$ be as defined in \eqref{hba_defn} and \eqref{htm_defn} respectively and let $\hat{\beta}_{A,j}$ and $\hat{\theta}_{M,j}$ be the $j^{th}$ coordinates of $\hba$ and $\htm$ respectively. Further let $\hat{\mathfrak{t}}_M$ and  $\hat{\mathfrak{r}}$ be the estimators for $\mathfrak{t}_M$ and $f^{\Sigma_B}_5(c_2)$ as given in Theorem \ref{fluc_est_var}. 
   Then under Assumptions \ref{ass:rows}-\ref{ass:path} and 
\begin{align}
\hat{\sigma}_j^2:=&\hat{\mathfrak{r}}\hat{\theta}^2_{M,j}+\psi_1(\kappa_1,c_1,c_2)\hat{\beta}^2_{A,j}+\psi_2(\kappa_1,c_1,c_2)\hat{\beta}^2_{A,j}(\hat{\theta}^2_{M,j}+\hat{\mathfrak{t}}_M)\label{est_var_sigma_j}
\end{align}
the following convergence holds: 
\begin{equation}\hat{\sigma}_{j}^2\xrightarrow{P}\sigma^2_{j}.\label{sigma_path_consis}\end{equation}
\end{theorem}

Consequently, $2\hat{\sigma}_j^2$ and $\hat{\sigma}_j^2$
consistently estimate the limiting variances in
Theorems~\ref{thm_plugin_path} and~\ref{thm_cf_path}, respectively.

\section{Numerical experiments\label{sec:num_exp}}

We illustrate the finite-sample behavior of the proposed estimators. Throughout this section, the total sample size is fixed at $m=2000$, so that each subsample has size $n=1000$. The mediator dimension is $p=\lfloor n c_1(1-c_2)\rfloor$, while the covariate dimension is $d=\lfloor n c_2\rfloor-1$. All results are based on $500$ independent Monte Carlo replications.

The numerical experiments examine the Gaussian approximations
predicted by Theorems~\ref{thm_plugin}, \ref{thm_cf},
\ref{thm_plugin_path}, and~\ref{thm_cf_path}, together
with the variance-consistency results of
Theorems~\ref{fluc_est_var} and~\ref{fluc_est_var_path}.
We use normal QQ-plots of the corresponding standardized
estimators to assess their finite-sample behavior.

\subsection{Simulation framework \label{sim_frame}}

We first describe the common data-generating mechanism and estimation procedure used throughout the numerical experiments.

For $i=1,\ldots,m$, we generate independent samples $(Y_i,M_i,C_i,A_i)_{i=1}^{2000}$ according to the LSEM \eqref{y_mod_vec} and \eqref{M_mod_vec}
where the errors are generated as described in Assumption \ref{ass:errors}. Further $A_i$ and the entries of $C_i$ are generated as independent standard Gaussian random variables.
The sample is split into two equal parts
$\mathfrak{D}_a$ and
$\mathfrak{D}_b$.

Using $\mathfrak D_a$, we compute the partial-ridge estimator
$\btm$ from~\eqref{btm_eqn_1}, with the tuning parameter
specified in Assumption~\ref{ass:tune}. In the experiments we use the limiting bias-correction factor
$\mathfrak c_1$ from~\eqref{fin_cin_conv} and set
$\htm=\mathfrak c_1\btm$. Appendix~\ref{sec:centering_correction}
justifies this replacement in the Gaussian setting.
Using $\mathfrak D_b$, we compute $\hba$ by ordinary
least squares as in~\eqref{hba_defn}. We then reverse
the fold assignments to obtain the corresponding
cross-fitted estimators.

\subsection{Global mediation effect \label{glob_exp}}

For the global mediation effect, the structural parameters are generated
once and subsequently held fixed across the Monte Carlo replications.
The vectors $\ba,\tm\in\mathbb R^p$ have independent entries distributed
as $X/\sqrt p$, where $X\sim\operatorname{Bernoulli}(1/2)$. The matrix
$\beta_C\in\mathbb R^{p\times d}$ has independent entries distributed as
$X/\sqrt{pd}$, and the vector $\theta_C\in\mathbb R^d$ has independent
entries distributed as $X/\sqrt d$, with all Bernoulli variables mutually
 independent.
 The scalar parameter $\theta_A$ is fixed at $5$.

 Recall the plug-in estimator of the global mediation effect in
Algorithm~\ref{alg:mediation} and its cross-fitted version in
Section~\ref{crossfit_defn}. Let $\widehat\sigma^2$ be the estimator
in \eqref{lim_var_est}. We use
$\sqrt{m}(\widehat\gamma-\gamma)/
(\sqrt{2}\widehat\sigma)$ for the plug-in estimator and
$\sqrt{m}(\widehat\gamma_{cf}-\gamma)/\widehat\sigma$
for the cross-fitted estimator.
The QQ-plots in Figure~\ref{fig:global_qq} show strong
agreement with the standard Gaussian distribution across
all three regimes. The near-linearity of the plots
suggests that the proposed standardization captures the
finite-sample fluctuations well. The same qualitative
behavior is observed for both the plug-in and cross-fitted
estimators.

\input{plots_main/figure_global.tex}

\subsection{Pathwise mediation effect}

For the pathwise experiment, the nonspecial coordinates of $\ba$ and
$\theta_M$ are generated according to the same Bernoulli scheme as in the
global experiment, while the first coordinates are fixed at 
$\beta_{A,1}=\theta_{M,1}=1.$
Thus $\eta_1=1$ throughout the experiment, while all remaining coordinates
remain of order $p^{-1/2}$. The parameters $\beta_C$, $\theta_C$, and
$\theta_A$ are generated as in the global experiment (cf. Section \ref{glob_exp}).

Recall the plug-in estimator in Algorithm~\ref{alg:mediation_path}
and its cross-fitted version in Section~\ref{sec:cross_path}.
Let $\widehat\sigma_1^2$ be the estimator in
\eqref{est_var_sigma_j}. We use
$\sqrt{m}(\widehat\eta_1-\eta_1)/
(\sqrt{2}\widehat\sigma_1)$ for the plug-in estimator and
$\sqrt{m}(\widehat\eta_{1,cf}-\eta_1)/\widehat\sigma_1$
for the cross-fitted estimator.
Figure~\ref{fig:pathwise_qq} reports the corresponding QQ-plots under the same three regimes used above. The QQ-plots of $\widehat{\eta}_1$ and $\widehat{\eta}_{1,cf}$ after proper centering and scaling  are consistent with the Gaussian approximation predicted by
the asymptotic theory. This indicates that the asymptotic Gaussian approximation remains accurate not only for the aggregate mediation effect, but also for the coordinate-level pathwise effect, where the target parameter depends on individual coordinates.

The limiting variance factors $\sigma^2$ and $\sigma_j^2$
depend nontrivially on $\kappa_1$. Their behavior and numerical
minimization over a grid are discussed in
Appendix~\ref{sec:num_exp_extra}.

\input{plots_main/figure_pathwise.tex}

\subsection{Robustness to non-Gaussian errors}

We finally examine whether the Gaussian approximation remains accurate when the Gaussian error variables are replaced by non-Gaussian variables with matching first two moments. This experiment is intended as an empirical robustness check rather than a rigorous universality result, which can be investigated in the future. We keep our simulation setting the same as above but change the error distributions. 
Specifically, we consider three error configurations:
\begin{align*}
\textbf{Case I:}\quad
& \epsilon_{ij}\sim X,\qquad
\delta_i\sim Y,\qquad
X,Y\text{ i.i.d.\ Rademacher},\\
\textbf{Case II:}\quad
& \epsilon_{ij}\sim X,\qquad
\delta_i\sim N(0,1),\qquad
X\text{ Rademacher},\\
\textbf{Case III:}\quad
& \epsilon_{ij}\sim N(0,1),\qquad
\delta_i\sim Y,\qquad
Y\text{ Rademacher}.
\end{align*}
All variables are generated independently across observations, and for each $i$, the coordinates $\epsilon_{ij}$ are mutually independent and independent of $\delta_i$.

Since the Rademacher distribution is centered with unit variance, each configuration matches the baseline Gaussian setting at the level of the first two moments.
Figure~\ref{fig:global_nongaussian_qq} reports QQ-plots
for the standardized global mediation-effect estimators
under the three non-Gaussian error configurations.
The plots remain close to linear across the cases
considered, suggesting that the Gaussian approximation
may extend beyond the Gaussian error model.
Corresponding results for the pathwise mediation effect
are presented in Appendix~\ref{sec:num_exp_extra}.

\input{plots_univ/figure_global_nongaussian.tex}

Overall, the simulations illustrate good agreement with
the Gaussian approximation in the high-dimensional
settings considered, for both global and pathwise
mediation effects. The non-Gaussian experiments provide
further exploratory evidence of robustness under the
error perturbations examined above.
\section{Proof Outline \label{sec:proof_outline}}

In this section, we discuss the main ideas behind the proof of the theorems stated in Section \ref{sec:main_results}. For the proof outline and the appendix calculations, we retain
the fold normalization $\sqrt{n}$, where $n=m/2$.
The main theorems are obtained from these calculations by
multiplying the centered statistics by $\sqrt{m/n}=\sqrt{2}$,
and hence multiplying their limiting variances by two.

As discussed above, a central challenge of our approach is
to debias $\tm$ without prior knowledge of the covariance
matrix of the full outcome-regression design.
Standard debiasing constructions often require covariance
information or structural assumptions on target-specific
precision directions; see
\cite{cai_guo,javan_mont3,BZ_dof}.
In our setting, the covariance matrix $\Sigma$ of
$(M_i^\top,A_i,C_i^\top)^\top$ depends on the unknown
mediator-model coefficients, including $\ba$, through
\eqref{M_mod_vec}.
Consequently, standard debiasing procedures are not
directly applicable under our assumptions.

We bypass this problem by using the partial ridge estimator \eqref{btm_eqn_1}. A crucial part of our debiasing technique is the alternate form of the partial estimator as shown in Lemma \ref{partial_ridge_lemma}, which states that the estimator $\btm$ has the following expansion:
\begin{equation}\breve{\theta}_M=\tm-\lambda_1H_{\lambda_1,a}^{-1}\tm+ H_{\lambda_1,a}^{-1}E_a^TP^{\perp}_{B_a}D_a.\label{part_ridge_simplified}\end{equation}

We note that the third term in the RHS has mean $0$ and the bias stems from the second term. Further the term $H^{-1}_{\lambda_1,a}=(\bea^T\bea+\lambda_1 I)^{-1}$ (a function of $E_a,B_a$) does not depend on unknown parameters (since we can assume $V=I$). Using the theory of deterministic equivalents, we characterize the asymptotic behaviour of the resolvent $H_{\lambda_1,a}^{-1}$ which helps in debiasing the estimator $\btm$ asymptotically.

\subsection{Deterministic Equivalent}

Deterministic equivalents provide a framework for approximating random matrices by non-random surrogates whose action on deterministic vectors and matrices is asymptotically equivalent. Unlike classical spectral methods, which primarily describe limiting eigenvalue distributions, this approach yields precise control of bilinear forms that are central to high-dimensional inference.

Formally, a matrix $\bar Q_n$ is a deterministic equivalent of a symmetric random matrix $Q_n$ if, for deterministic matrices $A_n$ with $\sup_n\|A_n\|<\infty$ and unit-norm vectors $a_n,b_n$,
\[
\frac{1}{n} \tr\!\big(A_n(Q_n-\bar Q_n)\big)\to 0,
\qquad
a_n^\top(Q_n-\bar Q_n)b_n \to 0,
\]
in probability or almost surely, according to the result under
consideration. The mode of convergence will be stated explicitly
whenever relevant.

This framework is particularly useful for analyzing resolvent-based quantities and deriving bias corrections and asymptotic normality results; see \cite{Liao_RMT, yin2023limiting}. We refer the reader to \cite{Liao_RMT} for a comprehensive account of the historical development and applications of deterministic equivalents. The analysis presented here highlights their emerging role as a core technical tool in high-dimensional causal inference.

In the present work, deterministic equivalents are used to analyze the bias and fluctuation behavior of partial ridge estimators in a proportional asymptotics regime. By expressing the estimator through resolvent-type functionals of sample covariance matrices (cf. \eqref{part_ridge_simplified}), we are able to derive approximate deterministic expressions for the bias terms and construct asymptotically bias-corrected estimators without requiring full knowledge of the design covariance structure---a perspective well-suited to the mediation setting studied here.

\subsection{Asymptotic debiasing using deterministic equivalents \label{asymp_debias}}

Defining $\kappa_1=\lambda_1/(n-d-1)$, the scaled resolvent is defined as
$$H^{-1}_{\kappa_1,a}:=(\bea^T\bea/(n-d-1)+\kappa_1I_p)^{-1},$$
so that $\lambda_1 H^{-1}_{\lambda_1,a}=\kappa_1H^{-1}_{\kappa_1,a}$. To state the deterministic equivalent for $H^{-1}_{\kappa_1,a}$, first define
\begin{equation}
r_n:=n-d-1,
\qquad\fin
:=
\frac{1}
{\kappa_1\{1+g^{0}_{2,r_n}(-\kappa_1)\}},
\label{finite_center_estimator}
\end{equation}
where $g^{0}_{2,r_n}$ is the finite-dimensional deterministic-equivalent
quantity defined in Appendix~\ref{sec:aux_lemma}.

The deterministic equivalent for $H^{-1}_{\kappa_1,a}$ is then given by $\fin I$ (cf. \cite{yin2023limiting}). 
Together with the relation \eqref{part_ridge_simplified}, this yields:
\begin{equation}
\btm\approx \tm-\kappa_1 \fin\tm+ H_{\lambda_1,a}^{-1}E_a^TP^{\perp}_{B_a}D_a =(1-\kappa_1\fin)\tm + H_{\lambda_1,a}^{-1}E_a^TP^{\perp}_{B_a}D_a \label{approx_btm}
\end{equation}
where the approximation $\approx$ is heuristic. Indeed, the deterministic-equivalent approximation holds at the level of traces and bilinear forms, but not at the matrix or vector level. We therefore use this heuristic only to motivate the estimator of the mediation effect; its asymptotic properties are established later directly through the corresponding bilinear forms. 
 Recall that the bias-correction factor $\cin$ and the estimator for $\tm$ were given by
\begin{equation}
\cin
:=
\left\{1-\kappa_1\fin\right\}^{-1}, \quad 
\htm=\cin\btm.
\label{cin_htm}
\end{equation}

Further, we define $(1-f_1(\kappa_1,c_1)\kappa_1)^{-1}=\mathfrak{c}_1$
and show in the appendix that
\begin{equation}
\fin\rightarrow f_1(\kappa_1,c_1),
\qquad
\cin\rightarrow
\mathfrak c_1:=
\{1-\kappa_1f_1(\kappa_1,c_1)\}^{-1}.\label{fin_cin_conv}
\end{equation}
Thus $f_1(\kappa_1,c_1)$ and $\mathfrak c_1$ will appear in the limiting variance
formulas, whereas the estimator itself uses the finite-$n$ corrections with $\fin$ and 
$\cin$. Thus, from \eqref{approx_btm} and \eqref{cin_htm},  we obtain:
\begin{equation}
\cin\btm\approx \tm + \cin H_{\lambda_1,a}^{-1}E_a^TP^{\perp}_{B_a}D_a \label{approx_debias}
\end{equation}
Here, for random vectors $X_n$ and $Y_n$, the notation
$X_n \approx Y_n$ means that
\[
    a_n^\top (X_n - Y_n) \xrightarrow{P} 0
\]
for every deterministic sequence of unit vectors $a_n$ of compatible
dimensions.
Since\\ $\cin H_{\lambda_1,a}^{-1}E_a^TP^{\perp}_{B_a}D_a$ has mean $0$, $\htm=\cin\btm$ serves as the \emph{asymptotically bias-corrected estimator} for $\tm$.

\subsection{Proof highlights of the main theorems}

 Next we outline the main steps in the proof of asymptotic normality of the plug-in (Theorem \ref{thm_plugin})  estimator of the global mediation effect. 
A natural starting point for proving the asymptotic normality of the plugin estimator is the simplified representation of the partial ridge estimator $\btm$ in \eqref{part_ridge_simplified}, combined with the deterministic equivalent in \eqref{approx_btm}. However, the approximation in \eqref{approx_btm} resulting from $H^{-1}_{\kappa_1,a}
    \approx
    f_1(\kappa_1,c_1)I$ is only heuristic.
More precisely, for
deterministic unit vectors $u$ and $v$, we have
\[
    u^\top H^{-1}_{\kappa_1,a}v
    =
    f_1(\kappa_1,c_1)\,u^\top v
    +
    o_p(1).
\]
Consequently, the debiasing step in \eqref{approx_debias} cannot be
justified by directly substituting the scalar deterministic equivalent
at the matrix level. To use the approximation at the level of bilinear forms we begin from the exact identity
\[
\cin \hba^\top \breve{\theta}_M
=
\cin \hba^\top \theta_M
-
\lambda_1 \cin \hba^\top H_{\lambda_1,a}^{-1} \theta_M
+
\cin \hba^\top H_{\lambda_1,a}^{-1} E_a^\top P_{B_a}^\perp D_a,
\]
and show that the following expansion holds:
\[
\sqrt{n}
\big(\hba^\top \hat{\theta}_M
-
\ba^\top \theta_M
\big)
=
\Delta_1 + \Delta_2 + \Delta_3 .
\]

The exact description of $\Delta_i$s and proof steps are given in the appendix. The three terms will be analyzed separately via conditioning arguments and deterministic equivalent techniques.  In this section we provide an overview of these conditioning techniques.

\medskip
\noindent
\textit{Term $\Delta_1 = \Delta_1(B_a,E_a,B_b,E_b,D_a)$.}
Conditioning on $(B_a,E_a,B_b,E_b)$, we observe that the  randomness arises from $D_a$, which is Gaussian. Hence $\Delta_1$ is conditionally Gaussian. Its asymptotic variance is obtained by  showing that the conditional variance  converges to a constant.

\medskip
\noindent
\textit{Term $\Delta_2 = \Delta_2(B_a,E_a,B_b,E_b)$.}
Conditioning on $(B_a,E_a,B_b)$, the randomness is driven by $E_b$. Gaussianity again yields conditional normality, and deterministic equivalents provide the limit of the conditional variance.

\medskip
\noindent
\textit{Term $\Delta_3 = \Delta_3(B_a,E_a)$.}
This term involves fluctuations of bilinear forms of resolvents. Its asymptotic normality follows from the central limit theorem for such forms established in \cite{yin2023limiting}.

\medskip

In the final step the three components are combined to yield the claimed Gaussian limit.

The proof of Theorem \ref{thm_cf} proceeds by the same decomposition into $\Delta$s as before, however the terms from the two permutations $\vartheta_0$ and $\vartheta_1$ are potentially dependent, hence a much finer analysis is needed. To mitigate this, we use the fact that $PE$ and $P^{\perp}E$ are independent where  $E$ is a matrix with i.i.d. Gaussian entries, $P$ a deterministic orthogonal projection matrix, and $P^{\perp}$ is its projection onto the orthogonal subspace.

Recall the definition of $\Delta_i$-s in the previous section. We use $\Delta_{i,\vartheta}$ to denote the  $\Delta_i$ for a given permutation $\vartheta$. In the case of cross-fitted estimator we use a different decomposition
\begin{align*}T_1:=\frac{1}{2}\sum_{\vartheta}\Delta_{1,\vartheta},\quad T_2:=\frac{1}{2}\sum_{\vartheta}\left(\Delta_{2,\vartheta}+\Delta_{3,\vartheta}\right)
\end{align*}
so that
\begin{equation*}
    \sqrt{n}(\hat{\gamma}_{cf}-\gamma)=T_1+T_2.
\end{equation*}

We again use conditioning arguments and deterministic equivalents
to analyze $T_1$ and $T_2$. Because the fold-swapped estimators reuse
the same observations, their joint analysis requires care. We defer
the details to Appendix~\ref{sec:proof_cf}.

\section{Discussion}

This paper develops a framework for inference on global and pathwise mediation effects in high-dimensional linear structural equation models under a proportional asymptotics regime without any sparsity assumption on the signals. By combining partial ridge regression with an asymptotic debiasing principle, we construct estimators that admit tractable limiting distributions even when both the mediator and covariate dimensions grow proportionally with the sample size. This work shows that the assumption of known covariance structure can be substantially weakened under additional structural constraints—here arising from a linear structural equation model. Our results thus contribute to the growing literature on high-dimensional causal inference by demonstrating that carefully exploiting model structure can mitigate otherwise restrictive covariance assumptions. 

From a broader perspective, our results contribute to the emerging literature at the intersection of high-dimensional statistics and causal inference. While high-dimensional debiasing techniques have been extensively studied in regression and generalized linear models, their adaptation to causal mediation settings—particularly under proportional asymptotics—has remained limited. Our analysis shows that structural equation modeling assumptions can be fruitfully exploited to overcome some of the intrinsic challenges posed by high dimensionality, and suggests a general strategy for importing tools from random matrix theory and high-dimensional statistics into causal inference.

Several limitations of the present work suggest directions for future research. First, our procedure requires $d+1<n$, where $n$ is the sample size of each fold. When the exposure-covariate dimension exceeds the fold sample size,
the corresponding Gram matrix is singular. If the design has full
row rank, projection onto its column space also leaves no
residualized mediators. Consequently, our partial-ridge construction and the associated debiasing corrections are no longer applicable. However, such high-dimensional covariate regimes are common in contemporary applications. Developing principled debiasing methodologies that remain valid when the exposure-covariate dimension exceeds the fold sample size therefore represents a compelling and practically relevant extension of this work.

A second restriction concerns the assumption of a known
mediator-error covariance matrix $V$. This assumption plays a crucial role in our deterministic equivalent approximations and, in turn, the asymptotic debiasing of the target parameters. Relaxing this assumption would require new techniques
to account for the unknown mediator-error covariance
in the debiasing and variance analysis. We conjecture that this problem is intrinsically challenging, as it is closely tied to the broader difficulty of debiasing under unknown and high-dimensional covariance structure—precisely the issue that motivates our structural modeling approach. A deeper understanding of this interaction remains an important open problem.

Finally, our analysis relies on Gaussian errors $E$ in Assumption \ref{ass:errors}, with rotational invariance playing a central role in both the construction of the initial plug-in estimators and the analysis of the cross-fitted debiased estimators. It is natural to conjecture that suitably formulated universality results may continue to hold beyond the Gaussian setting. However, establishing such results in the mediation model is highly nontrivial. In particular, the mediation model induces an intrinsic dependence structure in the design matrix, complicating a direct extension of existing universality results. While universality has been shown to hold under certain classes of dependent designs \cite{lahiry2024universality}, it is also known to fail in more general settings \cite{elkaroui_ptrf, han2023universality}. Developing a rigorous universality theory for the present framework-potentially via adaptations of Gaussian interpolation techniques-lies beyond the scope of this paper and is left for future work.

\section*{Acknowledgements} P.S.~gratefully acknowledges support from the Sloan Research Fellowship
in Mathematics. P.S.~and R.M.~were also supported by National Science
Foundation CAREER Awards DMS-2440824 and DMS-2338760, respectively.

\section*{Disclosure of AI Use}

GPT 5 and 6 were used to assist with reviewing the manuscript and copyediting. The authors  take full responsibility for everything in this manuscript.

\bibliographystyle{plain}
\bibliography{mediation}

\appendix

\section{Deterministic Equivalent and exact formula for the limiting variances \label{sec:var_details}}

\subsection{Deterministic Equivalents}
The following definition of deterministic equivalent is adapted from \cite{Liao_RMT}

\begin{definition}
 We say that $\bar{Q}_n\in \mathbb{R}^{n\times n}$ is a deterministic equivalent for the symmetric random matrix $Q_n\in \mathbb{R}^{n\times n}$ if, for every sequence of deterministic matrices $A_n\in\mathbb{R}^{n\times n}$ with $\sup_n\|A_n\|<\infty$ and
deterministic unit vectors $a_n,b_n\in\mathbb{R}^n$, we have, as $n\to\infty$,

\begin{align*}
    \frac{1}{n}\tr[A_n(Q_n-\bar{Q}_n)]&\rightarrow 0\\
    a_n^T(Q_n-\bar{Q}_n)b_n&\rightarrow 0
\end{align*}
where the convergence is either in probability or almost sure.
\end{definition}

\begin{theorem}[Mar\v{c}enko--Pastur deterministic equivalent;
Theorem~2.4 of \cite{Liao_RMT}]
\label{mp_detequ}
Let
$
X=[X_1,\ldots,X_\n]\in\mathbb R^{\p\times \n},
$
where the columns $X_i$ are i.i.d., and the coordinates of each
$X_i$ are independent with mean zero and variance one. Assume, for
some $\varepsilon>0$, that
\[
\sup_{\p}\max_{1\le j\le \p}
\E |X_{ij}|^{4+\varepsilon}<\infty.
\]
Suppose that
$
\frac \p\n\rightarrow c\in(0,\infty)$. For a fixed $z\in\mathbb C\setminus[0,\infty)$, define the resolvent
\[
Q_\n(z)
:=
\left(
\frac1\n XX^\top-zI_\p
\right)^{-1}.
\]
Then the deterministic equivalent of 
$
Q_\n(z)$ is given by $m_c(z)I_\p,$
where $m_c(z)$ is the Stieltjes transform of the
Mar\v{c}enko-Pastur distribution with aspect ratio $c$.

More explicitly, for every sequence of deterministic matrices
$A_\p\in\mathbb R^{\p\times \p}$ satisfying
$\sup_\p\|A_\p\|<\infty$, and every pair of deterministic vectors
$u_\p,v_\p\in\mathbb R^\p$ satisfying
\[
\sup_\p\|u_\p\|<\infty,
\qquad
\sup_\p\|v_\p\|<\infty,
\]
we have
\[
\frac1\p
\operatorname{Tr}
\left[
A_\p
\left\{
Q_\n(z)-m_c(z)I_\p
\right\}
\right]
\xrightarrow{\mathrm{a.s.}}0
\]
and
\[
u_\p^\top
\left\{
Q_\n(z)-m_c(z)I_\p
\right\}
v_\p
\xrightarrow{\mathrm{a.s.}}0.
\]

The function $m_c(z)$ is characterized as the unique solution of
\begin{equation}
zc\,m_c(z)^2
-
(1-c-z)m_c(z)
+
1
=
0
\label{eq:MP_fixed_point}
\end{equation}
satisfying the Stieltjes-transform condition
\[
\Im(z)\Im(m_c(z))>0
\qquad\text{when }\Im(z)\neq0,
\]
and
\[
m_c(z)>0
\qquad\text{when }z<0.
\]
\end{theorem}

For $z=-\kappa$ with $\kappa>0$, the Stieltjes-transform condition
selects the unique positive solution of \eqref{eq:MP_fixed_point},
namely
\begin{equation}
m_c(-\kappa)
=
\frac{
-(1-c+\kappa)
+
\sqrt{(1-c+\kappa)^2+4c\kappa}
}{
2c\kappa
}
=
f_1(\kappa,c).
\label{eq:MP_m_f1}
\end{equation}

\subsection{Limiting variance of the global mediation effect}
To define the limiting variance factors $\sigma^2$  we need to define several quantities related to the deterministic equivalent of the resolvent of the sample covariance matrix.

For a matrix $X\in\mathbb R^{\n\times\p}$ satisfying the
assumptions of Theorem~\ref{mp_detequ} after transposition, and for
$z>0$, the deterministic equivalent of
$(X^\top X/\n+zI_\p)^{-1}$ is $f_1(z,c)I_\p$, where
$c=\lim_{\n\to\infty}\p/\n$ and

\[
f_1(z,c)
=\frac{-(1-c+z)+\sqrt{(1-c+z)^2+4cz}}{2cz},
\qquad z>0.
\]

We introduce a collection of auxiliary functions defined in terms of $f_1(z,c)$ that will be used to characterize
the limiting variance of the proposed mediation effect estimators. These functions
are defined in two stages.

\vspace{0.3em}

\paragraph{Intermediate auxiliary functions.}
The first group consists of technical transforms derived from the function
$f_1(z,c)$. These functions do not appear directly in the final variance expressions,
but serve as intermediate quantities in the construction of the objects defined below. For the resolvent notation below, let $z<0$. We set

\begin{align*}
\underline f_1(z,c)
    &:= -\frac{1-c}{z}+c f_1(-z,c), \\
g_1(z,c)
    &:= c f_1(-z,c), \\
g_2(z,c)
    &:= -z^{-1}\bigl(f_1(-z,c)\bigr)^{-1}-1, \\
d(z,c)
    &:= \frac{1}{z^2}
       \left(\frac{g_1(z,c)}{\partial_z g_1(z,c)}+z\right)
       \left(\frac{g_2(z,c)}{\partial_z g_2(z,c)}+z\right).
\end{align*}

\paragraph{Variance-defining functions.}
Using the intermediate quantities above, we now define the functions that enter
explicitly into the asymptotic variance characterization. For $z<0$, define

\begin{align}
h_1(z,c)
    &:= \frac{c}{z^4\bigl(1-d(z,c)\bigr)}
       \left[
       \frac{g_2(z,c)}{\partial_z g_1(z,c)}
       + \frac{z\,\partial_z g_2(z,c)}{\partial_z g_1(z,c)}
       \right],
            \label{h-kernels}\\
f_2(\kappa,c)
    & := f_1(\kappa,c)+\kappa\,\partial_{\kappa} f_1(\kappa,c),\nonumber
\end{align}
which will appear repeatedly in the subsequent variance expressions.

Recall the definitions of $c_1, c_2$, and $\kappa_1$ as follows:

$$c_1=\lim_{n\rightarrow \infty}\frac{p}{n-d-1},\quad c_2=\lim_{n\rightarrow \infty}\frac{d+1}{n},\quad \kappa_1=\frac{\lambda_1}{n-d-1}.$$

Using  the definition of the constants $\mathfrak{b}_A, \mathfrak{t}_M,\mathfrak{g}$ (see Assumption \ref{ass:signal}) appearing as limits of the signal strengths, the limiting value of the function of the covariance matrix $f^{\Sigma_B}_5(c_2)$, and the functions $h_1(.,.)$ and $f_2(.,.)$ as defined in \eqref{h-kernels}, we are ready to define the limiting variance factor of the plugin estimator for the global mediation effect:

\begin{align}
\sigma^2=&(1-\kappa_1f_1(\kappa_1,c_1))^{-2}\left[ f_2(\kappa_1,c_1)(c_1 f^{\Sigma_B}_5(c_2)+(1-c_2)^{-1} \mathfrak{b}_A)\right.\nonumber\\
&+(1-2\kappa_1f_1(\kappa_1,c_1)-\kappa^2_1\partial_{\kappa}f_1(\kappa_1,c_1))f^{\Sigma_B}_5(c_2) \mathfrak{t}_M\nonumber\\
&\left.+\kappa^2_1\frac{(f_1(\kappa_1,c_1)\kappa_1)^{4}}{c_1(1-c_2)}\left(h_1(-\kappa_1,c_1)\mathfrak{g}^2+h_1(-\kappa_1,c_1)\mathfrak{b}_A\mathfrak{t}_M\right)\right].\label{lim_var_defn_explicit}
\end{align}

For $\kappa_1>0$ and $c_1>0$, Lemma~\ref{func_resolv_2} gives
\[
f_2(\kappa_1,c_1)
=\int\frac{x}{(x+\kappa_1)^2}\,d\mu_{MP,c_1}(x)>0.
\]
Moreover, the remaining contributions in~\eqref{lim_var_defn_explicit} are nonnegative. In particular, the term involving $\mathfrak b_A$ is strictly positive because $\mathfrak b_A>0$ and $c_2<1$. Therefore $\sigma^2>0$. We now write the same variance factor in abbreviated form:

\begin{align*}
\sigma^2(\kappa_1,c_1,c_2,f^{\Sigma_B}_5(c_2),\mathfrak{b}_A,\mathfrak{t}_M,\mathfrak{g})=&\varphi_1(\kappa_1,c_1,c_2)f^{\Sigma_B}_5(c_2)+\varphi_2(\kappa_1,c_1,c_2)\mathfrak{b}_A\\
+\varphi_3(\kappa_1,c_1,c_2)f^{\Sigma_B}_5(c_2)\mathfrak{t}_M
&+\varphi_4(\kappa_1,c_1,c_2)(\mathfrak{b}_A\mathfrak{t}_M+\mathfrak{g}^2),
\end{align*}

with the functions $\varphi_j()$s given by

\begin{align*}
    \varphi_1(\kappa_1,c_1,c_2)&=(1-\kappa_1f_1(\kappa_1,c_1))^{-2}f_2(\kappa_1,c_1)c_1\\
    \varphi_2(\kappa_1,c_1,c_2)&=(1-\kappa_1f_1(\kappa_1,c_1))^{-2}\frac{f_2(\kappa_1,c_1)}{1-c_2}\\
    \varphi_3(\kappa_1,c_1,c_2)&=(1-\kappa_1f_1(\kappa_1,c_1))^{-2}(1-2\kappa_1f_1(\kappa_1,c_1)-\kappa^2_1\partial_{\kappa}f_1(\kappa_1,c_1))\\
    \varphi_4(\kappa_1,c_1,c_2)&=(1-\kappa_1f_1(\kappa_1,c_1))^{-2}\kappa^2_1\frac{(f_1(\kappa_1,c_1)\kappa_1)^{4}}{c_1(1-c_2)}h_1(-\kappa_1,c_1)\\
\end{align*}

\subsection{Limiting variance of the pathwise mediation effect}
The limiting variance factor in this case also depends on the functions $h_1(.,.)$ and $f_1(.,.)$ as in the previous section. However,  instead of the limiting signal strength $\mathfrak{b}_A$, the coordinates of the signal $\beta_{A,j}$ and $\theta_{M,j}$ appear in the limiting variance factor. The formula is given by:

\begin{align}
\sigma^2_j
:=& (1-\kappa_1f_1(\kappa_1,c_1))^{-2}\Big((1-c_2)^{-1} f_2(\kappa_1,c_1)\baj^2+(1-\kappa_1f_1(\kappa_1,c_1))^2f^{\Sigma_B}_5(c_2) \tmj^2\Big.\nonumber\\
&\Big.+\kappa_1^2\frac{\baj^2(f_1(\kappa_1,c_1)\kappa_1)^{4}}{c_1(1-c_2)}h_1(-\kappa_1,c_1)(\tmj^2+\mathfrak{t}_M)\Big)\label{lim_var_sigma_j_explicit}
\end{align}

 The abbreviated expression is:

\begin{align*}
&\sigma_j^2(\kappa_1,c_1,c_2,f^{\Sigma_B}_5(c_2),\beta_{A,j},\theta_{M,j},\mathfrak{t}_M)\\
=&f^{\Sigma_B}_5(c_2)\theta^2_{M,j}+\psi_1(\kappa_1,c_1,c_2)\beta^2_{A,j}+\psi_2(\kappa_1,c_1,c_2)\beta^2_{A,j}(\theta^2_{M,j}+\mathfrak{t}_M)
\end{align*}
with the functions $\psi_k()$s given by

\begin{align*}
    \psi_1(\kappa_1,c_1,c_2)&=(1-\kappa_1f_1(\kappa_1,c_1))^{-2}(1-c_2)^{-1} f_2(\kappa_1,c_1)\\
    \psi_2(\kappa_1,c_1,c_2)&=(1-\kappa_1f_1(\kappa_1,c_1))^{-2}\kappa_1^2\frac{(f_1(\kappa_1,c_1)\kappa_1)^{4}}{c_1(1-c_2)}h_1(-\kappa_1,c_1)
\end{align*}

Note that, under Assumption~\ref{ass:path}, $\sigma_j^2>0$. Indeed, $f_5^{\Sigma_B}(c_2)>0$ makes the first term strictly positive when $\theta_{M,j}\ne0$. If $\theta_{M,j}=0$, the same assumption gives $\beta_{A,j}\ne0$, and the term $\psi_1(\kappa_1,c_1,c_2)\beta_{A,j}^2$ is strictly positive because $\psi_1>0$. 

\section{Proofs of Main Theorems \label{sec:proofs}}

In this appendix we prove the main theorems. The technical lemmas are proved in Appendix~\ref{sec:tech_lemma}.

\medskip

\textbf{Normalization convention.}
Throughout this appendix, $n=m/2$ and all intermediate
fluctuations retain the $\sqrt{n}$ normalization. In particular,
$\sigma^2$ and $\sigma_j^2$ denote the fold-normalized plug-in
variance factors, and the notation
$\sigma_{cf}^2=\sigma^2/2$ and
$\sigma_{j,cf}^2=\sigma_j^2/2$ is used only for the corresponding
fold-normalized cross-fitted variances. For any centered
statistic $T_m$ and any $v\geq0$,
\[
\sqrt{n}\,T_m\xrightarrow{d}N(0,v)
\quad\Longleftrightarrow\quad
\sqrt{m}\,T_m\xrightarrow{d}N(0,2v).
\]
Thus the four fold-normalized limits proved below are equivalent
to the four total-sample limits stated in
Theorems~\ref{thm_plugin}, \ref{thm_cf},
\ref{thm_plugin_path}, and~\ref{thm_cf_path}.
\subsection{Proofs for the plug-in estimator \label{sec:proof_plugin}}

In this subsection we discuss the proof of Theorem \ref{thm_plugin}. The proof proceeds by decomposing $\sqrt{n}(\hat{\gamma}-\gamma)$ into three components and analyzing their conditional distributions. We sequentially condition on the sigma-fields generated by $\{B_a,E_a,B_b,E_b\}$ and $\{B_a,E_a,B_b\}$. In the final step, we study the remaining fluctuation as a function of $E_a$, which yields the desired asymptotic normal distribution.

\subsubsection{Decomposition of $\hat{\gamma}$}

Recall that
$$\htm=\cin(\tm-\lambda_1H_{\lambda_1,a}^{-1}\tm+ H_{\lambda_1,a}^{-1}E_a^TP^{\perp}_{B_a}D_a)$$
and 
$$\hba^T=\ba^T+e^T_1(B_b^TB_b)^{-1}B_b^TE_b.$$
We start with a following decomposition: 

\begin{equation}\sqrt{n}(\hat{\gamma}-\gamma)=\Delta_1+\Delta_2+\Delta_3\label{gamma_decomp}\end{equation} with the terms on the RHS of the above equation given by:

\begin{align*}\Delta_1&=\sqrt{n}\cin(\ba^T+e^T_1(B_b^TB_b)^{-1}B_b^TE_b)H_{\lambda_1,a}^{-1}E_a^TP^{\perp}_{B_a}D_a,\\
    \Delta_2&=\sqrt{n}\cin e^T_1(B_b^TB_b)^{-1}B_b^TE_b(\tm-\lambda_1H_{\lambda_1,a}^{-1}\tm),\\
\Delta_3&=\sqrt{n}(\cin\ba^T(I-\lambda_1H_{\lambda_1,a}^{-1})\tm-\ba^T\tm).
\end{align*}

Define the sigma-fields
\[
\mathcal F_{1,n}:=\sigma(B_a,E_a,B_b),
\qquad
\mathcal F_{2,n}:=\sigma(B_a,E_a,B_b,E_b).
\]
We use these same definitions throughout the global and pathwise
proofs. Conditioning on a random matrix or tuple is shorthand for
conditioning on the sigma-field it generates.

Exploiting the Gaussianity of $D_a$, we show that $\Delta_1$, conditional on $\mathcal F_{2,n}$, is Gaussian with variance converging to a deterministic limit. Similarly, using the Gaussianity of the entries of $E_b$, we establish that $\Delta_2$, conditional on $\mathcal F_{1,n}$, is Gaussian with variance converging to another deterministic limit. Finally, we rewrite $\Delta_3$ as an appropriately centered and scaled bilinear form involving a resolvent. Invoking the fluctuation result for bilinear forms proved in \cite{yin2023limiting}, we show that $\Delta_3$ converges in distribution to a Gaussian random variable. Combining the convergence results for the terms $\Delta_i$ with the decomposition \eqref{gamma_decomp} completes the proof of Theorem~\ref{thm_plugin}.

\subsubsection{Conditional distribution of $\Delta_1$}

In this section we consider the  distribution of $\Delta_1$ conditioned on $\mathcal F_{2,n}$ and the limiting value of the variance. Since $D_a\sim N(0,I_n)$, it follows that
\[
\Delta_1\mid\mathcal F_{2,n}
\sim N\bigl(0,V_1(B_a,E_a,B_b,E_b)\bigr).
\]
with \begin{equation}V_1(B_a,E_a,B_b,E_b)=n\cin^2\|\bea H_{\lambda_1,a}^{-1}(\ba +E_b^TB_b(B_b^TB_b)^{-1}e_1)\|^2,\label{V_1_squared_form}\end{equation}
where we have used the notation $\bea=P^{\perp}_{B_a}E_a$. Before analyzing $V_1(B_a,E_a,B_b,E_b)$ further, we introduce proper scaling for the resolvent as follows:
\begin{align*}
  H_{\kappa_1,a}^{-1}&:=(n-d-1)H_{\lambda_1,a}^{-1}=(\bea^T\bea/(n-d-1)+\kappa_1I)^{-1}.
\end{align*}
We show that $V_1(B_a,E_a,B_b,E_b)$ converges to some deterministic quantity $V_4$. Since the proof is involved, we structure the entire proof as a collection of lemmas which has a hierarchical structure. First, we define the following terms:
\begin{align}
   V_2(B_a,E_a,B_b): = &n^{-1}\alpha\mathfrak{c}^2_1e^T_1\left(\frac{B_b^TB_b}{n}\right)^{-1}e_1 \tr\left[H_{\kappa_1,a}^{-1}\frac{\bea^T\bea}{n-d-1} H_{\kappa_1,a}^{-1}\right]\nonumber\\
   +&\alpha\mathfrak{c}^2_1 \ba^TH_{\kappa_1,a}^{-1}\frac{\bea^T\bea}{n-d-1} H_{\kappa_1,a}^{-1} \ba,\nonumber\\
   V_3(B_a,E_a):=& n^{-1}\alpha\mathfrak{c}^2_1f^{\Sigma_B}_5(c_2) \tr\left[H_{\kappa_1,a}^{-1}\frac{\bea^T\bea}{n-d-1} H_{\kappa_1,a}^{-1}\right]+\alpha\mathfrak{c}^2_1 \ba^TH_{\kappa_1,a}^{-1}\frac{\bea^T\bea}{n-d-1} H_{\kappa_1,a}^{-1} \ba,\nonumber\\
   V_4:=& c_1\mathfrak{c}^2_1 f^{\Sigma_B}_5(c_2)f_2(\kappa_1,c_1)+\alpha \mathfrak{c}^2_1f_2(\kappa_1,c_1)\mathfrak{b}_A,\label{V4_def}
\end{align}
where $\alpha:=\lim_{n\rightarrow \infty}\frac{n}{n-d-1}=(1-c_2)^{-1}$. We have the following lemmas: 
\begin{lemma}{\label{V1_Lemma}}
$$V_1(B_a,E_a,B_b,E_b)=V_2(B_a,E_a,B_b)+o_p(1).$$
\end{lemma}

\begin{lemma}{\label{V2_lemma}}
    $$V_2(B_a,E_a,B_b)=V_3(B_a,E_a)+o_p(1).$$
\end{lemma}

\begin{lemma}{\label{V3_Lemma}}
    $$V_3(B_a,E_a)=V_4+o_p(1).$$
\end{lemma}

We note that although $\cin$ appears in the term $V_1(B_a,E_a,B_b,E_b)$, it can be replaced by the limit $\mathfrak{c}_1$ in the subsequent expressions for $V_2(B_a,E_a,B_b)$, $V_3(B_a,E_a)$, and $V_4$.

\subsubsection{Conditional distribution of $\Delta_2$}

Next we consider the distribution of $\Delta_2$ conditioned on $\mathcal F_{1,n}$.
Since the entries of $E_b$ are i.i.d Gaussian, we have
\[
E_b(\tm-\lambda_1H_{\lambda_1,a}^{-1}\tm)
\mid\mathcal F_{1,n}
\sim N\left(0,
\|\tm-\lambda_1H_{\lambda_1,a}^{-1}\tm\|^2 I_n\right).
\]
 so that we have,
\[
\Delta_2\mid\mathcal F_{1,n}
\sim N\bigl(0,W_2(B_a,E_a,B_b)\bigr).\]
with 
$$W_2(B_a,E_a,B_b)=n||(\tm-\lambda_1H_{\lambda_1,a}^{-1}\tm)||^2\cin^2e_1^T(B_b^TB_b)^{-1}e_1.$$
Similar to the term $V_1(B_a,E_a,B_b,E_b)$, $W_2(B_a,E_a,B_b)$ also converges to a deterministic quantity - we call it $W_4$. We will show that $W_2(B_a,E_a,B_b)$ converges in probability to $W_4$ in two steps.
First, define the following quantities:
\begin{align*}
    W_3(B_b)&=n\mathfrak{c}_1^2(1-2\kappa_1f_1(\kappa_1,c_1)-\kappa^2_1\partial_{\kappa} f_1(\kappa_1,c_1))\|\tm\|^2e_1^T(B_b^TB_b)^{-1}e_1\\
    W_4&=\mathfrak{c}_1^2(1-2\kappa_1f_1(\kappa_1,c_1)-\kappa^2_1\partial_{\kappa}f_1(\kappa_1,c_1))f^{\Sigma_B}_5(c_2) \mathfrak{t}_M.
\end{align*}

The next lemma shows the convergence of $W_2(B_a,E_a,B_b)$.
\begin{lemma}{\label{W_approx}}
    We have the following relations:
    \begin{align}
    W_2(B_a,E_a,B_b)&=W_3(B_b)+o_p(1)\label{W_23}.\\
        W_3(B_b)&=W_4+o_p(1).\label{W_34}
    \end{align}
\end{lemma}

\subsubsection{Asymptotic distribution of $\Delta_3$}
Recall that

\begin{align*}\Delta_3&=\sqrt{n}(\cin \ba^T(I-\lambda_1H_{\lambda_1,a}^{-1})\tm-\ba^T\tm)
\end{align*}

We show the convergence of $\Delta_3$, by simplifying $\Delta_3$ to a centered and scaled version of $\ba^TH_{\lambda_1,a}^{-1}\tm$. We note that this is a bilinear form of the  resolvent $H_{\lambda_1,a}^{-1}$. Using a central limit theorem for bilinear forms of resolvents (cf. Lemma \ref{MP_fluc}), we obtain the following: 

\begin{lemma}{\label{delta3_conv_lemma}}

\begin{equation}\Delta_3\xrightarrow{d} N(0,\mathcal{V})
\label{Delta_3_conv}\end{equation}

where
\begin{align}\mathcal{V}=&\kappa^2_1\mathfrak{c}^2_1\frac{(f_1(\kappa_1,c_1)\kappa_1)^{4}}{c_1(1-c_2)}\left(h_1(-\kappa_1,c_1)\mathfrak{g}^2+h_1(-\kappa_1,c_1)\mathfrak{b}_A\mathfrak{t}_M\right).\label{fluc_var}
\end{align}

\end{lemma}

\subsubsection{Asymptotic distribution of $\hat{\gamma}$}

\begin{proof}[Proof of Theorem \ref{thm_plugin}]
Recall the decomposition \eqref{gamma_decomp} we obtain

$$\sqrt{n}(\hat{\gamma}-\gamma)=\Delta_1+\Delta_2+\Delta_3$$

We first combine $\Delta_2$ and $\Delta_3$. Since $\Delta_3$
depends only on $(B_a,E_a)$, it is $\mathcal F_{1,n}$-measurable.
Moreover, conditional on $\mathcal F_{1,n}$, the only randomness in
$\Delta_2$ is due to $E_b$. Therefore,
\[
\Delta_2\mid\mathcal F_{1,n}
\sim
N\!\left(
0,
W_2(B_a,E_a,B_b)
\right).
\]
By Lemma~\ref{W_approx}, $ 
W_2(B_a,E_a,B_b)\xrightarrow{P}W_4,
$
while Lemma~\ref{delta3_conv_lemma} gives
$
\Delta_3\xrightarrow{d}N(0,\mathcal V).
$
Applying Lemma~\ref{lem:cond_gauss_shift_sum}, with
\[
Z_n=(B_a,E_a,B_b),
\qquad
X_n=\Delta_2,
\qquad
Y_n=\Delta_3,
\]
we obtain
\[
\Delta_{2,3}
:=
\Delta_2+\Delta_3
\xrightarrow{d}
N(0,W_4+\mathcal V).
\]

We next combine $\Delta_1$ and $\Delta_{2,3}$. Since both $\Delta_2$
and $\Delta_3$ are functions of $(B_a,E_a,B_b,E_b)$,
$\Delta_{2,3}$ is $\mathcal F_{2,n}$-measurable. Conditional on
$\mathcal F_{2,n}$, the only randomness in $\Delta_1$ is due to
$D_a$. Hence
\[
\Delta_1\mid\mathcal F_{2,n}
\sim
N\!\left(
0,
V_1(B_a,E_a,B_b,E_b)
\right).
\]
By Lemmas~\ref{V1_Lemma}-\ref{V3_Lemma},
$
V_1(B_a,E_a,B_b,E_b)\xrightarrow{P}V_4.
$
A second application of Lemma~\ref{lem:cond_gauss_shift_sum}, now
with
\[
Z_n=(B_a,E_a,B_b,E_b),
\qquad
X_n=\Delta_1,
\qquad
Y_n=\Delta_{2,3},
\]
yields
\[
\Delta_1+\Delta_{2,3}
\xrightarrow{d}
N(0,V_4+W_4+\mathcal V).
\]
Since
\[
\sqrt{n}(\hat{\gamma}-\gamma)
=
\Delta_1+\Delta_{2,3},
\]
we conclude that
\[
\sqrt{n}(\hat{\gamma}-\gamma)
\xrightarrow{d}
N(0,\sigma^2),
\qquad
\sigma^2:=V_4+W_4+\mathcal V.
\]
Substituting the expressions for $V_4$, $W_4$, and $\mathcal V$
gives the variance formula stated in
\eqref{lim_var_defn_explicit}.

\end{proof}

\subsection{Proofs for the Cross-fitted estimator \label{sec:proof_cf}}

In this section we discuss the proof of the cross-fitted estimator i.e. Theorem \ref{thm_cf}, which uses the same $\Delta$
decomposition as the plug-in proof, but the fold-swapped estimators
reuse the same observations. Their joint behavior therefore requires
an analysis of the projection structure. We use the conditional
independence of $P_{B_s}E_s$ and $P_{B_s}^{\perp}E_s$ given $B_s$
for each $s\in\{a,b\}$, which follows from Gaussianity and
orthogonality. We decompose the centered and scaled estimator into
$T_1$ and $T_2$ and combine them by sequential conditioning.

\subsubsection{Decomposition of the cross fitted estimator}

 Let us recall the definition of $\Delta$s.

\begin{align*}\Delta_1&=\sqrt{n}\cin (\ba^T+e^T_1(B_b^TB_b)^{-1}B_b^TE_b)H_{\lambda_1,a}^{-1}E_a^TP^{\perp}_{B_a}D_a\\
    \Delta_2&=\sqrt{n}\cin e^T_1(B_b^TB_b)^{-1}B_b^TE_b(\tm-\lambda_1H_{\lambda_1,a}^{-1}\tm).\\
\Delta_3&=\sqrt{n}(\cin \ba^T(I-\lambda_1H_{\lambda_1,a}^{-1})\tm-\ba^T\tm)
\end{align*}
Further let,
\begin{align*}T_1:=\frac{1}{2}\sum_{\vartheta}\Delta_{1,\vartheta},\quad T_2:=\frac{1}{2}\sum_{\vartheta}\left(\Delta_{2,\vartheta}+\Delta_{3,\vartheta}\right)
\end{align*}
so that
\begin{equation}\sqrt{n}(\hat{\gamma}_{cf}-\gamma)=T_1+T_2. \label{cf_decomp}\end{equation} 

We will discuss the convergence of $T_1$ and $T_2$ separately. Finally, combining the analyses yields the CLT for the cross-fitted estimator, as in the plug-in case.

\subsubsection{Conditional distribution of $T_1$}

The analysis of $T_1$ proceeds by observing that since $D_a$ and $D_b$ are independent, conditioned on $\mathcal F_{2,n}$, the two terms $\Delta_{1,\vartheta_0}$ and $\Delta_{1,\vartheta_1}$ are independent. The limit of the variance proceeds as in the Lemmas \ref{V1_Lemma}, \ref{V2_lemma}, and $\ref{V3_Lemma}$. They are collected into the single lemma stated below, whose proof is deferred to the appendix.
\begin{lemma}\label{T1_conv}
Recall that $\mathcal F_{2,n}:=\sigma(B_a,E_a,B_b,E_b).$
Then, we have,
\[
T_1\mid\mathcal F_{2,n}
\sim
N(0,V_{T_1,n}),
\]
where
\[
V_{T_1,n}
:=
V_{T_1}(B_a,E_a,B_b,E_b)
\xrightarrow{P}\frac{V_4}{2}.
\]
\end{lemma}

\subsubsection{Asymptotic normality of $T_2$}

Write $\mathcal X_s=(B_s,E_s)$ for $s\in\{a,b\}$.
Then $T_2$ has the representation
\begin{align*}
T_2=\frac12[g(\mathcal{X}_a)+g(\mathcal{X}_b)+h^T(\mathcal{X}_b)l(\mathcal{X}_a)+h^T(\mathcal{X}_a)l(\mathcal{X}_b)]
\end{align*}
where 
\begin{align*}
    g(\mathcal{X}_a)&=\sqrt{n}(\cin \ba^T(I-\lambda_1H_{\lambda_1,a}^{-1})\tm-\ba^T\tm)\\
    l(\mathcal{X}_a)&=(\tm-\lambda_1H_{\lambda_1,a}^{-1}\tm)\\
h^T(\mathcal{X}_b)&=\sqrt{n}\cin  e_1^T(B^T_bB_b)^{-1}B_b^TE_b
\end{align*}
Recall that $\bea=P^{\perp}_{B_a} E_a$ and define $\breve{E}_a:=P_{B_a}E_a$. We note that $h^T(\mathcal{X}_a)$ can be rewritten as 
$$h^T(\mathcal{X}_a)=\sqrt{n}\cin  e_1^T(B^T_aB_a)^{-1}B_a^T\breve{E}_a$$
and since $\bea$ and $\breve{E}_a$ are independent conditioned on $B_a$, $h^T(\mathcal{X}_a)$ is independent of $\bea$ conditioned on $B_a$. Thus we observe that 

\begin{align*}
    T_2&=S_1+S_2\\
    S_1&=\frac12[ h^T(\mathcal{X}_b)l(\mathcal{X}_a)+h^T(\mathcal{X}_a)l(\mathcal{X}_b)]\\
    S_2&=\frac12[g(\mathcal{X}_a)+g(\mathcal{X}_b)]
\end{align*}

Using the conditional independence of $\bar E_s$ and
$\breve E_s$ given $B_s$, for $s\in\{a,b\}$, together with
independence of the two folds, we prove the following lemma:

\begin{lemma}\label{S_1_conv}
Let $\mathcal G_n
:=
\sigma(B_a,B_b,\bar E_a,\bar E_b).
$
Then $S_2$ is $\mathcal G_n$-measurable and
\[
S_1\mid\mathcal G_n
\sim
N(0,V_{S_1,n}),
\]
where
\[
V_{S_1,n}
:=
V_{S_1}(\bar E_a,\bar E_b,B_a,B_b)
\xrightarrow{P}\frac{W_4}{2}.
\]
\end{lemma}

Next we analyze $S_2$. 
\begin{align*}
    S_2=\frac12[\sqrt{n}(\cin  \ba^T(I-\lambda_1H_{\lambda_1,a}^{-1})\tm-\ba^T\tm)+\sqrt{n}(\cin  \ba^T(I-\lambda_1H_{\lambda_1,b}^{-1})\tm-\ba^T\tm)]
\end{align*}

Using the CLT for bilinear form similar to Lemma \ref{delta3_conv_lemma}, we have the following lemma:

\begin{lemma}{\label{S_2_conv_lemma}}

\begin{equation}S_2\xrightarrow{d} N(0,\mathcal{V}/2)\label{S_2_conv}
\end{equation}
where $\mathcal{V}$ is as given in \eqref{fluc_var}.

\end{lemma}

Since $S_2$ is $\mathcal G_n$-measurable,
Lemma~\ref{S_1_conv} gives
\[
S_1\mid\mathcal G_n
\sim N(0,V_{S_1,n}),
\qquad
V_{S_1,n}\xrightarrow{P}\frac{W_4}{2}.
\]
Moreover, by Lemma~\ref{S_2_conv_lemma},
$
S_2\xrightarrow{d}
N\left(0,\mathcal V/2\right)$. 
Therefore, applying
Lemma~\ref{lem:cond_gauss_shift_sum} with
$X_n=S_1$ and $Y_n=S_2$ yields
\begin{equation}
T_2=S_1+S_2
\xrightarrow{d}
N\left(0,\frac12(W_4+\mathcal V)\right).
\label{T_2_mid_conv}\end{equation}

Finally, we combine the above lemmas, to characterize the asymptotic distribution of the cross-fitted estimator.

\subsubsection{Asymptotic distribution of $\hat{\gamma}_{cf}$}

\begin{proof}[Proof of Theorem~\ref{thm_cf}]

Recall from \eqref{cf_decomp} that
\[
\sqrt{n}(\hat{\gamma}_{cf}-\gamma)=T_1+T_2.
\]

By \eqref{T_2_mid_conv}, we have
\[
T_2\xrightarrow{d}
N\left(0,\frac12(W_4+\mathcal V)\right).
\]
Now let $
\mathcal F_{2,n}:=\sigma(B_a,E_a,B_b,E_b).$
The random variable $T_2$ is $\mathcal F_{2,n}$-measurable.
Moreover, by Lemma~\ref{T1_conv},
\[
T_1\mid\mathcal F_{2,n}
\sim
N(0,V_{T_1,n}),
\qquad
V_{T_1,n}\xrightarrow{P}\frac{V_4}{2}.
\]
Therefore another application of
Lemma~\ref{lem:cond_gauss_shift_sum}, with
$X_n=T_1$ and $Y_n=T_2$, gives
\[
T_1+T_2
\xrightarrow{d}
N\left(
0,
\frac12(V_4+W_4+\mathcal V)
\right).
\]
Since
$\sqrt{n}(\hat{\gamma}_{cf}-\gamma)=T_1+T_2,$
we obtain
\[
\sqrt{n}(\hat{\gamma}_{cf}-\gamma)
\xrightarrow{d}
N(0,\sigma_{cf}^2),
\qquad
\sigma_{cf}^2
:=
\frac12(V_4+W_4+\mathcal V).
\]

Since the plug-in variance satisfies
$
\sigma^2=V_4+W_4+\mathcal V,
$
it follows that
$
\sigma_{cf}^2=\frac{\sigma^2}{2}.
$
Substituting the expressions for $V_4$, $W_4$, and $\mathcal V$
gives the stated variance formula.
\end{proof}

 \subsection{Proofs for the estimation of limiting variance \label{sec:var_est}}
 In this section we prove Theorem~\ref{fluc_est_var} by consistently
estimating the variance factor $\sigma^2$ in
\eqref{lim_var_defn}. The limiting variances in
Theorems~\ref{thm_plugin} and~\ref{thm_cf} are then estimated by
$2\hat\sigma^2$ and $\hat\sigma^2$, respectively. Recall that

\begin{align*}
\sigma^2=&(1-\kappa_1f_1(\kappa_1,c_1))^{-2}\left[ f_2(\kappa_1,c_1)(c_1 f^{\Sigma_B}_5(c_2)+\alpha \mathfrak{b}_A)\right.\\
&+(1-2\kappa_1f_1(\kappa_1,c_1)-\kappa^2_1\partial_{\kappa}f_1(\kappa_1,c_1))f^{\Sigma_B}_5(c_2) \mathfrak{t}_M\\
&\left.+\kappa^2_1\frac{(f_1(\kappa_1,c_1)\kappa_1)^{4}}{c_1(1-c_2)}\left(h_1(-\kappa_1,c_1)\mathfrak{g}^2+h_1(-\kappa_1,c_1)\mathfrak{b}_A\mathfrak{t}_M\right)\right]
\end{align*}

Hence if we can estimate $\mathfrak{t}_M, \mathfrak{b}_A, \mathfrak{g},$ and $f^{\Sigma_B}_5(c_2)$, we can plugin the estimates in the above expression and obtain $\hat{\sigma}^2$.

\subsubsection{Estimator for $\|\tm\|^2$}
We recall the ridge estimator for $\tm$ given in \eqref{ridge_decomp}

$$\breve{\theta}_M=\tm-\lambda_1H_{\lambda_1,a}^{-1}\tm+ H_{\lambda_1,a}^{-1}E_a^TP^{\perp}_{B_a}D_a.$$
 The asymptotic behaviour of $\|\breve{\theta}_M\|^2$ can be characterized as follows:

 \begin{lemma}{\label{tm2_lemma}}
  $$\|\btm\|^2=(1-2\kappa_1f_1(\kappa_1,c_1)-\kappa^2_1\partial_{\kappa}f_1(\kappa_1,c_1))\mathfrak{t}_M+c_1f_2(\kappa_1,c_1)+o_p(1)$$
  
 \end{lemma}

Hence defining

$$\hat{\mathfrak{t}}_M=(1-2\kappa_1f_1(\kappa_1,c_1)-\kappa^2_1\partial_{\kappa}f_1(\kappa_1,c_1))^{-1}(\|\btm\|^2-c_1f_2(\kappa_1,c_1)),$$ it follows that $\hat{\mathfrak{t}}_M$ is a consistent estimator for $\mathfrak{t}_M$.

\subsubsection{Estimator for $f^{\Sigma_B}_5(c_2) $}

Since $e^T_1(B_b^TB_b/n)^{-1}e_1$ converges to $f^{\Sigma_B}_5(c_2)$, we can estimate $f^{\Sigma_B}_5(c_2)$ by 

$$\hat{\mathfrak{r}}=e^T_1(B_b^TB_b/n)^{-1}e_1.$$

\subsubsection{Estimator for $\|\ba\|^2$}
 Recall OLS estimator for $\ba$

\begin{align*}\hba^T=&\ba^T+e^T_1(B_b^TB_b)^{-1}B_b^TE_b.\\
\end{align*}
We have the following lemma:

\begin{lemma}{\label{ba2_lemma}}
$$\|\hba\|^2=\mathfrak{b}_A+c_1(1-c_2)f^{\Sigma_B}_5(c_2) +o_p(1)$$
\end{lemma}

Hence defining

$$\hat{\mathfrak{b}}_A=\|\hba\|^2-c_1(1-c_2)\hat{\mathfrak{r}} ,$$ it follows that $\hat{\mathfrak{b}}_A$ is a consistent estimator for $\mathfrak{b}_A$.

\subsubsection{Estimator for $\gamma^2$}

Since we already have a consistent estimator $\hat{\gamma}$ for $\gamma$, we can just use $\hat{\gamma}^2$ as an estimator for $\gamma^2$ and hence for $\mathfrak{g}^2$.

\begin{proof}[Proof of Theorem \ref{fluc_est_var}]

Using the consistent estimators for  $\mathfrak{t}_M, \mathfrak{b}_A, \mathfrak{g},$ and $f^{\Sigma_B}_5(c_2)$ we observe that

\begin{align*}
\hat{\sigma}^2=&(1-\kappa_1f_1(\kappa_1,c_1))^{-2}\left[c_1 \hat{\mathfrak{r}} f_2(\kappa_1,c_1)+\alpha f_2(\kappa_1,c_1)\hat{\mathfrak{b}}_A\right.\\
&+(1-2\kappa_1f_1(\kappa_1,c_1)-\kappa^2_1\partial_{\kappa}f_1(\kappa_1,c_1))\hat{\mathfrak{r}} \hat{\mathfrak{t}}_M\\
&\left.+\kappa^2_1\frac{(f_1(\kappa_1,c_1)\kappa_1)^{4}}{c_1(1-c_2)}\left(h_1(-\kappa_1,c_1)\hat{\gamma}^2+h_1(-\kappa_1,c_1)\hat{\mathfrak{b}}_A\hat{\mathfrak{t}}_M\right)\right]
\end{align*}
is a consistent estimator of $\sigma^2$.

Note that we may assume without loss of generality that
$\hat{\sigma}^2>0$. Indeed, since
\[
\hat{\sigma}^2\xrightarrow{P}\sigma^2
\qquad\text{with}\qquad
\sigma^2>0,
\]
we have $\mathbb{P}(\hat{\sigma}^2\leq 0)\to0$. More formally,  define
$
\hat{\sigma}_+^2:=\max\{\hat{\sigma}^2,n^{-1}\}
$.
Then
\[
\mathbb{P}(\hat{\sigma}_+^2\neq\hat{\sigma}^2)
=
\mathbb{P}(\hat{\sigma}^2<n^{-1})
\rightarrow0.
\]
Hence replacing $\hat{\sigma}^2$ by $\hat{\sigma}_+^2$ in the
studentized statistic does not affect its asymptotic distribution.
    
\end{proof}

\subsection{Proofs for the pathwise mediator effect \label{sec:pathwise_proof_plugin}}

In this section we prove Theorem \ref{thm_plugin_path}. Similar to the proof of Theorem \ref{thm_plugin}, the proof proceeds by decomposing $\sqrt{n}(\hej-\ej)$ into three terms and considering conditional distributions. We will condition on $\{B_a,E_a,B_b,E_b\}$ and  $\{B_a, E_a, B_b\}$ as before. The final step will be to consider the fluctuation as function of $\{E_a\}$ to obtain the desired normal distribution.
\subsubsection{Decomposition of $\hej$}

Recall that
$$\htm=\cin  (\tm-\lambda_1H_{\lambda_1,a}^{-1}\tm+ H_{\lambda_1,a}^{-1}E_a^TP^{\perp}_{B_a}D_a)$$
and 
$$\hba^T=\ba^T+e^T_1(B_b^TB_b)^{-1}B_b^TE_b.$$
We start with a following decomposition, similar to the one used for $\gamma$: 

\begin{equation}\sqrt{n}(\hej-\ej)=\Gamma_1+\Gamma_2+\Gamma_3\label{eta_decomp}\end{equation} with the terms on the RHS of the above equation given by:

\begin{align*}\Gamma_1&=\sqrt{n}\cin  (\ba^T+e^T_1(B_b^TB_b)^{-1}B_b^TE_b)e_je^T_jH_{\lambda_1,a}^{-1}E_a^TP^{\perp}_{B_a}D_a,\\
    \Gamma_2&=\sqrt{n}\cin  e^T_1(B_b^TB_b)^{-1}B_b^TE_be_je_j^T(\tm-\lambda_1H_{\lambda_1,a}^{-1}\tm),\\
\Gamma_3&=\sqrt{n}(\cin  \ba^Te_je_j^T(I-\lambda_1H_{\lambda_1,a}^{-1})\tm-\beta_{A,j}\theta_{M,j}).
\end{align*}

\subsubsection{Conditional distribution of \texorpdfstring{$\Gamma_1$}{Gamma_1}}

 Recall that $\mathcal F_{2,n}=\sigma(B_a,E_a,B_b,E_b).$
$D_a\sim N(0,I_n)$, conditional on $\mathcal F_{2,n}$, it follows that
$$\Gamma_1| \mathcal{F}_{2,n}\sim N(0,\tV_1(B_a,E_a,B_b,E_b))$$
with \begin{equation}\tV_1(B_a,E_a,B_b,E_b)=n\cin^2 e^T_jH_{\lambda_1,a}^{-1}\bea^T\bea H_{\lambda_1,a}^{-1}e_j(\baj +e_j^TE_b^TB_b(B_b^TB_b)^{-1}e_1)^2\label{tV_1_squared_form}.\end{equation}
where we have used $\bea=P^{\perp}_{B_a}E_a$. We recall that the resolvent can be scaled as follows:
\begin{align*}
  H_{\kappa_1,a}^{-1}&:=(n-d-1)H_{\lambda_1,a}^{-1}=(\bea^T\bea/(n-d-1)+\kappa_1I)^{-1}\\
\end{align*}
Next define the following quantities:
\begin{align}
   \tV_3(B_a,E_a)=& \alpha\mathfrak{c}^2_1 e_j^TH_{\kappa_1,a}^{-1}\frac{\bea^T\bea}{n-d-1} H_{\kappa_1,a}^{-1} e_j\baj^2,\nonumber\\
   \tV_4=& \alpha\mathfrak{c}^2_1 f_2(\kappa_1,c_1)\baj^2.\label{tV4_def}
\end{align}
where $\alpha:=\lim_{n\rightarrow \infty}\frac{n}{n-d-1}$. We have the following lemma:
\begin{lemma}{\label{tV1_Lemma}}
\begin{align*}
        \tV_1(B_a,E_a,B_b,E_b)&=\tV_3(B_a,E_a)+o_p(1)\\
        \tV_3(B_a,E_a)&=\tV_4+o_p(1)\\
\end{align*}
\end{lemma}

\subsubsection{Conditional distribution of $\Gamma_2$}

Recall that 
\begin{equation}\Gamma_2=\sqrt{n}\cin  e_1^T(B^T_bB_b)^{-1}B_b^TE_be_je_j^T(\tm-\lambda_1H_{\lambda_1,a}^{-1}\tm)
\end{equation}
and $\mathcal F_{1,n}=\sigma(B_a,E_a,B_b)$.
Since the entries of $E_b$ are
independent standard Gaussian random variables, so  we have,
\begin{equation}\Gamma_2|\mathcal{F}_{1,n}\sim N(0, \tW_2(B_a,E_a,B_b))
\end{equation}

with 
$$\tW_2(B_a,E_a,B_b)=n(e_j^T(I-\lambda_1H_{\lambda_1,a}^{-1})\tm)^2\cin  ^2e_1^T(B_b^TB_b)^{-1}e_1.$$
Now, define the following quantities:
\begin{align*}
    \tW_3(B_b)&=n\mathfrak{c}_1^2(1-\kappa_1f_1(\kappa_1,c_1))^2\tmj^2e_1^T(B_b^TB_b)^{-1}e_1\\
    \tW_4&=\mathfrak{c}_1^2(1-\kappa_1f_1(\kappa_1,c_1))^2f^{\Sigma_B}_5(c_2) \tmj^2.
\end{align*}

\begin{lemma}{\label{tW_approx}}
   The following relations hold:
    \begin{align}
    \tW_2(B_a,E_a,B_b)&=\tW_3(B_b)+o_p(1)\label{tW_23}.\\
        \tW_3(B_b)&=\tW_4+o_p(1).\label{tW_34}
    \end{align}
\end{lemma}

\subsubsection{Asymptotic distribution of $\Gamma_3$}
Recall that 

\begin{align*}\Gamma_3&=\sqrt{n}(\cin  \ba^Te_je^T_j(I-\lambda_1H_{\lambda_1,a}^{-1})\tm-\ba^Te_je_j^T\tm)
\end{align*}

We show the convergence of $\Gamma_3$, by simplifying $\Gamma_3$ to a centered and scaled version of $e_j^TH_{\lambda_1,a}^{-1}\tm$. We note that this is a bilinear form of the  resolvent $H_{\lambda_1,a}^{-1}$. Once again the CLT for bilinear forms of resolvents (\cite{yin2023limiting}), can be used to prove the following lemma:

\begin{lemma}{\label{gamma3_conv_lemma}}

\begin{equation}\Gamma_3\xrightarrow{d} N(0,\ctV),
\label{Gamma_3_conv}\end{equation}
where,
\begin{align}\ctV=&\kappa_1^2\mathfrak{c}_1^2\frac{\baj^2(f_1(\kappa_1,c_1)\kappa_1)^{4}}{c_1(1-c_2)}\left(h_1(-\kappa_1,c_1)\tmj^2+h_1(-\kappa_1,c_1)\mathfrak{t}_M\right)\label{fluc_var_new}
\end{align}

\end{lemma}

\subsubsection{Asymptotic distribution of $\hej$}

\begin{proof}[Proof of Theorem~\ref{thm_plugin_path}]

Recall from \eqref{eta_decomp} that
\[
\sqrt{n}(\hej-\ej)
=
\Gamma_1+\Gamma_2+\Gamma_3.
\]

 We first combine $\Gamma_2$ and $\Gamma_3$.
Since $\Gamma_3$ depends only on $(B_a,E_a)$,
it is $\mathcal F_{1,n}$-measurable. Moreover, conditional
on $\mathcal F_{1,n}$, the only randomness in $\Gamma_2$
is due to $E_b$. Hence
\[
\Gamma_2\mid\mathcal F_{1,n}
\sim
N\!\left(
0,
\tilde W_2(B_a,E_a,B_b)
\right).
\]
By Lemma~\ref{tW_approx},
\[
\tilde W_2(B_a,E_a,B_b)
\xrightarrow{P}
\tilde W_4.
\]
On the other hand, Lemma~\ref{gamma3_conv_lemma} gives
\[
\Gamma_3
\xrightarrow{d}
N(0,\tilde{\mathcal V}).
\]
Therefore, applying
Lemma~\ref{lem:cond_gauss_shift_sum} with
$
X_n=\Gamma_2,
Y_n=\Gamma_3,
$
we obtain
\[
\Gamma_{2,3}
:=
\Gamma_2+\Gamma_3
\xrightarrow{d}
N(0,\tilde W_4+\tilde{\mathcal V}).
\]

Next, $\Gamma_{2,3}$ is $\mathcal F_{2,n}$-measurable.
Conditional on $\mathcal F_{2,n}$, the only randomness in
$\Gamma_1$ is due to $D_a$. Therefore
\[
\Gamma_1\mid\mathcal F_{2,n}
\sim
N\!\left(
0,
\tilde V_1(B_a,E_a,B_b,E_b)
\right).
\]
By Lemma~\ref{tV1_Lemma},
\[
\tilde V_1(B_a,E_a,B_b,E_b)
\xrightarrow{P}
\tilde V_4.
\]
A second application of
Lemma~\ref{lem:cond_gauss_shift_sum}, now with
$
X_n=\Gamma_1,
Y_n=\Gamma_{2,3},
$
gives
\[
\Gamma_1+\Gamma_{2,3}
\xrightarrow{d}
N\!\left(
0,
\tilde V_4+\tilde W_4+\tilde{\mathcal V}
\right).
\]
Consequently,
\[
\sqrt{n}(\hej-\ej)
\xrightarrow{d}
N(0,\sigma_j^2),
\]
where
\[
\sigma_j^2
:=
\tilde V_4+\tilde W_4+\tilde{\mathcal V}.
\]
Substituting the expressions for
$\tilde V_4$, $\tilde W_4$, and
$\tilde{\mathcal V}$ yields the stated variance formula.
\end{proof}

\subsection{Proofs for the Cross-fitted estimator for pathwise mediator effect \label{sec:pathwise_proof_cf}}

The proof of Theorem~\ref{thm_cf_path} follows the same
conditioning scheme as the global cross-fitted proof.
Lemmas~\ref{tT1_conv}, \ref{tS_1_conv}, and
\ref{tS_2_conv_lemma} are the corresponding coordinatewise
versions. We state these components and explain how they combine;
their proofs follow from the same projection and conditional-variance
arguments. Recall the definitions of the $\Gamma$ terms.

\begin{align*}\Gamma_1&=\sqrt{n}\cin  (\ba^T+e^T_1(B_b^TB_b)^{-1}B_b^TE_b)e_je^T_jH_{\lambda_1,a}^{-1}E_a^TP^{\perp}_{B_a}D_a\\
    \Gamma_2&=\sqrt{n}\cin  e^T_1(B_b^TB_b)^{-1}B_b^TE_be_je^T_j(\tm-\lambda_1H_{\lambda_1,a}^{-1}\tm).\\
\Gamma_3&=\sqrt{n}(\cin  \ba^Te_je^T_j(I-\lambda_1H_{\lambda_1,a}^{-1})\tm-\beta_{A,j}\theta_{M,j})
\end{align*}
Further let,
\begin{align*}\tT_1&:=\frac{1}{2}\sum_{\vartheta}\Gamma_{1,\vartheta}\\
\tT_2&:=\frac{1}{2}\sum_{\vartheta}\left(\Gamma_{2,\vartheta}+\Gamma_{3,\vartheta}\right)
\end{align*}
so that
\begin{equation}
\sqrt{n}(\hat{\eta}_{j,cf}-\eta_j)
=
\tT_1+\tT_2.
\label{eta_cf_decomp}
\end{equation}
The conditional distribution of $\tT_1$ follows in the same way as the conditional
distribution of $T_1$ in the proof of Theorem \ref{thm_cf}. 
\begin{lemma}{\label{tT1_conv}}
$$\tT_1|\mathcal{F}_{2,n}\sim N(0, V_{\tT_1}(B_a,E_a,B_b,E_b))$$
with 
$$V_{\tT_1}(B_a,E_a,B_b,E_b)=\frac{1}{2}\tV_4+o_p(1)$$
where $\tV_4$ is defined as in \eqref{tV4_def}.
\end{lemma}

Next we note that $\tT_2$ is of the form:
\begin{align*}
\tT_2=\frac12[g(\mathcal{X}_a)+g(\mathcal{X}_b)+h(\mathcal{X}_b)l(\mathcal{X}_a)+h(\mathcal{X}_a)l(\mathcal{X}_b)]
\end{align*}
where 
\begin{align*}
    g(\mathcal{X}_a)&=\sqrt{n}(\cin  \ba^Te_je_j^T(I-\lambda_1H_{\lambda_1,a}^{-1})\tm-\ba^Te_je_j^T\tm)\\
    l(\mathcal{X}_a)&=e_j^T(\tm-\lambda_1H_{\lambda_1,a}^{-1}\tm)\\
h(\mathcal{X}_b)&=\sqrt{n}\cin  e_1^T(B^T_bB_b)^{-1}B_b^TE_be_j
\end{align*}
Recall that $\bea=P^{\perp}_{B_a} E_a$ and define $\breve{E}_a:=P_{B_a}E_a$. We note that $h(\mathcal{X}_a)$ can be rewritten as 
$$h(\mathcal{X}_a)=\sqrt{n}\cin  e_1^T(B^T_aB_a)^{-1}B_a^T\breve{E}_ae_j$$
and since $\bea$ and $\breve{E}_a$ are independent conditioned on $B_a$, $h(\mathcal{X}_a)$ is independent of $\bea$ conditioned on $B_a$. Thus we observe that 

\begin{align*}
    \tT_2&=\tS_1+\tS_2\\
    \tS_1&=\frac12[ h(\mathcal{X}_b)l(\mathcal{X}_a)+h(\mathcal{X}_a)l(\mathcal{X}_b)]\\
    \tS_2&=\frac12[g(\mathcal{X}_a)+g(\mathcal{X}_b)]
\end{align*}

\begin{lemma}\label{tS_1_conv}
Recall that
$
\mathcal G_n
:=
\sigma(B_a,B_b,\bar E_a,\bar E_b)$.
Then $\tS_2$ is $\mathcal G_n$-measurable and
$
\tS_1\mid\mathcal G_n
\sim
N(0,V_{\tS_1,n}),
$
where
\[
V_{\tS_1,n}
:=
V_{\tS_1}(\bar E_a,\bar E_b,B_a,B_b)
\xrightarrow{P}
\frac{\tW_4}{2}.
\]
\end{lemma}

Next we analyze $\tS_2$, which has the following expansion:
\begin{align*}
    \tilde{S}_2=&\frac12\left[\sqrt{n}(\cin  \ba^Te_je_j^T(I-\lambda_1H_{\lambda_1,a}^{-1})\tm-\ba^Te_je_j^T\tm)\right.\\
    &\left.+\sqrt{n}(\cin  \ba^Te_je_j^T(I-\lambda_1H_{\lambda_1,b}^{-1})\tm-\ba^Te_je_j^T\tm)\right]
\end{align*}

We note the appearance of the bilinear forms of the resolvent in the above expansion and using the CLT, we prove the following lemma:

\begin{lemma}{\label{tS_2_conv_lemma}}

\begin{equation}\tilde S_2\xrightarrow{d} N(0,\ctV/2)\label{tS_2_conv}
\end{equation}
where $\ctV$ is as given in \eqref{fluc_var_new}.

\end{lemma}

Since $\tS_2$ is $\mathcal G_n$-measurable,
Lemma~\ref{tS_1_conv} gives
\[
\tS_1\mid\mathcal G_n
\sim
N(0,V_{\tS_1,n}),
\qquad
V_{\tS_1,n}\xrightarrow{P}\frac{\tW_4}{2}.
\]
On the other hand, by Lemma~\ref{tS_2_conv_lemma},
$
\tS_2\xrightarrow{d}
N\left(0,\ctV/2\right).
$
Therefore, applying
Lemma~\ref{lem:cond_gauss_shift_sum} with
$X_n=\tS_1$ and $Y_n=\tS_2$, we obtain
\[
\tT_2=\tS_1+\tS_2
\xrightarrow{d}
N\left(0,\frac{\tW_4+\ctV}{2}\right).
\]

Finally we are ready to show the asymptotic distribution of $\hat{\eta}_{j,cf}$

\begin{proof}[Proof of Theorem~\ref{thm_cf_path}]
Recall from \eqref{eta_cf_decomp} that
\[
\sqrt{n}(\hat{\eta}_{j,cf}-\eta_j)
=
\tT_1+\tT_2.
\]

As shown above,
\[
\tT_2
\xrightarrow{d}
N\left(0,\frac{\tW_4+\ctV}{2}\right).
\]

Since $\tT_2$ is a function only of
$(B_a,E_a,B_b,E_b)$, it is $\mathcal F_{2,n}$-measurable.
Moreover, by Lemma~\ref{tT1_conv},
\[
\tT_1\mid\mathcal F_{2,n}
\sim
N(0,V_{\tT_1,n}),
\qquad
V_{\tT_1,n}=V_{\tT_1}(B_a,E_a,B_b,E_b)
\xrightarrow{P}
\frac{\tV_4}{2}.
\]
Therefore, applying
Lemma~\ref{lem:cond_gauss_shift_sum} with
$X_n=\tT_1$ and $Y_n=\tT_2$, we obtain
\[
\tT_1+\tT_2
\xrightarrow{d}
N\left(
0,
\frac{\tV_4+\tW_4+\ctV}{2}
\right).
\]
Hence
\[
\sqrt{n}(\hat{\eta}_{j,cf}-\eta_j)
\xrightarrow{d}
N(0,\sigma_{j,cf}^2),
\]
where
\[
\sigma_{j,cf}^2
:=
\frac12(\tV_4+\tW_4+\ctV).
\]
Since
$
\sigma_j^2
=
\tV_4+\tW_4+\ctV,
$
it follows immediately that
$
\sigma_{j,cf}^2=\frac{\sigma_j^2}{2}.
$
\end{proof}

\subsection{Proofs for the limiting variance of the pathwise mediator effect \label{sec:var_est_path}}

In this section we prove Theorem~\ref{fluc_est_var_path} by
consistently estimating the variance factor $\sigma_j^2$ in
\eqref{lim_var_sigma_j}. The limiting variances in
Theorems~\ref{thm_plugin_path} and~\ref{thm_cf_path} are then
estimated by $2\hat\sigma_j^2$ and $\hat\sigma_j^2$,
respectively. Recall that

\begin{align*}
\sigma^2_j
=& (1-\kappa_1f_1(\kappa_1,c_1))^{-2}\Big(\alpha f_2(\kappa_1,c_1)\baj^2+(1-\kappa_1f_1(\kappa_1,c_1))^2f^{\Sigma_B}_5(c_2) \tmj^2\Big.\\
&\Big.+\kappa_1^2\frac{\baj^2(f_1(\kappa_1,c_1)\kappa_1)^{4}}{c_1(1-c_2)}\left(h_1(-\kappa_1,c_1)\tmj^2+h_1(-\kappa_1,c_1)\mathfrak{t}_M\right)\Big).\\
\end{align*}
Since we already have estimators for $\mathfrak{t}_M$ and $f^{\Sigma_B}_5(c_2)$ from Theorem \ref{fluc_est_var}, only estimators for $\tmj$ and $\baj$ are required.

\subsubsection{Estimator for $\baj$}
 Recall OLS estimator for $\ba$

\begin{align*}\hba^T=&\ba^T+e^T_1(B_b^TB_b)^{-1}B_b^TE_b.\\
\end{align*}
We will show that the $j^{th}$ element of $\hba$ i.e. $\hat{\beta}_{A,j}$ is a consistent estimator for $\beta_{A,j}$.

Indeed we have that

$$\hat{\beta}_{A,j}=\beta_{A,j}+e^T_1(B_b^TB_b)^{-1}B_b^TE_be_j$$

Since the entries of $E_b$ are i.i.d. standard normal, conditioned on $B_b$, $e^T_1(B_b^TB_b)^{-1}B_b^TE_be_j$ has mean $0$ and variance $e_1^T(B_b^TB_b)^{-1}e_1$. Since $e_1^T\left(\frac{B_b^TB_b}{n}\right)^{-1}e_1=O_p(1)$, it follows that $e^T_1(B_b^TB_b)^{-1}B_b^TE_be_j=o_p(1)$ and hence $\hat{\beta}_{A,j}$ is a consistent estimator of $\baj$.

\subsubsection{Estimator for $\tmj$}

Recall that
$$\htm=\cin  (\tm-\lambda_1H_{\lambda_1,a}^{-1}\tm+ H_{\lambda_1,a}^{-1}E_a^TP^{\perp}_{B_a}D_a)$$
so that

$$\hat{\theta}_{M,j}=\cin  e_j^T(I-\lambda_1H_{\lambda_1,a}^{-1})\tm+ \cin  e^T_jH_{\lambda_1,a}^{-1}E_a^TP^{\perp}_{B_a}D_a$$

Using the deterministic equivalent relation i.e. Theorem \ref{mp_detequ} we have 
$$e_j^T(I-\lambda_1H_{\lambda_1,a}^{-1})\tm=(1-\kappa_1 f_1(\kappa_1,c_1))\tmj+o_p(1),$$ so that using the fact that $\mathfrak{c}_1=(1-\kappa_1 f_1(\kappa_1,c_1))^{-1}$ and $\cin \rightarrow \mathfrak{c}_1$ we obtain
$$\cin e_j^T(I-\lambda_1H_{\lambda_1,a}^{-1})\tm\xrightarrow{P}\tmj.$$ On the other hand

\begin{align*}\mathrm{Var}(e^T_jH_{\lambda_1,a}^{-1}E_a^TP^{\perp}_{B_a}D_a|E_a, B_a)&=e_j^TH_{\lambda_1,a}^{-1}\bea^T\bea H_{\lambda_1,a}^{-1}e_j\\
&=\frac{1}{n-d-1}e_j^TH_{\kappa_1,a}^{-1}\frac{\bea^T\bea}{n-d-1} H_{\kappa_1,a}^{-1}e_j\end{align*}

Since $e_j^TH_{\kappa_1,a}^{-1}\frac{\bea^T\bea}{n-d-1} H_{\kappa_1,a}^{-1}e_j=O_p(1)$, we have $\mathfrak{c}_1e^T_jH_{\lambda_1,a}^{-1}E_a^TP^{\perp}_{B_a}D_a=o_p(1)$, and hence $\hat{\theta}_{M,j}$ is a consistent estimator of $\tmj$.

\begin{proof}[Proof of Theorem \ref{fluc_est_var_path}]

Using the consistent estimators for  $\mathfrak{t}_M, \beta_{A,j}, \theta_{M,j},$ and $f^{\Sigma_B}_5(c_2)$ we observe that

\begin{align*}
\hat{\sigma}^2_j
=& (1-\kappa_1f_1(\kappa_1,c_1))^{-2}\Big(\alpha f_2(\kappa_1,c_1)\hat{\beta}_{A,j}^2+(1-\kappa_1f_1(\kappa_1,c_1))^2\hat{\mathfrak{r}} \hat{\theta}_{M,j}^2\Big.\\
&\Big.+\kappa_1^2\frac{\hat{\beta}_{A,j}^2(f_1(\kappa_1,c_1)\kappa_1)^{4}}{c_1(1-c_2)}\left(h_1(-\kappa_1,c_1) \hat{\theta}_{M,j}^2+h_1(-\kappa_1,c_1)\hat{\mathfrak{t}}_M\right)\Big).\\
\end{align*}
is a consistent estimator of $\sigma^2_j$.
Finally we observe that analogous to Theorem \ref{fluc_est_var}, we can assume w.l.o.g. that $\hat{\sigma}^2_j>0$.
\end{proof}

\section{Proofs of the technical Lemmas \label{sec:tech_lemma}}

\begin{proof}[Proof of Lemma \ref{partial_ridge_lemma}]
 Since 
$$Y_a=A_a\theta_A+M_a\theta_M+C_a\theta_C+D_a, \quad M_a=B_a\varrho^T+E_a$$  the estimator $\btm$ can be rewritten as:
\begin{align}\btm&=(E_a^TP^{\perp}_{B_a} E_a+\lambda_1 I)^{-1}(E_a^TP^{\perp}_{B_a}E_a \tm+ E_a^TP^{\perp}_{B_a}D_a)\nonumber\\
 & =\tm-\lambda_1 (E_a^TP^{\perp}_{B_a} E_a+\lambda_1 I)^{-1}\tm+ (E_a^TP^{\perp}_{B_a} E_a+\lambda_1 I)^{-1}E_a^TP^{\perp}_{B_a}D_a,\label{ridge_decomp}
 \end{align}
where the first equality follows from
$P^{\perp}_{B_a}B_a=0.$ To simplify further we note that $E_a^TP^{\perp}_{B_a}E_a$ can be written as $\bea^T\bea$ with $\bea=P^{\perp}_{B_a}E_a$. We emphasize that $\bea$ is an $n\times p$ projected Gaussian matrix and does not itself have independent rows. However, its Gram matrix admits an i.i.d.\ Gaussian representation. Further conditioned on $B_a$, let $r_a=n-\operatorname{rank}(B_a)$ and let $U_a\in\mathbb R^{n\times r_a}$ have orthonormal columns spanning the range of $P^{\perp}_{B_a}$, so that
$$U_aU_a^T=P^{\perp}_{B_a},\qquad U_a^TU_a=I_{r_a}.$$
Define
$\widetilde{E}_a:=U_a^TE_a.$
Then, conditioned on $B_a$, $\widetilde{E}_a$ is an $(n-\operatorname{rank}(B_a))\times p$ matrix with i.i.d.\ standard Gaussian entries, by rotation invariance of the Gaussian and the independence of $E_a$ and $B_a$. Moreover,
$$\bea=P^{\perp}_{B_a}E_a=U_aU_a^TE_a=U_a\widetilde{E}_a,$$
and therefore
$\widetilde{E}^T_a\widetilde{E}_a=\bea^T\bea.$
Thus, although $\bea$ itself does not have i.i.d.\ rows, its Gram matrix has, conditionally on $B_a$, exactly the same Wishart representation as the Gram matrix of an $(n-\operatorname{rank}(B_a))\times p$ matrix with i.i.d.\ standard Gaussian entries.

Let $H_{\lambda_1,a}=\bea^T\bea+\lambda_1I$ so that
$$\breve{\theta}_M=\tm-\lambda_1H_{\lambda_1,a}^{-1}\tm+ H_{\lambda_1,a}^{-1}E_a^TP^{\perp}_{B_a}D_a.$$

\end{proof}

\begin{proof}[Proof of Lemma \ref{pathwise_med_lemma}]
Recall the LSEM

\begin{align}\stY&=\theta_M^T\stM+\theta_A\stA+\theta_C^T \stC+\delta\label{y_mod_vec_new}\\
\stM&=\beta_A \stA+\beta_C \stC+\epsilon\label{M_mod_vec_new}
\end{align}

Under the intervention $\mathrm{do}(\stA=a)$, the structural equation for $\stA$ is replaced
by the constant $\stA := a$, and the mediator becomes
\[
\stM = \beta_A a + \beta_C \stC + \epsilon
\quad \implies \quad
\mathbb{E}[\stM \mid \mathrm{do}(\stA=a)] = \beta_A a.
\]
Substituting in the equation for $\stY$,
\[
\stY = \theta_A a + \theta_C^\top \stC + \tm^\top \stM + \delta,
\]
and using $\mathbb{E}[\stC]=0$ and $\mathbb{E}[\delta]=0$, we obtain
\[
\mathbb{E}[\stY \mid \mathrm{do}(\stA=a)]
= \theta_A a + \theta_M^\top (\beta_A a)
= (\theta_A + \theta_M^\top \beta_A)\, a
\]

so that 

\begin{equation}\frac{\partial}{\partial a}\,
\mathbb{E}[\stY \mid \mathrm{do}(\stA=a)]=\tm^T\ba+\theta_A \label{deriv_1}
\end{equation}

On the other hand when we compute expectation under $\mathrm{do}(\stM^j = m_j,\stA = a)$, the intervention replaces the structural equations for $\stA$ and $\stM^j$ with
\[
\stA := a, \qquad \stM^j := m_j.
\]
For $k \neq j$, the mediator equation remains unchanged:
\[
\stM^k = \beta_{A,k} \stA + \beta_{C,k}^\top \stC + \epsilon_{k},
\]
so, using $\mathbb{E}[\stC]=0$ and $\mathbb{E}[\epsilon]=0$,
\[
\mathbb{E}[\stM^k \mid \mathrm{do}(\stM^j=m_j, \stA=a)]
= \beta_{A,k} a,
\]
while for the intervened component,
\[
\mathbb{E}[\stM^j \mid \mathrm{do}(\stM^j=m_j, \stA=a)] = m_j.
\]
Collecting these in vector form,
\[
\mathbb{E}[\stM \mid \mathrm{do}(\stM^j=m_j, \stA=a)]
= \beta_A a + (m_j - \beta_{A,j} a)\, e_j,
\]
where $e_j$ is the $j$th standard basis vector. Substituting into the structural equation for $Y$ and taking expectations,
\[
\mathbb{E}[\stY \mid \mathrm{do}(\stM^j=m_j, \stA=a)]
= \theta_A a + \theta_M^\top \mathbb{E}[\stM \mid \mathrm{do}(\stM^j=m_j, \stA=a)].
\]
Therefore,
\[
\mathbb{E}[\stY \mid \mathrm{do}(\stM^j=m_j, \stA=a)]
= (\theta_A + \theta_M^\top \beta_A)\, a
  + \theta_{M,j}\bigl(m_j - \beta_{A,j} a\bigr).
\]

This implies that 
 
 \begin{equation}\frac{\partial}{\partial a}\,
\mathbb{E}\!\left[\stY \,\middle|\, \doo(\stM^j=m_j;\stA=a)\right]=\theta_A+\sum_{k\neq j}\theta_{M,k}\beta_{A,k}\label{deriv_2}\end{equation}

From \eqref{deriv_1} and \eqref{deriv_2}, we obtain $\eta_j=\theta_{M,j}\beta_{A,j}$. We note that this quantity is nothing but $\ba^T P_j \tm$ where $P_j=e_je_j^T$. Thus the natural  estimator for $\eta_j$ is $\hat{\eta}_j=\hba^TP_j\htm$.

\end{proof}

\begin{proof}[Proof of Lemma \ref{V1_Lemma}]

Expanding \eqref{V_1_squared_form}, we obtain
  \begin{align}
    &V_1(B_a,E_a,B_b,E_b)\\
    = &n\cin ^2 \ba^TH_{\lambda_1,a}^{-1}\bea^T\bea H_{\lambda_1,a}^{-1} \ba+2n\cin ^2\ba^TH_{\lambda_1,a}^{-1}\bea^T\bea H_{\lambda_1,a}^{-1}E_b^TB_b(B^T_bB_b)^{-1}e_1\nonumber\\
    & +n\cin ^2 e^T_1(B^T_bB_b)^{-1}B_b^TE_b H_{\lambda_1,a}^{-1}\bea^T\bea H_{\lambda_1,a}^{-1}E_b^TB_b(B_b^TB_b)^{-1}e_1\nonumber\\
    = 
    &\alpha\mathfrak{c}^2_1 \ba^TH_{\kappa_1,a}^{-1}\frac{\bea^T\bea}{n-d-1} H_{\kappa_1,a}^{-1} \ba+2\alpha\mathfrak{c}^2_1\ba^TH_{\kappa_1,a}^{-1}\frac{\bea^T\bea}{n-d-1} H_{\kappa_1,a}^{-1}\frac{E_b^TB_b}{n}\left(\frac{B_b^TB_b}{n}\right)^{-1}e_1\nonumber\\
    & +\alpha\mathfrak{c}^2_1 e^T_1\left(\frac{B_b^TB_b}{n}\right)^{-1}\frac{B_b^TE_b}{n}H_{\kappa_1,a}^{-1}\frac{\bea^T\bea}{n-d-1} H_{\kappa_1,a}^{-1}\frac{E_b^TB_b}{n}\left(\frac{B_b^TB_b}{n}\right)^{-1}e_1\nonumber+o_p(1)\label{V_1_decomp}\\
\end{align}
where $\alpha:=\lim_{n\rightarrow \infty}\frac{n}{n-d-1}$ and we have used the fact that $\cin \rightarrow \mathfrak{c}_1$. Since conditioned on $B_b$, $\frac{E_b^TB_b}{n}\left(\frac{B_b^TB_b}{n}\right)^{-1}e_1\sim n^{-1/2}N(0, e^T_1\left(\frac{B_b^TB_b}{n}\right)^{-1}e_1 I_p)$.

Note that the following holds 
\begin{align*} e^T_1\left(\frac{B_b^TB_b}{n}\right)^{-1}e_1&=O_p(1)\\
\|H_{\kappa_1,a}^{-1}\frac{\bea^T\bea}{n-d-1} H_{\kappa_1,a}^{-1}\|&=O_p(1)
\end{align*}
and conditional on $\mathcal F_{1,n}$, application of Lemma \ref{quad_form_lemma} yields:

\begin{align}
    &\alpha\mathfrak{c}^2_1 e^T_1\left(\frac{B_b^TB_b}{n}\right)^{-1}\frac{B_b^TE_b}{n}H_{\kappa_1,a}^{-1}\frac{\bea^T\bea}{n-d-1} H_{\kappa_1,a}^{-1}\frac{E_b^TB_b}{n}\left(\frac{B_b^TB_b}{n}\right)^{-1}e_1\nonumber\\
    =& n^{-1}\alpha\mathfrak{c}^2_1e^T_1\left(\frac{B_b^TB_b}{n}\right)^{-1}e_1 \tr\left[H_{\kappa_1,a}^{-1}\frac{\bea^T\bea}{n-d-1} H_{\kappa_1,a}^{-1}\right]+o_p(1)\label{quad_V1}
\end{align}

On the other hand 

\begin{equation}
2\alpha\mathfrak{c}^2_1\ba^TH_{\kappa_1,a}^{-1}\frac{\bea^T\bea}{n-d-1} H_{\kappa_1,a}^{-1}\frac{E_b^TB_b}{n}\left(\frac{B_b^TB_b}{n}\right)^{-1}e_1=o_p(1).\label{lin_V1}
\end{equation}
Conditionally on $\mathcal F_{1,n}$, the left-hand side is a centered
Gaussian random variable. Its conditional variance is of the form
$O(n^{-1}R_n)$, with $R_n=O_p(1)$. Hence the left-hand side is
$o_p(1)$. Indeed the variance of the LHS conditioned on $\mathcal F_{1,n}$ is given by
$$4n^{-1}\alpha^2\mathfrak{c}^4_1\ba^TH_{\kappa_1,a}^{-1}\frac{\bea^T\bea}{n-d-1} H_{\kappa_1,a}^{-2}\frac{\bea^T\bea}{n-d-1}H_{\kappa_1,a}^{-1}\ba e^T_1\left(\frac{B_b^TB_b}{n}\right)^{-1}e_1$$

Now $e^T_1\left(\frac{B_b^TB_b}{n}\right)^{-1}e_1=O_p(1)$ (see Lemma \ref{V2_lemma}). The fact that $$\ba^TH_{\kappa_1,a}^{-1}\frac{\bea^T\bea}{n-d-1} H_{\kappa_1,a}^{-2}\frac{\bea^T\bea}{n-d-1}H_{\kappa_1,a}^{-1}\ba=O_p(1)$$ follows from standard functional calculus arguments (see proof of Lemma \ref{func_resolv_2} for example).

Using \eqref{lin_V1}, \eqref{quad_V1}, and \eqref{V_1_decomp} it follows that:
$$V_1(B_a,E_a,B_b,E_b)=V_2(B_a,E_a,B_b)+o_p(1).$$
\end{proof}

\begin{proof}[Proof of Lemma \ref{V2_lemma}]

We start by noting that
\begin{equation}e^T_1\left(\frac{B_b^TB_b}{n}\right)^{-1}e_1 =f^{\Sigma_B}_5(c_2)+o_p(1)\label{f1k2_lim}\end{equation}
 by Assumption \ref{ass:rows}. This yields: 
\begin{align}&n^{-1}\alpha\mathfrak{c}^2_1e^T_1\left(\frac{B_b^TB_b}{n}\right)^{-1}e_1 \tr\left[H_{\kappa_1,a}^{-1}\frac{\bea^T\bea}{n-d-1} H_{\kappa_1,a}^{-1}\right]\nonumber\\
= &n^{-1}\alpha\mathfrak{c}^2_1f^{\Sigma_B}_5(c_2) \tr\left[H_{\kappa_1,a}^{-1}\frac{\bea^T\bea}{n-d-1} H_{\kappa_1,a}^{-1}\right]+o_p(1).\label{V2_first}\end{align}

so that
$$V_2(B_a,E_a,B_b)=V_3(B_a,E_a)+o_p(1)$$

\end{proof}

\begin{proof}[Proof of Lemma \ref{V3_Lemma}]
Let
\[
r_n:=n-d-1,
\qquad
\mathcal R_{a,n}:=
\{\operatorname{rank}(B_a)=d+1\}.
\]
By Assumption~\ref{ass:rows},
$
\mathbb P(\mathcal R_{a,n})\to1.
$
For this argument, $P_{B_a}$ denotes the orthogonal projection onto
the column space of $B_a$, which agrees with
$B_a(B_a^TB_a)^{-1}B_a^T$ on $\mathcal R_{a,n}$.

Define
\[
S_{a,n}:=\frac{\bea^T\bea}{r_n},
\qquad
\phi_{\kappa_1}(x):=\frac{x}{(x+\kappa_1)^2},
\quad x\ge0.
\]
Then
\[
\phi_{\kappa_1}(S_{a,n})
=
H_{\kappa_1,a}^{-1}
\frac{\bea^T\bea}{n-d-1}
H_{\kappa_1,a}^{-1}.
\]
Since $\kappa_1>0$, the function $\phi_{\kappa_1}$ is bounded
and continuous on $[0,\infty)$.

Let $Z_n$ be an $r_n\times p$ matrix with i.i.d.\ standard Gaussian
entries, independent of $(B_a,E_a)$, and set
\[
W_n:=\frac1{r_n}Z_n^TZ_n.
\]
Conditionally on $B_a$, the matrix $P_{B_a}^{\perp}$ is a
deterministic orthogonal projection of rank
$n-\operatorname{rank}(B_a)$. Since $E_a$ is independent of $B_a$
and has i.i.d.\ standard Gaussian entries, its rotational invariance
implies that, whenever $\operatorname{rank}(B_a)=d+1$, the
conditional law of $S_{a,n}$ is the law of $W_n$, independently
of the particular realization of $B_a$.

More precisely, for every Borel set $\mathcal A$ of symmetric
$p\times p$ matrices,
\[
\mathbf 1_{\mathcal R_{a,n}}\,
\mathbb P(S_{a,n}\in\mathcal A\mid B_a)
=
\mathbf 1_{\mathcal R_{a,n}}\,
\mathbb P(W_n\in\mathcal A)
\qquad\text{a.s.}
\]
Taking expectations gives
\[
\mathbb P(S_{a,n}\in\mathcal A,\mathcal R_{a,n})
=
\mathbb P(\mathcal R_{a,n})\,
\mathbb P(W_n\in\mathcal A).
\]
Consequently,
\begin{equation}
\mathcal L(S_{a,n}\mid\mathcal R_{a,n})
=
\mathcal L(W_n).
\label{eq:V3_conditional_Wishart}
\end{equation}

Since $p/r_n\to c_1$, the Lemma \ref{func_resolv_2} yields
\[
\frac1p\tr\{\phi_{\kappa_1}(W_n)\}
\xrightarrow{P}
f_2(\kappa_1,c_1),
\]
where
\[
f_2(\kappa_1,c_1)
=
\int\frac{x}{(x+\kappa_1)^2}\,d\mu_{MP,c_1}(x)
=
f_1(\kappa_1,c_1)
+
\kappa_1\partial_\kappa f_1(\kappa_1,c_1).
\]
Moreover, since $\ba$ is deterministic and
$\|\ba\|^2\to\mathfrak b_A<\infty$, Lemma~\ref{func_resolv_2}
gives
\[
\ba^T\phi_{\kappa_1}(W_n)\ba
-
f_2(\kappa_1,c_1)\|\ba\|^2
\xrightarrow{P}0.
\]

To transfer both convergences to $S_{a,n}$, define, for every
symmetric positive semidefinite matrix $S$,
\[
D_n(S):=
\max\left\{
\left|
\frac1p\tr\{\phi_{\kappa_1}(S)\}
-
f_2(\kappa_1,c_1)
\right|,
\left|
\ba^T\phi_{\kappa_1}(S)\ba
-
f_2(\kappa_1,c_1)\|\ba\|^2
\right|
\right\}.
\]
The preceding two convergences imply
\[
D_n(W_n)\xrightarrow{P}0.
\]
For every $\epsilon>0$, \eqref{eq:V3_conditional_Wishart} gives
\begin{align*}
\mathbb P(D_n(S_{a,n})>\epsilon)
&\le
\mathbb P(\mathcal R_{a,n}^c)
+
\mathbb P(\mathcal R_{a,n})\,
\mathbb P(D_n(W_n)>\epsilon)\\
&\to0.
\end{align*}
Thus, under the original probability law,
\[
\frac1p\tr\{\phi_{\kappa_1}(S_{a,n})\}
=
f_2(\kappa_1,c_1)+o_p(1)
\]
and
\[
\ba^T\phi_{\kappa_1}(S_{a,n})\ba
=
f_2(\kappa_1,c_1)\|\ba\|^2+o_p(1).
\]

Since
\[
\frac pn
=
\frac{p}{r_n}\frac{r_n}{n}
\to c_1\alpha^{-1},
\]
we obtain
\begin{align}
&n^{-1}
\tr\left[
H_{\kappa_1,a}^{-1}
\frac{\bea^T\bea}{n-d-1}
H_{\kappa_1,a}^{-1}
\right]\nonumber\\
&\qquad=
\frac pn\,
\frac1p\tr\{\phi_{\kappa_1}(S_{a,n})\}
=
c_1\alpha^{-1}f_2(\kappa_1,c_1)+o_p(1).
\label{f_2_conv}
\end{align}
Consequently,
\begin{align}
V_{3,1}(B_a,E_a)
&:=
\alpha\mathfrak c_1^2 n^{-1}f^{\Sigma_B}_5(c_2)
\tr\left[
H_{\kappa_1,a}^{-1}
\frac{\bea^T\bea}{n-d-1}
H_{\kappa_1,a}^{-1}
\right]\nonumber\\
&=
c_1\mathfrak c_1^2
f^{\Sigma_B}_5(c_2)f_2(\kappa_1,c_1)+o_p(1).
\label{V3_trace}
\end{align}

Similarly,
\begin{align}
V_{3,2}(B_a,E_a)
&:=
\alpha\mathfrak c_1^2
\ba^TH_{\kappa_1,a}^{-1}
\frac{\bea^T\bea}{n-d-1}
H_{\kappa_1,a}^{-1}\ba\nonumber\\
&=
\alpha\mathfrak c_1^2
f_2(\kappa_1,c_1)\|\ba\|^2+o_p(1).
\label{V3_quad}
\end{align}
Finally, since $\|\ba\|^2\to\mathfrak b_A$, combining
\eqref{V3_trace} and \eqref{V3_quad} yields
\[
V_3(B_a,E_a)
=
c_1\mathfrak c_1^2
f^{\Sigma_B}_5(c_2)f_2(\kappa_1,c_1)
+
\alpha\mathfrak c_1^2
f_2(\kappa_1,c_1)\mathfrak b_A
+
o_p(1)
=
V_4+o_p(1).
\]
\end{proof}

\begin{remark}
In the remaining lemmas, whenever an argument involves the projected
Gram matrix from a single fold $s\in\{a,b\}$, we work on
\[
\mathcal R_{s,n}
:=
\{\operatorname{rank}(B_s)=d+1\}.
\]
By Assumption~\ref{ass:rows},
\[
\mathbb P(\mathcal R_{s,n})\to1.
\]
For arguments involving both folds simultaneously, we work on
\[
\mathcal R_{ab,n}
:=
\mathcal R_{a,n}\cap\mathcal R_{b,n}.
\]
Since
\[
\mathbb P(\mathcal R_{ab,n}^c)
\le
\mathbb P(\mathcal R_{a,n}^c)
+
\mathbb P(\mathcal R_{b,n}^c)
\to0,
\]
this localization does not affect any limiting distribution. On the
corresponding rank events,
\[
\operatorname{rank}(P_{B_s}^\perp)=n-d-1,
\]
and the projected Gaussian Gram matrices have the standard Wishart
representation with $n-d-1$ degrees of freedom. Since the arguments are identical to that used in
Lemma~\ref{V3_Lemma}, we omit
repeated rank-localization arguments henceforth.
\end{remark}

\begin{proof}[Proof of Lemma \ref{W_approx}]

To prove \eqref{W_23} we note that,

\begin{align*}\|\tm-\lambda_1H_{\lambda_1,a}^{-1}\tm\|^2&=\|\tm\|^2-2\kappa_1\tm^TH_{\kappa_1,a}^{-1}\tm+\kappa^2_1\tm^TH_{\kappa_1,a}^{-2}\tm\\
&=(1-2\kappa_1f_1(\kappa_1,c_1)-\kappa^2_1\partial_{\kappa}f_1(\kappa_1,c_1))\|\tm\|^2+o_p(1)\end{align*}
where the last equality follows from Lemma \ref{func_resolv_2}.  On the other hand, we observe that \eqref{W_34} follows by Assumption \ref{ass:rows}, and Assumption \ref{ass:signal}.
\end{proof}

\begin{proof}[Proof of Lemma \ref{delta3_conv_lemma}]

Recalling  that by the definition of $\cin$ in \eqref{htm_defn}, $\cin =(1-\fin\kappa_1)^{-1}$ holds. Thus it is enough to study the convergence of 

\begin{align*}&\sqrt{n}\cin (\lambda_1\ba^T H^{-1}_{\lambda_1,a}\tm-\fin\kappa_1\ba^T\tm)\\
=&\sqrt{n}\kappa_1\cin (\ba^T (\bea^T\bea/(n-d-1)+\kappa_1 I)^{-1}\tm-\fin\ba^T\tm)\\
=&\kappa_1\cin \sqrt{\frac{n}{p}}\sqrt{p}(\ba^T (\bea^T\bea/(n-d-1)+\kappa_1 I)^{-1}\tm-\fin\ba^T\tm)\end{align*} Using the fact that $\lim_{n\rightarrow \infty}p/n=c_1(1-c_2)$, $\cin \rightarrow \mathfrak{c}_1$, Lemma \ref{MP_fluc}, and Assumption \ref{ass:signal} we obtain
$$\sqrt{n}(\cin\ba^T(\tm-\lambda_1H_{\lambda_1,a}^{-1}\tm)-\ba^T\tm)\xrightarrow{d} N(0, \mathcal{V}).$$

\end{proof}

\begin{proof}[Proof of Lemma \ref{T1_conv}]
Denote $\mathcal{X}_a=(B_a,E_a)$ and similarly define $\mathcal{X}_b$. Note that $T_1$ can be rewritten as:

$$T_1=\frac{1}{2}(m(\mathcal{X}_a, \mathcal{X}_b, D_a)+m(\mathcal{X}_b,\mathcal{X}_a,D_b))$$
where 
\begin{align*}
m(\mathcal{X}_a,\mathcal{X}_b,D_a)& =\sqrt{n}\cin(\ba^T+e^T_1(B_b^TB_b)^{-1}B_b^TE_b)H_{\lambda_1,a}^{-1}E_a^TP^{\perp}_{B_a}D_a.
\end{align*}
Conditional on $\mathcal F_{2,n}$, the variables
$m(\mathcal X_a,\mathcal X_b,D_a)$ and
$m(\mathcal X_b,\mathcal X_a,D_b)$ are independent centered
Gaussian random variables, because $D_a$ and $D_b$ are independent
standard Gaussian vectors and are independent of $\mathcal F_{2,n}$.
Since $m(\mathcal X_a,\mathcal X_b,D_a)=\Delta_1$, we have

 $$\mathrm{Var}(m(\mathcal{X}_a,\mathcal{X}_b,D_a)|\mathcal F_{2,n})=V_1(B_a,E_a,B_b,E_b)=V_4+o_p(1)$$ 
where we have used Lemmas \ref{V1_Lemma}, \ref{V2_lemma}, and \ref{V3_Lemma}.
By symmetry $\mathrm{Var}(m(\mathcal{X}_b,\mathcal{X}_a,D_b)|\mathcal F_{2,n})=V_4+o_p(1)$ and hence we have

\[
T_1\mid\mathcal F_{2,n}
\sim N\bigl(0,V_{T_1}(B_a,E_a,B_b,E_b)\bigr).
\]
where \begin{align*}V_{T_1}(B_a,E_a,B_b,E_b)&=\frac{1}{4}(\mathrm{Var}(m(\mathcal{X}_a,\mathcal{X}_b,D_a)|\mathcal F_{2,n})+\mathrm{Var}(m(\mathcal{X}_b,\mathcal{X}_a,D_b)|\mathcal F_{2,n}))\\
&=\frac{1}{2}V_4+o_p(1).\end{align*}
\end{proof}

\begin{proof}[Proof of Lemma \ref{S_1_conv}]
Conditional on $(B_a,B_b)$, the four projected Gaussian blocks
$\bar E_a$, $\breve E_a$, $\bar E_b$, and $\breve E_b$
are mutually independent. Consequently, conditional on
$\mathcal G_n=\sigma(B_a,B_b,\bar E_a,\bar E_b)$,
$\breve E_a$ and $\breve E_b$ remain independent centered
Gaussian matrices. The functions $g(\mathcal X_s)$
and $l(\mathcal X_s)$ depend on $E_s$ only through $\bar E_s$;
in particular, $S_2$ is $\mathcal G_n$-measurable. We have

\[
h^\top(\mathcal X_b)l(\mathcal X_a)\mid\mathcal G_n
\sim N\bigl(0,W_2(B_a,E_a,B_b)\bigr).
\]

with $$W_2(B_a,E_a,B_b)=n||(\tm-\lambda_1H_{\lambda_1,a}^{-1}\tm)||^2\cin^2e_1^T(B_b^TB_b)^{-1}e_1$$

We recall the following expressions and identities from the proof of Lemma \ref{W_approx}

    \begin{align*}
    W_3(B_b)&=n\mathfrak{c}_1^2(1-2\kappa_1f_1(\kappa_1,c_1)-\kappa^2_1\partial_{\kappa} f_1(\kappa_1,c_1))\|\tm\|^2e_1^T(B_b^TB_b)^{-1}e_1\\
    W_4&=\mathfrak{c}_1^2(1-2\kappa_1f_1(\kappa_1,c_1)-\kappa^2_1\partial_{\kappa}f_1(\kappa_1,c_1))f^{\Sigma_B}_5(c_2) \mathfrak{t}_M\\
    W_2(B_a,E_a,B_b)&=W_3(B_b)+o_p(1)\\
        W_3(B_b)&=W_4+o_p(1)
    \end{align*}

Thus observing that $h^T(\mathcal{X}_b)l(\mathcal{X}_a)$ and $h^T(\mathcal{X}_a)l(\mathcal{X}_b)$ are independent conditioned on $\mathcal G_n$, we obtain the desired result.
\end{proof}

\begin{proof}[Proof of Lemma \ref{S_2_conv_lemma}]
As in the proof of Lemma \ref{delta3_conv_lemma}, 
we have

$$\sqrt{n}(\cin\ba^T(\tm-\lambda_1H_{\lambda_1,a}^{-1}\tm)-\ba^T\tm)\xrightarrow{d} N(0, \mathcal{V})$$

where

\begin{align*}\mathcal{V}=&\kappa^2_1\mathfrak{c}^2_1\frac{(f_1(\kappa_1,c_1)\kappa_1)^{4}}{c_1(1-c_2)}\left(h_1(-\kappa_1,c_1)\mathfrak{g}^2+h_1(-\kappa_1,c_1)\mathfrak{b}_A\mathfrak{t}_M\right)\end{align*}

Since $\sqrt{n}(\cin\ba^T(\tm-\lambda_1H_{\lambda_1,a}^{-1}\tm)-\ba^T\tm)$ and $\sqrt{n}(\cin\ba^T(\tm-\lambda_1H_{\lambda_1,b}^{-1}\tm)-\ba^T\tm)$ are independent \eqref{S_2_conv} follows.

\end{proof}

\begin{proof}[Proof of Lemma \ref{tm2_lemma}]
\begin{align*}
    \|\btm\|^2=&\|\tm-\lambda_1H_{\lambda_1,a}^{-1}\tm\|^2+2 (\tm-\lambda_1H_{\lambda_1,a}^{-1}\tm)^TH_{\lambda_1,a}^{-1}\bea^TD_a+D^T_a \bea H_{\lambda_1,a}^{-2}\bea^TD_a\\
    =&\|\tm-\kappa_1H_{\kappa_1,a}^{-1}\tm\|^2+2(n-d-1)^{-1} (\tm-\kappa_1H_{\kappa_1,a}^{-1}\tm)^TH_{\kappa_1,a}^{-1}\bea^T D_a\\
    &+(n-d-1)^{-2}D^T_a \bea H_{\kappa_1,a}^{-2}\bea^T D_a\\
\end{align*}

For the next two noise terms, condition on $(B_a,E_a)$,
so that $D_a$ remains an independent standard Gaussian vector.
We first consider the third term:

\begin{align}
&(n-d-1)^{-2}D^T_a \bea H_{\kappa_1,a}^{-2}\bea^TD_a\nonumber\\
    =&(n-d-1)^{-2}\tr[\bea H_{\kappa_1,a}^{-2}\bea^T]+o_p(1)\nonumber\\
     =&(n-d-1)^{-2}\tr[H_{\kappa_1,a}^{-1}\bea^T\bea H_{\kappa_1,a}^{-1} ]\nonumber +o_p(1)\\
     =& c_1f_2(\kappa_1,c_1)+o_p(1)
\end{align}
with the first equality following from Lemma \ref{quad_form_lemma} and the last equality follows from \eqref{f_2_conv}. Next we observe that

$$2(n-d-1)^{-1} (\tm-\kappa_1H_{\kappa_1,a}^{-1}\tm)^TH_{\kappa_1,a}^{-1}\bea^T D_a=o_p(1)$$
since the LHS has mean $0$ and the conditional variance is of the form $O(n^{-1}R_n)$, with $R_n=O_p(1)$. Finally we consider term $\|\tm-\lambda_1H_{\lambda_1,a}^{-1}\tm\|^2$ for which we have

$$\|\tm-\lambda_1H_{\lambda_1,a}^{-1}\tm\|^2=(1-2\kappa_1f_1(\kappa_1,c_1)-\kappa^2_1\partial_{\kappa}f_1(\kappa_1,c_1))\|\tm\|^2+o_p(1).$$
using Lemma \ref{func_resolv_2}. Using Assumption \ref{ass:signal}, one thus obtains
$$\|\btm\|^2=(1-2\kappa_1f_1(\kappa_1,c_1)-\kappa^2_1\partial_{\kappa}f_1(\kappa_1,c_1))\mathfrak{t}_M+c_1f_2(\kappa_1,c_1)+o_p(1).$$
\end{proof}

\begin{proof}[Proof of Lemma \ref{ba2_lemma}]
We have
\begin{align*}
    \|\hba\|^2=&\|\ba\|^2+2 e^T_1(B_b^TB_b)^{-1}B_b^TE_b\ba+e^T_1(B_b^TB_b)^{-1}B_b^TE_bE_b^TB_b(B_b^TB_b)^{-1}e_1\\
    =&\|\ba\|^2+2 n^{-1}e^T_1(B_b^TB_b/n)^{-1}B_b^TE_b\ba\\
    &+n^{-2}e^T_1(B_b^TB_b/n)^{-1}B_b^TE_bE_b^TB_b(B_b^TB_b/n)^{-1}e_1
\end{align*}

Since conditioned on $B_b$, $E_b^TB_b(B_b^TB_b)^{-1}e_1\sim N(0, e_1^T(B_b^TB_b)^{-1}e_1 I_p)$, we have:

\begin{align}
&n^{-2}e^T_1(B_b^TB_b/n)^{-1}B_b^TE_bE_b^TB_b(B_b^TB_b/n)^{-1}e_1\nonumber\\
    =&pn^{-1}e^T_1(B_b^TB_b/n)^{-1}e_1 +o_p(1)\nonumber\\
     =& c_1(1-c_2)f^{\Sigma_B}_5(c_2) +o_p(1)
\end{align}

The first equality follows from Lemma \ref{quad_form_lemma}. The second equality follows from
$p/n=c_1(1-c_2)+o(1)$ and Assumption \ref{ass:rows}.

Next we observe that 
$$2 n^{-1}e^T_1(B_b^TB_b/n)^{-1}B_b^TE_b\ba=o_p(1)$$
Indeed, conditional on $B_b$, the left-hand side is a centered
Gaussian random variable, and direct computation shows that its
conditional variance is of the form $O(n^{-1}R_n)$, with
$R_n=O_p(1)$. Consequently, the left-hand side is $o_p(1)$.

Using Assumption \ref{ass:signal}, we have

$$\|\hba\|^2=\mathfrak{b}_A+c_1(1-c_2)f^{\Sigma_B}_5(c_2) +o_p(1)$$

\end{proof}

\begin{proof}[Proof of Lemma \ref{tV1_Lemma}]

Expanding \eqref{tV_1_squared_form}, we obtain
  \begin{align}
    &\tV_1(B_a,E_a,B_b,E_b)\\
   =&n\cin^2e^T_jH_{\lambda_1,a}^{-1}\bea^T\bea H_{\lambda_1,a}^{-1}e_j (\baj^2+2\baj e_j^TE_b^TB_b(B_b^TB_b)^{-1}e_1+(e_j^TE_b^TB_b(B_b^TB_b)^{-1}e_1)^2)\\ 
    = 
    &\alpha\mathfrak{c}^2_1 e_j^TH_{\kappa_1,a}^{-1}\frac{\bea^T\bea}{n-d-1} H_{\kappa_1,a}^{-1} e_j\left(\baj^2+2\baj e^T_j\frac{E_b^TB_b}{n}\left(\frac{B_b^TB_b}{n}\right)^{-1}e_1\right.\nonumber\\
    &+\left.\left(e_j^T\frac{E_b^TB_b}{n}\left(\frac{B_b^TB_b}{n}\right)^{-1}e_1\right)^2\right)
    +o_p(1)\label{tV_1_decomp}
\end{align}
where $\alpha:=\lim_{n\rightarrow \infty}\frac{n}{n-d-1}$ and we have used $\cin\rightarrow \mathfrak{c}_1$. Since conditioned on $B_b$, $\frac{E_b^TB_b}{n}\left(\frac{B_b^TB_b}{n}\right)^{-1}e_1\sim n^{-1/2}N(0, e^T_1\left(\frac{B_b^TB_b}{n}\right)^{-1}e_1 I_p)$, we have, conditioned on $B_b$,
$$e^T_j\frac{E_b^TB_b}{n}\left(\frac{B_b^TB_b}{n}\right)^{-1}e_1 \sim n^{-1/2}N(0, e^T_1\left(\frac{B_b^TB_b}{n}\right)^{-1}e_1)$$

Since  $e^T_1\left(\frac{B_b^TB_b}{n}\right)^{-1}e_1=O_p(1)$ and $e_j^TH_{\kappa_1,a}^{-1}\frac{\bea^T\bea}{n-d-1} H_{\kappa_1,a}^{-2}\frac{\bea^T\bea}{n-d-1}H_{\kappa_1,a}^{-1}e_j=O_p(1)$  it follows that:
$$\tV_1(B_a,E_a,B_b,E_b)=\tV_3(B_a,E_a)+o_p(1).$$

where 

$$\tV_3(B_a,E_a)=\alpha\mathfrak{c}^2_1 e_j^TH_{\kappa_1,a}^{-1}\frac{\bea^T\bea}{n-d-1} H_{\kappa_1,a}^{-1} e_j\baj^2$$

Finally note that $$e_j^TH_{\kappa_1,a}^{-1}\frac{\bea^T\bea}{n-d-1} H_{\kappa_1,a}^{-1} e_j=f_2(\kappa_1,c_1)+o_p(1)$$

Thus defining $\tV_4=\alpha\mathfrak{c}^2_1 f_2(\kappa_1,c_1)\baj^2$

$$\tV_3(B_a,E_a)=\tV_4+o_p(1).$$

\end{proof}

\begin{proof}[Proof of Lemma \ref{tW_approx}]
Recall that 
$$\tW_2(B_a,E_a,B_b)=n(e_j^T(I-\lambda_1H_{\lambda_1,a}^{-1})\tm)^2\cin^2e_1^T(B_b^TB_b)^{-1}e_1.$$

Note that  $$e_j^T(I-\lambda_1H_{\lambda_1,a}^{-1})\tm=(1-\kappa_1f_1(\kappa_1,c_1))\theta_{M,j}+o_p(1)$$

so that using $\cin \rightarrow \mathfrak{c}_1$, we obtain 
\begin{align*}
    \tW_2(B_a,E_a,B_b)&=n\mathfrak{c}_1^2(1-\kappa_1f_1(\kappa_1,c_1))^2\tmj^2e_1^T(B_b^TB_b)^{-1}e_1+o_p(1)\\
    &= \tW_3(B_b)+o_p(1)
\end{align*}

On the other hand since $ne_1^T(B_b^TB_b)^{-1}e_1=f^{\Sigma_B}_5(c_2)+o_p(1)$,
we have 
$$\tW_3(B_b)=\tW_4+o_p(1).$$
\end{proof}

\begin{proof}[Proof of Lemma \ref{gamma3_conv_lemma}]
Recalling that $\cin=(1-\fin\kappa_1)^{-1}$ holds, it is enough to study the convergence of 

\begin{align*}&\sqrt{n}\baj\cin(\lambda_1e_j^T H^{-1}_{\lambda_1,a}\tm-\fin\kappa_1e_j^T\tm)\\
=&\sqrt{n}\baj\kappa_1\cin(e_
j^T (\bea^T\bea/(n-d-1)+\kappa_1 I)^{-1}\tm-\fin e_j^T\tm),
\end{align*}
using Lemma \ref{MP_fluc}  and  $\cin \rightarrow \mathfrak{c}_{1}$, we obtain
$$\sqrt{n}(\cin\ba^Te_je_j^T(\tm-\lambda_1H_{\lambda_1,a}^{-1}\tm)-\ba^Te_je_j^T\tm)\xrightarrow{d} N(0, \ctV).$$

\end{proof}

\section{Auxiliary lemmas\label{sec:aux_lemma}}

In this section we prove results for general $\n\times \p$ random matrices. The suitable dimension can be replaced as required.

\subsection{Results on deterministic equivalents}

\begin{lemma}[Adapted from Lemma B.26 of \cite{bai2010spectral}]
\label{quad_form_lemma}
Let $X\in\mathbb R^\p$ have independent, mean-zero,
unit-variance entries satisfying $\E|X_i|^k\le\nu_k$
for the relevant moment orders $k$. Then, for any
deterministic $A\in\mathbb R^{\p\times\p}$ and any
integer $k\ge1$ for which these moments exist,
\[
\E\bigl|X^\top AX-\tr(A)\bigr|^k
\le C_k\left[
\bigl(\nu_4\tr(AA^\top)\bigr)^{k/2}
+\nu_{2k}\bigl(\tr(AA^\top)\bigr)^{k/2}
\right],
\]
where $C_k$ is independent of $\p$. In particular, if
$\|A\|\le1$ and the entries of $X$ have uniformly
bounded eighth moments, then
\[
\E\bigl|X^\top AX-\tr(A)\bigr|^4
\le C\p^2
\]
for a constant $C$ independent of $\p$. Consequently,
when $\p\to\infty$ through positive integers,
\[
\frac{X^\top AX}{\p}-\frac{\tr(A)}{\p}
\xrightarrow{\mathrm{a.s.}}0.
\]
\end{lemma}

\textbf{Fluctuations of bilinear form}
\begin{lemma}{\label{MP_fluc}}
    Let $E$ be a $\n\times \p$ matrix with i.i.d. standard Gaussian entries and $\p/\n \rightarrow c$, for $c\in (0,\infty)$. Further suppose $\ba\in \mathbb R^\p$, $\tm\in \mathbb R^\p$ are sequences of deterministic vectors for which the Assumption \ref{ass:signal} holds. Then for a fixed $\kappa>0$, the following holds:
    $$\sqrt{\p}(\ba^T (E^TE/\n+\kappa I)^{-1}\theta_M-\fingen\ba^T\theta_M)\xrightarrow{d}N(0,f_0(\kappa))$$

where $\fingen$ is defined by \eqref{eq:f1n_def} and
\begin{equation}
f_0(\kappa)
=(\kappa f_1(\kappa,c))^4h_1(-\kappa,c)
\left(\mathfrak g^2+\mathfrak b_A\mathfrak t_M\right).
\label{f_0_def}
\end{equation}
In the proof below, the dependence on the fixed aspect ratio $c$
is suppressed in the functions $g_1$, $g_2$, $h_1$, $h_2$, and $d$.

\end{lemma}

\begin{proof}[Proof of Lemma \ref{MP_fluc}]

We employ Theorem 2.3 of \cite{yin2023limiting} to derive this CLT.
To establish a direct correspondence between their notation and ours,
we first recall their basic setup.

The paper considers random vectors $y$ with elliptically correlated
structure
\[
y=\rho\Gamma\mathbf{u}+\mu,
\]
where $\rho\geq0$ is a random variable describing the radius and
$\mathbf{u}\sim U(S^{\p-1})$ is uniformly distributed on the unit
sphere and independent of $\rho$. By setting
\[
\Gamma=I_\p,\qquad \mu=0,\qquad \rho^2\sim\chi^2_\p,
\]
we recover a random vector with i.i.d.\ standard Gaussian entries.
The authors consider the sample covariance matrix
\[
S_{0\n}
=
\frac{1}{\n}X_\n X_\n^\top,
\]
where the columns of $X_\n$ are i.i.d.\ and
$\mathbf{x}_1\sim\rho_1\mathbf{u}_1$. Thus, in the present Gaussian
specialization, $X_\n^\top$ has the same distribution as $E$ and
$S_{0\n}$ has the same distribution as $E^\top E/\n$.

\textbf{Assumption 1.}
In the notation of \cite{yin2023limiting},
\[
\E[\rho_1^2]=\p,
\qquad
\E[\rho_1^4]=m_\p=\nu_\p+\p^2.
\]
If $\rho_1^2\sim\chi^2_\p$, then
\[
\nu_\p=2\p,
\qquad
m_\p=\p^2+2\p.
\]

\textbf{Assumption 2.}
The matrix $\Gamma_\n$ is a $\p\times\p$ non-random matrix with
uniformly bounded spectral norm. As $\n\to\infty$, the ESD of
\[
\Sigma_\n=\Gamma_\n\Gamma_\n^\top,
\]
denoted by $H_{1\n}$, converges weakly to a proper distribution
$H_1$. Furthermore, the distribution function $H_{2\n}$ of
$m_\p^{-1/2}\rho_1^2$ converges weakly to a distribution $H_2$
with bounded support.

In our specialization $\Gamma_\n=I_\p$, so $H_{1\n}=H_1$ is the
point mass at $1$. Moreover,
\[
m_\p^{-1/2}\rho_1^2
=
(\p^2+2\p)^{-1/2}\rho_1^2
\xrightarrow{P}1,
\]
and hence $H_2$ is also the point mass at $1$.

The authors work with the normalized sample covariance matrix
\[
S_\n
=
\sqrt{\frac{\p^2}{m_\p}}\,S_{0\n}.
\]
In the Gaussian case,
\[
a_\p
:=
\sqrt{\frac{\p^2}{m_\p}}
=
\left(1+\frac{2}{\p}\right)^{-1/2}
=
1+O(\p^{-1}).
\]
We will ultimately replace $S_\n=a_\p S_{0\n}$ by $S_{0\n}$.
This replacement is valid even at the $\sqrt{\p}$ fluctuation scale.
Indeed, for every fixed $\kappa>0$, the resolvent identity gives
\begin{align*}
&\left\|
(a_\p S_{0\n}+\kappa I_\p)^{-1}
-
(S_{0\n}+\kappa I_\p)^{-1}
\right\|\\
&\qquad\leq
|a_\p-1|\,
\|(a_\p S_{0\n}+\kappa I_\p)^{-1}\|\,
\|S_{0\n}\|\,
\|(S_{0\n}+\kappa I_\p)^{-1}\|\\
&\qquad\leq
\kappa^{-2}|a_\p-1|\,\|S_{0\n}\|.
\end{align*}
Since $\|S_{0\n}\|=O_P(1)$ and $|a_\p-1|=O(\p^{-1})$, it follows
that
\begin{equation}
\sqrt{\p}\left\|
(S_\n+\kappa I_\p)^{-1}
-
(S_{0\n}+\kappa I_\p)^{-1}
\right\|
=o_P(1).
\label{eq:Sn_S0n_resolvent}
\end{equation}

For unit vectors $\pi_{1,\n}$ and $\pi_{2,\n}$,
\cite{yin2023limiting} considers the bilinear forms
\[
\pi_{1,\n}^\top\Upsilon(S_\n,z)\pi_{2,\n},
\qquad
\Upsilon(S_\n,z):=(S_\n-zI_\p)^{-1}.
\]
After specializing $\Sigma_\n=I_\p$, Theorem 2.2 of
\cite{yin2023limiting} states that
\begin{equation}
\pi_{1,\n}^\top\Upsilon(S_\n,z)\pi_{2,\n}
+
\frac{1}{z(1+g^0_{2\n}(z))}
\pi_{1,\n}^\top\pi_{2,\n}
\xrightarrow{a.s.}0.
\label{asconv}
\end{equation}
Here $g^0_{1\n}$ and $g^0_{2\n}$ are the finite-$\n$ deterministic
functions obtained from the system
\begin{align*}
zg^0_{1\n}(z)
&=
-c_\n
\int\frac{x}{1+g^0_{2\n}(z)x}\,dH_{1\n}(x),\\
zg^0_{2\n}(z)
&=
-\int\frac{x}{1+g^0_{1\n}(z)x}\,dH_{2\n}(x),
\end{align*}
where
$
c_\n:=\frac{\p}{\n}.
$
For $z>0$, define the finite-$\n$ deterministic equivalent
\begin{equation}
f_{1,\n}(z)
:=
\frac{1}{z(1+g^0_{2\n}(-z))}.
\label{eq:f1n_def}
\end{equation}
Then \eqref{asconv}, evaluated at $-z$, becomes
\begin{equation}
\pi_{1,\n}^\top\Upsilon(S_\n,-z)\pi_{2,\n}
-
f_{1,\n}(z)\pi_{1,\n}^\top\pi_{2,\n}
\xrightarrow{a.s.}0.
\label{eq:asconv_centered}
\end{equation}

We next identify the limit of this deterministic equivalent. Let
$g_1$ and $g_2$ denote the limiting solutions obtained by replacing
$c_\n,H_{1\n},H_{2\n}$ by $c,H_1,H_2$. Since both $H_1$ and $H_2$
are point masses at $1$, the limiting system reduces to
\begin{align*}
zg_1(z)
&=
-\frac{c}{1+g_2(z)},\\
zg_2(z)
&=
-\frac{1}{1+g_1(z)}.
\end{align*}
For $z>0$, define
\[
h(z):=z(1+g_2(-z)).
\]
The preceding equations imply
\[
h(z)-z=\frac{h(z)}{h(z)+c}.
\]
Hence, if
\[
b(z):=h(z)^{-1},
\]
then
\[
czb(z)^2+(1-c+z)b(z)-1=0.
\]
This is precisely the Mar\v{c}enko-Pastur fixed-point equation
\eqref{fixed_pt_1}. Since $b(z)>0$ for $z>0$, it follows that
\begin{equation}
b(z)=f_1(z,c).
\label{eq:b_f1_identification}
\end{equation}
In particular,
\[
f_{1,\n}(\kappa)
\rightarrow
f_1(\kappa,c).
\]
The important point, however, is that the finite-$\n$ quantity
$f_{1,\n}(\kappa)$, rather than its limit $f_1(\kappa,c)$, is the
centering appearing in the CLT below.

We now turn to the fluctuation result. Theorem 2.3 of
\cite{yin2023limiting} considers the two-dimensional process
\begin{align*}
M_{\n}(z)
=
\begin{pmatrix}
M_{\n1}(z)\\
M_{\n2}(z)
\end{pmatrix}
=
\begin{pmatrix}
\sqrt{\p}\left(
\pi_{1,\n}^\top\Upsilon(S_\n,z)\pi_{2,\n}
+
\frac{1}{z(1+g^0_{2\n}(z))}
\pi_{1,\n}^\top\pi_{2,\n}
\right)\\
\sqrt{\p}\left(
\pi_{3,\n}^\top\Upsilon(S_\n,z)\pi_{4,\n}
+
\frac{1}{z(1+g^0_{2\n}(z))}
\pi_{3,\n}^\top\pi_{4,\n}
\right)
\end{pmatrix}.
\end{align*}
For $j,k\in\{1,2,3,4\}$, define
\begin{align*}
r_{jk}(z_1,z_2)
&=
\lim_{\n\rightarrow\infty}
\frac{\pi_{j,\n}^\top\pi_{k,\n}}
{(1+g^0_{2\n}(z_1))(1+g^0_{2\n}(z_2))},\\
r_{jk}(z)
&=
\lim_{\n\rightarrow\infty}
\frac{\pi_{j,\n}^\top\pi_{k,\n}}
{(1+g^0_{2\n}(z))^2}.
\end{align*}
Then $M_\n(z)$ converges weakly to a two-dimensional zero-mean
Gaussian process $M(z)$ with covariance function
\begin{align*}
\operatorname{Cov}(M_1(z_1),M_2(z_2))
=&\,
h_1(z_1,z_2)r_{14}(z_1,z_2)r_{23}(z_1,z_2)\\
&+
h_1(z_1,z_2)r_{13}(z_1,z_2)r_{24}(z_1,z_2)\\
&+
h_2(z_1,z_2)r_{12}(z_1)r_{34}(z_2),
\end{align*}
where
\begin{align*}
h_1(z_1,z_2)
=&
\frac{
c\{z_1g_2(z_1)-z_2g_2(z_2)\}
}{
z_1^2z_2^2
\{g_1(z_1)-g_1(z_2)\}
\{1-d(z_1,z_2)\}
},\\
h_2(z_1,z_2)
=&
\frac{
cg_2'(z_1)g_2'(z_2)
\{\underline m(z_1)g_2(z_2)
-\underline m(z_2)g_2(z_1)\}
}{
g_2(z_1)g_2(z_2)
\{g_1(z_1)-g_1(z_2)\}
},\\
d(z_1,z_2)
=&
\frac{1}{z_1z_2}
\frac{z_1g_1(z_1)-z_2g_1(z_2)}
{g_1(z_1)-g_1(z_2)}
\frac{z_1g_2(z_1)-z_2g_2(z_2)}
{g_2(z_1)-g_2(z_2)}.
\end{align*}

The variance of a single bilinear form is obtained by taking
$\pi_{3,\n}=\pi_{1,\n}$ and $\pi_{4,\n}=\pi_{2,\n}$ and then
letting $z_1\to z_2=z$. Denoting the corresponding diagonal limits
by $h_1(z)$, $h_2(z)$, and $d(z)$, we obtain
\begin{equation}
\operatorname{Var}(M_1(z))
=
h_1(z)r_{12}^2(z)
+
h_1(z)r_{11}(z)r_{22}(z)
+
h_2(z)r_{12}^2(z).
\label{eq:YinZhou_variance}
\end{equation}
The diagonal functions are therefore
\begin{align*}
h_1(z)
&=
\frac{c}{z^4(1-d(z))}
\left[
\frac{g_2(z)}{g_1'(z)}
+
\frac{zg_2'(z)}{g_1'(z)}
\right],\\
h_2(z)
&=
\frac{
c(g_2'(z))^2
\{\underline m'(z)g_2(z)
-\underline m(z)g_2'(z)\}
}{
g_1'(z)(g_2(z))^2
},\\
d(z)
&=
\frac{1}{z^2}
\left(
\frac{g_1(z)}{g_1'(z)}+z
\right)
\left(
\frac{g_2(z)}{g_2'(z)}+z
\right).
\end{align*}

In the present Gaussian identity-covariance specialization,
the companion Stieltjes transform $\underline m$ coincides with
$g_2$. Indeed, the Mar\v{c}enko-Pastur identities give
\[
g_1(z)=c\,m(z),
\qquad
\underline m(z)
=
-\frac{1-c}{z}+c\,m(z)
=
g_2(z).
\]
Consequently,
\[
\underline m(z_1)g_2(z_2)
-
\underline m(z_2)g_2(z_1)
=0,
\]
and hence the two-variable kernel $h_2(z_1,z_2)$ vanishes
identically. In particular,
\begin{equation}
h_2(z)=0.
\label{eq:h2_zero_gaussian}
\end{equation}

We now apply the theorem to the deterministic vectors $\ba$ and
$\tm$. Suppose first that their limiting squared norms
$\mathfrak b_A$ and $\mathfrak t_M$ are strictly positive, and set
\[
\pi_{1,\n}
=
\frac{\ba}{\|\ba\|},
\qquad
\pi_{2,\n}
=
\frac{\tm}{\|\tm\|}.
\]
Then, at $z=-\kappa$,
\begin{align*}
r_{12}(-\kappa)
&=
(1+g_2(-\kappa))^{-2}
\frac{\mathfrak g}
{\sqrt{\mathfrak b_A\mathfrak t_M}},\\
r_{11}(-\kappa)
&=
(1+g_2(-\kappa))^{-2},\\
r_{22}(-\kappa)
&=
(1+g_2(-\kappa))^{-2}.
\end{align*}
Using \eqref{eq:YinZhou_variance} and
\eqref{eq:h2_zero_gaussian}, Theorem 2.3 therefore yields
\begin{align}
&\sqrt{\p}
\left[
\pi_{1,\n}^\top
(S_\n+\kappa I_\p)^{-1}
\pi_{2,\n}
-
f_{1,\n}(\kappa)
\pi_{1,\n}^\top\pi_{2,\n}
\right]
\nonumber\\
&\qquad\xrightarrow{d}
N\left(
0,
h_1(-\kappa)
(1+g_2(-\kappa))^{-4}
\left[
1+
\frac{\mathfrak g^2}
{\mathfrak b_A\mathfrak t_M}
\right]
\right).
\label{eq:YinZhou_normalized_CLT}
\end{align}

Multiplying \eqref{eq:YinZhou_normalized_CLT} by
$\|\ba\|\,\|\tm\|$ and using
\[
\|\ba\|^2\rightarrow\mathfrak b_A,
\qquad
\|\tm\|^2\rightarrow\mathfrak t_M,
\qquad
\ba^\top\tm\rightarrow\mathfrak g,
\]
we obtain by Slutsky's theorem
\begin{align}
&\sqrt{\p}
\left[
\ba^\top
(S_\n+\kappa I_\p)^{-1}
\tm
-
f_{1,\n}(\kappa)\ba^\top\tm
\right]
\nonumber\\
&\qquad\xrightarrow{d}
N\left(
0,
h_1(-\kappa)
(1+g_2(-\kappa))^{-4}
\left[
\mathfrak g^2+
\mathfrak b_A\mathfrak t_M
\right]
\right).
\label{eq:YinZhou_unnormalized_CLT}
\end{align}
It remains to consider the case in which
$\mathfrak b_A=0$ or $\mathfrak t_M=0$. By the Cauchy-Schwarz
inequality, necessarily $\mathfrak g=0$.

Consider an arbitrary subsequence. If along a further subsequence
either $\ba=0$ or $\tm=0$ exactly, then the centered bilinear form is
identically zero on that subsequence. Otherwise both vectors are
nonzero, and we may define
\[
\pi_{1,\n}=\frac{\ba}{\|\ba\|},
\qquad
\pi_{2,\n}=\frac{\tm}{\|\tm\|}.
\]
Passing, if necessary, to a further subsequence, we may assume that
$\pi_{1,\n}^\top\pi_{2,\n}$ converges. The normalized resolvent CLT
above implies that the corresponding centered normalized bilinear form
is $O_p(1)$. Multiplying by
\[
\|\ba\|\,\|\tm\|\rightarrow0
\]
shows that the unnormalized centered bilinear form converges to zero in
probability. Since every subsequence admits a further subsequence with
this property, \eqref{eq:YinZhou_unnormalized_CLT} remains valid, with
zero limiting variance, whenever
$\mathfrak b_A=0$ or $\mathfrak t_M=0$.

Finally, \eqref{eq:Sn_S0n_resolvent} shows that $S_\n$ may be
replaced by $S_{0\n}=E^\top E/\n$ in
\eqref{eq:YinZhou_unnormalized_CLT} without affecting the
$\sqrt{\p}$-scale limit. Therefore
\begin{align}
&\sqrt{\p}
\left[
\ba^\top
\left(\frac{E^\top E}{\n}+\kappa I_\p\right)^{-1}
\tm
-
\fingen\ba^\top\tm
\right]
\nonumber\\
&\qquad\xrightarrow{d}
N(0,f_0(\kappa)),
\label{eq:MP_fluc_finite_center}
\end{align}
where
\[
f_0(\kappa)
=
h_1(-\kappa)
(1+g_2(-\kappa))^{-4}
\left(
\mathfrak g^2+
\mathfrak b_A\mathfrak t_M
\right).
\]
From \eqref{eq:b_f1_identification},
\[
\kappa(1+g_2(-\kappa))
=
\frac{1}{f_1(\kappa,c)},
\]
and therefore
\[
(1+g_2(-\kappa))^{-4}
=
\bigl(\kappa f_1(\kappa,c)\bigr)^4.
\]
Consequently,
\begin{equation}
f_0(\kappa)
=
\bigl(\kappa f_1(\kappa,c)\bigr)^4
h_1(-\kappa)
\left(
\mathfrak g^2+
\mathfrak b_A\mathfrak t_M
\right),
\end{equation}
which is the claimed limiting variance.

\end{proof}

\begin{remark}

We emphasize that the centering in
\eqref{eq:MP_fluc_finite_center} is the finite-$\n$ deterministic
equivalent $f_{1,\n}(\kappa)$ defined in
\eqref{eq:f1n_def}. Since
\[
f_{1,\n}(\kappa)\rightarrow f_1(\kappa,c),
\]
the latter determines the limiting deterministic equivalent, but
mere convergence is not sufficient to replace
$f_{1,\n}(\kappa)$ by $f_1(\kappa,c)$ inside a
$\sqrt{\p}$-scaled statistic. Such a replacement is valid whenever
\[
\sqrt{\p}\,
\left|
f_{1,\n}(\kappa)-f_1(\kappa,c)
\right|
\rightarrow0.
\]
By retaining $f_{1,\n}(\kappa)$ as the centering throughout the
finite-$\n$ debiasing and fluctuation arguments, no additional rate
assumption is required.
Further we note that 
\[\fin:=\fingen\mid_{\n=n-d-1,\p=p, \kappa=\kappa_1}.\]
\end{remark}

To study convergence of quadratic forms of inverse of covariance matrix, we need some more technical tools. 

\medskip

\textbf{Multiplicative free convolution.}

\medskip 

Let $\mu$ be a probability measure on $[0,\infty)$ and $\gamma > 0$. The measure
\[
\nu = \mathrm{MP}(\gamma) \boxtimes \mu
\]
is defined as the unique probability measure whose Stieltjes transform $g_\nu$ satisfies, for all $z \in \mathbb{C}_+$,
\[
g_\nu(z)
=
\int \frac{1}{x(1 - \gamma - \gamma z g_\nu(z)) - z} \, d\mu(x).
\]

\medskip

\textbf{Assumptions on the design}

\medskip

We next give sufficient conditions for the inverse-Gram
deterministic equivalent
\[
e_1^\top S_{\mathfrak n}^{-1}e_1
=
\frac{1}{1-c}
e_1^\top\Sigma_B^{-1}e_1+o_P(1)
\]
and for asymptotic full column rank. These results cover the two
design classes following Remark~\ref{rem:ass_1}. The independent-coordinate
case is proved below; the concentration and lower-edge conditions
for the nonlinear design are verified in Appendix~\ref{sec:lip_conc}.

These verifications establish the inverse-Gram and
rank conditions in Assumption~\ref{ass:rows}.
The additional population limit
\[
e_1^\top\Sigma_B^{-1}e_1\rightarrow\tau
\]
remains a separate requirement of that assumption.

\medskip

Let $B \in \mathbb{R}^{\n \times \p}$ be a random matrix and define the sample covariance matrix
\[
S_\n := \frac{1}{\n} B^\top B, 
\quad 
\Sigma_B := \mathbb{E}[S_\n], \quad G_{S_\n}(z) := (S_\n - z I_\p)^{-1}, \quad z \in (-\infty, 0).
\]

We assume that:

\begin{enumerate}
\item[(i)] $c_\n:=\p/\n \to c \in (0,1)$;
    
     \item[(ii)] $\Sigma_B$ is deterministic, symmetric, and uniformly well-conditioned:
    \[
    0 < \underline{\sigma} \le \lambda_{\min}(\Sigma_B)
    \le \lambda_{\max}(\Sigma_B) \le \overline{\sigma} < \infty;
    \]
    \item[(iii)] The rows of $B$ are i.i.d.\ random vectors, and the matrix
$B$ is $\propto \mathcal{E}_2(1)$-concentrated with respect to the
Frobenius norm (see Appendix~\ref{sec:lip_conc}). Further, define the
whitened design matrix
\[
\mathbf Z:=B\Sigma_B^{-1/2}.
\]
We assume that there exists a deterministic constant $a_0>0$ such that
\[
\mathbb P\left(
\lambda_{\min}\left(
\frac1\n \mathbf Z^\top\mathbf Z
\right)\ge a_0
\right)\rightarrow1.
\]
   \item[(iv)] The rows of $B$ are i.i.d.\ copies of
$\Sigma_B^{1/2}Z$, where the coordinates of $Z$ are independent and identically
distributed according to a common law, independent of $\n$ and $\p$,
with mean zero, variance one,
and finite fourth moment.
\end{enumerate}

Now, define
\[
\nu_\n := \mathrm{MP}(c_\n) \boxtimes \mu_{\Sigma_B},
\qquad
\check{\nu}_\n := (1-c_\n)\delta_0 + c_\n \nu_\n,
\]
and
\[
G(z) := \bigl(-z\, g_{\check{\nu}_\n}(z)\, \Sigma_B - z I_{\p} \bigr)^{-1}.
\]

We have the following theorem adapted from Theorem 2.3 and Proposition 2.4 of \cite{chouard2022quantitative}.

\begin{theorem}
Let $z_\n \in (-\infty,0)$ be a deterministic sequence such that $|z_\n|$ is bounded and
$|z_\n|^{-7} = O(\n)$. Assume that the Assumptions (i), (ii), and (iii) hold.
Then for any deterministic matrix $A \in \mathbb{R}^{\p \times \p}$,
\begin{equation}\bigl|
\operatorname{Tr}(G_{S_\n}(z_\n) A - G(z_\n) A)
\bigr|
\le
\|A\|_F \,
O\!\left( \frac{|z_\n|^{-11/2} (\log \n)^{1/2}}{\n^{1/2}} \right)
\quad \text{a.s.}.\label{dep_det_eq}\end{equation}
\end{theorem}

The above theorem helps us to find the asymptotic limits of the quadratic forms (and bilinear forms in general) of the inverse of the sample covariance matrix when $c<1$.

\begin{lemma}\label{func_resolv_3}
Let $B_\n \in \mathbb{R}^{\n\times \p}$ have i.i.d.\ rows distributed as a random vector $x_\n \in \mathbb{R}^\p$, and define the sample covariance matrix
\[
S_\n := \frac{1}{\n} B_\n^\top B_\n,
\qquad
\Sigma_B := \mathbb{E}[x_\n x_\n^\top].
\]
Assume that the Assumptions (i) and (ii) hold. Further assume that either Assumption (iii) or (iv) holds. Then the following holds: \begin{equation}P(\rm{rank}(B_\n)=\p)\rightarrow 1.\label{full_rank_gen}\end{equation} Further for any deterministic unit vector $u_\n \in \mathbb{R}^\p$,
\begin{equation}
u_\n^\top S_\n^{-1} u_\n
=
\frac{1}{1-c}\, u_\n^\top \Sigma_B^{-1} u_\n + o_P(1).\label{inv_app}
\end{equation}
In particular,
\[
e_1^\top S_\n^{-1} e_1
=
\frac{1}{1-c}\, e_1^\top \Sigma_B^{-1} e_1 + o_P(1).
\]
\end{lemma}

\begin{proof}

\textbf{Case-I:} Assumptions (i),(ii), and (iii) hold:
\vskip 0.2in

Note that the rows of $B$ are of the form $\Sigma_{B}^{1/2}Z$ where $Z$ is isotropic  random vector.
Using Assumptions (ii) and (iii) we can give a lower bound to the smallest singular value of $B$. 
In particular, there exists a deterministic constant $a>0$ such that
\begin{equation}
\mathbb P\left(
\lambda_{\min}\left(
\frac1\n B^\top B
\right)\ge a
\right)
\rightarrow1. \label{singular_low}
\end{equation}
Thus \eqref{full_rank_gen} follows.
Next fix $\alpha \in (0,1/11)$ and define $t_\n := \n^{-\alpha}$ and $z_\n := -t_\n < 0$. By the resolvent identity $A^{-1}-B^{-1}=-A^{-1}(A-B)B^{-1}$ we have:
\[
S_\n^{-1} - (S_\n + t_\n I)^{-1}
=
t_\n S_\n^{-1}(S_\n + t_\n I)^{-1}.
\]
Hence on the event $\{\lambda_{\min}(S_\n) \ge a\}$,
\[
\| S_\n^{-1} - (S_\n + t_\n I)^{-1} \|
\le
\frac{t_\n}{a(a + t_\n)}.
\]

Therefore, on the event
$\{\lambda_{\min}(S_\n)\ge a\}$, for any deterministic unit vector
$u_\n$,
\[
\left|
u_\n^\top S_\n^{-1}u_\n
-
u_\n^\top(S_\n+t_\n I)^{-1}u_\n
\right|
\le
\frac{t_\n}{a(a+t_\n)}.
\]
Since
\[
\mathbb P\{\lambda_{\min}(S_\n)<a\}\to0
\]
by \eqref{singular_low}, and
$t_\n/[a(a+t_\n)]\to0$, it follows that
\begin{equation}
u_\n^\top S_\n^{-1}u_\n
=
u_\n^\top(S_\n+t_\n I)^{-1}u_\n+o_P(1).
\label{inv_app_new}
\end{equation}

On the other hand using  anisotropic deterministic equivalent \eqref{dep_det_eq} and noting that $\frac{|z_\n|^{-11/2} (\log \n)^{1/2}}{\n^{1/2}}=o(1)$ for $|z_\n|=\n^{-\alpha}$, $\alpha \in(0,1/11)$, we obtain

\begin{equation}u_\n^\top (S_\n + t_\n I)^{-1} u_\n=u_\n^\top G(-t_\n) u_\n +o_P(1).\label{inv_app_2}
\end{equation}

Next to deal with the term $u_\n^\top G(-t_\n) u_\n$, we first note that 
\[
G(-t)
=
\big(t\, g_{\check{\nu}_\n}(-t)\Sigma_B + t I_\p \big)^{-1}.
\]
where the following relation holds (c.f. Lemma 6.1 of \cite{chouard2022quantitative})
\[
t g_{\check{\nu}_\n}(-t)
=
1 - c_\n + c_\n t g_{\nu_\n}(-t).
\]

For $t>0$, recall that the definition of Stieltjes transform yields
\[
t_\n g_{\nu_\n}(-t_\n)
=
\int \frac{t_\n}{x+t_\n}\,d\nu_\n(x).
\]
Note that $\exists \varepsilon>0$ such that $\nu_\n([0,\varepsilon])=0$ for large enough $\n$ (see Lemma \ref{hard_edge_free_conv}), so that 
\[
0\le t_\n g_{\nu_\n}(-t_\n)
=
\int_{\varepsilon}^ \infty \frac{t_\n}{x+t_\n}\,d\nu_\n(x)\le \frac{t_\n}{\varepsilon}.
\]
Thus we have,
\[
t_\n g_{\check{\nu}_\n}(-t_\n)
=
(1-c_\n)+c_\n\,t_\n g_{\nu_\n}(-t_\n)\to 1-c.
\]
Define
\[
\alpha_\n:=t_\n g_{\check{\nu}_\n}(-t_\n),
\qquad
A_\n:=t_\n I_\p+\alpha_\n \Sigma_B, \qquad A:=(1-c)\Sigma_B
\]

 Since $\alpha_\n\to 1-c>0$, there exists $N$ such that for all $\n\ge N$, $\alpha_\n\ge (1-c)/2$.
This yields
\[
\lambda_{\min}(A_\n)=\lambda_{\min}(t_\n I_\p+\alpha_\n\Sigma_B)\ge \frac{1-c}{2}\underline{\sigma}.
\]
By the same argument we have
\[
\lambda_{\min}(A)\ge (1-c)\underline{\sigma}>0,
\]
The resolvent identity then gives
\begin{align*}
\|G(-t_\n)-\frac{1}{1-c}\Sigma^{-1}_B\|
&\leq \|A_\n^{-1}-A^{-1}\|\\
&\le
\|A^{-1}\|\,
\|A-A_\n\|\,
\|A_\n^{-1}\|\\
&\leq \frac{2}{((1-c)\underline{\sigma})^2} \|A-A_\n\|
\to 0.
\end{align*}
Thus we have
\[
\left|
u_\n^\top G(-t_\n)u_\n
-
\frac{1}{1-c}u_\n^\top \Sigma_B^{-1}u_\n
\right|
\le
\|u_\n\|^2
\left\|
G(-t_\n)-\frac{1}{1-c}\Sigma_B^{-1}
\right\|
\rightarrow 0.
\]
This proves
\begin{equation}
u_\n^\top G(-t_\n)u_\n
=
\frac{1}{1-c}\,u_\n^\top \Sigma_B^{-1}u_\n+o(1). \label{inv_app_3}
\end{equation}

Combining \ref{inv_app_new}, \ref{inv_app_2}, and \ref{inv_app_3} yields \ref{inv_app} which was to be shown.\\ 

\textbf{Case-II:} Assumptions (i),(ii), and (iv) hold:
\vskip 0.2in

Set
$A_\n:=\frac1\n Z^\top Z$, so that
$
S_\n
=
\Sigma_B^{1/2}A_\n\Sigma_B^{1/2}
$. By the Bai-Yin theorem \cite{bai2010spectral}, since $\p/\n\to c<1$ and the entries have finite fourth
moment,
\[
\lambda_{\min}(A_\n)
\to
(1-\sqrt c)^2
\qquad\text{a.s.}.
\]

Using Assumption (ii), \eqref{full_rank_gen} follows.
Also since
$
S_\n^{-1}
=
\Sigma_B^{-1/2}A_\n^{-1}\Sigma_B^{-1/2}
$
we have
\[
a^\top S_\n^{-1}a
=
b^\top A_\n^{-1}b,
\qquad
b:=\Sigma_B^{-1/2}a.
\]

Define
$
u:=\frac{b}{\|b\|}.
$
so that
\[
b^\top A_\n^{-1}b
=
\|b\|^2 u^\top A_\n^{-1}u.
\]
It remains to prove
\[
u^\top A_\n^{-1}u
\to
\frac{1}{1-c}.
\]

Let
\[
A_\n=\sum_{j=1}^\p \lambda_j v_jv_j^\top
\]
be the spectral decomposition of $A_\n$. Then
\[
u^\top A_\n^{-1}u
=
\sum_{j=1}^\p \frac{|\langle u,v_j\rangle|^2}{\lambda_j}.
\]
Define the eigenvector-weighted empirical spectral distribution
\[
F_{\n,u}(x)
:=
\sum_{j=1}^\p |\langle u,v_j\rangle|^2 \mathbf 1_{\{\lambda_j\le x\}}.
\]
Then
\[
u^\top A_\n^{-1}u
=
\int \frac1x\,dF_{\n,u}(x).
\]

By Theorem 1 of \cite{bai_eigen}, for any deterministic unit vector $u=u_\p$,
the eigenvector-weighted empirical spectral distribution $F_{\n,u}$ converges
weakly almost surely to the Marchenko-Pastur law $F_c$ with parameter $c$.
That is,
\[
F_{\n,u}\Rightarrow F_c
\qquad\text{a.s.}
\]

We would like to apply this convergence to the function $x\mapsto 1/x$. This
function is not bounded on $(0,\infty)$, however we can apply the convergence result by using the fact the lower edge of the spectrum is almost surely bounded away from $0$.
Once again by the Bai-Yin theorem \cite{bai2010spectral}, since $\p/\n\to c<1$ and the entries have finite fourth
moment,
\[
\lambda_{\min}(A_\n)
\to
(1-\sqrt c)^2
\qquad\text{a.s.}
\]
In particular, for any fixed
\[
0<\delta<(1-\sqrt c)^2,
\]
we have eventually almost surely
\[
\lambda_{\min}(A_\n)\ge \delta.
\]
Therefore, outside a null set $\Omega$ for large enough $\n$, the support of $F_{\n,u}$ is contained in
$[\delta,\infty)$.

Choose a bounded continuous function $\varphi_\delta$ such that
\[
\varphi_\delta(x)=\frac1x
\qquad\text{for }x\ge \delta.
\]
Outside the same null set $\Omega$, for all sufficiently large $\n$,
\[
\int \frac1x\,dF_{\n,u}(x)
=
\int \varphi_\delta(x)\,dF_{\n,u}(x).
\]
By weak convergence of $F_{\n,u}$ to $F_c$,
\[
\int \varphi_\delta(x)\,dF_{\n,u}(x)
\to
\int \varphi_\delta(x)\,dF_c(x).
\]
Since the Marchenko-Pastur law with parameter $c<1$ has support
$
[(1-\sqrt c)^2,(1+\sqrt c)^2],
$
and $\delta<(1-\sqrt c)^2$, we also have
\[
\int \varphi_\delta(x)\,dF_c(x)
=
\int \frac1x\,dF_c(x).
\]
The inverse first moment of the Marchenko-Pastur law is
\[
\int \frac1x\,dF_c(x)
=
\frac{1}{1-c}.
\]
Hence
\[
u^\top A_\n^{-1}u
=
\int \frac1x\,dF_{\n,u}(x)
\to
\frac{1}{1-c}
\qquad\text{a.s.}
\]

Consequently,
\[
a^\top S_\n^{-1}a
=
\|b\|^2 u^\top A_\n^{-1}u
=
\frac{\|b\|^2}{1-c}+o_p(1).
\]
Since
$
\|b\|^2
=
a^\top\Sigma_B^{-1}a,
$
we obtain
\[
a^\top S_\n^{-1}a
=
\frac{1}{1-c}a^\top\Sigma_B^{-1}a+o_p(1).
\]
This proves the result.
\end{proof}

\begin{lemma}\label{hard_edge_free_conv}
Let $\mu_{\Sigma_{\p}}$ denote the empirical spectral distribution of a positive definite matrix $\Sigma_{\p}$
such that 
\[
\inf \mathrm{supp}(\mu_{\Sigma_{\p}})\ge \underline{\sigma}>0,
\]
and let $c_\n:=\p/\n\rightarrow c\in(0,1)$. Then the multiplicative free convolution
$
\mu_{\Sigma_{\p}} \boxtimes \mathrm{MP}_{c_\n}
$
has support bounded away from $0$. More precisely,
\[
\liminf_\n \inf \mathrm{supp}(\mu_{\Sigma_{\p}} \boxtimes \mathrm{MP}_{c_\n})> 0.
\]
\end{lemma}

\begin{proof}
Fix $\n$ and $\p$. Choose a sequence of deterministic positive definite matrices $\Sigma_{q(\m)}$ whose empirical spectral distributions converge weakly to $\mu_{\Sigma_\p}$. Note that this can always be done if $q(\m)$ is chosen as multiples of $\p$. Note that in this case
\[
\lambda_{\min}(\Sigma_{q(\m)})\ge \underline{\sigma}
\]
for all $q(\m)$. Let $X_{q(\m)}\in\mathbb R^{q(\m)\times \m}$ have i.i.d.\ entries with mean $0$, variance $1$, and finite fourth moment. Also assume $q(\m)/\m\to c_\n\in(0,1)$. Define
\[
S_{q(\m)}:=\Sigma_{q(\m)}^{1/2}\frac{X_{q(\m)}X_{q(\m)}^\top}{\m}\Sigma_{q(\m)}^{1/2}.
\]
By the standard results on convergence of ESD for covariance matrices \cite{silverstein95} the empirical spectral distribution of $S_{q(\m)}$ converges weakly almost surely to
\[
\mu_{\Sigma_\p}\boxtimes \mathrm{MP}_{c_\n}.
\]

Thus we obtain
\[
\lambda_{\min}(S_{q(\m)})
\ge
\lambda_{\min}(\Sigma_{q(\m)})\,
\lambda_{\min}\!\left(\frac{X_{q(\m)}X_{q(\m)}^\top}{\m}\right)
\ge
\underline{\sigma}\,\lambda_{\min}\!\left(\frac{X_{q(\m)}X_{q(\m)}^\top}{\m}\right).
\]
Because $q(\m)/\m\to c_\n<1$, the standard results on extremes of eigenvalues of covariance matrices (c.f. \cite{bai2010spectral}) yields
\[
\lambda_{\min}\!\left(\frac{X_{q(\m)}X_{q(\m)}^\top}{\m}\right)\to (1-\sqrt c_\n)^2
\qquad\text{a.s.}
\]
Hence
\[
\liminf_{\m\to\infty}\lambda_{\min}(S_{q(\m)})\ge \underline{\sigma}(1-\sqrt c_\n)^2
\qquad\text{a.s.}
\]

Fix $\varepsilon>0$. Then almost surely, for all sufficiently large $\m$, all eigenvalues of $S_{q(\m)}$ lie in
\[
[\underline{\sigma}(1-\sqrt c_\n)^2-\varepsilon,\infty).
\]

Since for every $\varepsilon>0$ the empirical spectral measures of $S_{q(\m)}$ eventually assign zero mass to $(-\infty,\underline{\sigma}(1-\sqrt{c_\n})^2-\varepsilon)$, weak convergence implies that 
\[
\operatorname{supp}(\mu_{\Sigma_\p}\boxtimes \mathrm{MP}_{c_\n})
\subseteq
[\underline{\sigma}(1-\sqrt c_\n)^2,\infty).
\]
Therefore
\[
\inf \operatorname{supp}(\mu_{\Sigma_\p} \boxtimes \mathrm{MP}_{c_\n})\ge \underline{\sigma}(1-\sqrt c_\n)^2>0.
\]
Since $c_\n\to c<1$, this lower bound is bounded away from zero for all sufficiently large $\n$, and the claim follows.
\end{proof}

\begin{lemma}[Bilinear forms of spectral functions of a Gaussian Wishart matrix]
\label{func_resolv_2}
Let $B\in\mathbb R^{\n\times \p}$ have i.i.d.\ $N(0,1)$ entries and set
\[
S:=\frac1\n B^\top B,
\qquad
\frac\p\n\rightarrow c\in(0,\infty).
\]
Let $u_\p,v_\p\in\mathbb R^\p$ be deterministic vectors satisfying
\[
\sup_\p\|u_\p\|<\infty,
\qquad
\sup_\p\|v_\p\|<\infty.
\]
Then, for every bounded continuous function
$\phi:[0,\infty)\to\mathbb R$,
\begin{equation}
u_\p^\top\phi(S)v_\p
-
\frac1\p\tr\{\phi(S)\}\,u_\p^\top v_\p
\xrightarrow{P}0.
\label{eq:isotropic_spectral_function}
\end{equation}
Consequently,
\begin{equation}
u_\p^\top\phi(S)v_\p
=
\left(\int\phi(x)\,d\mu_{MP,c}(x)\right)
u_\p^\top v_\p+o_P(1),
\label{eq:isotropic_spectral_MP}
\end{equation}
where $\mu_{MP,c}$ denotes the Mar\v{c}enko-Pastur law with
aspect ratio $c$.

In particular, for every fixed $\kappa>0$, if
$
G_\kappa:=(S+\kappa I_\p)^{-1},
$
then
\begin{align}
u_\p^\top G_\kappa v_\p
&=
f_1(\kappa,c)\,u_\p^\top v_\p+o_P(1),
\label{eq:G_general}\\
u_\p^\top G_\kappa^2v_\p
&=
-\partial_\kappa f_1(\kappa,c)\,
u_\p^\top v_\p+o_P(1),
\label{eq:G2_general}\\
u_\p^\top G_\kappa S G_\kappa v_\p
&=
f_2(\kappa,c)\,
u_\p^\top v_\p+o_P(1),
\label{eq:GSG_general}
\end{align}
where
\[
f_1(\kappa,c)
=
\int\frac{1}{x+\kappa}\,d\mu_{MP,c}(x),
\qquad
f_2(\kappa,c)
=
f_1(\kappa,c)
+
\kappa\,\partial_\kappa f_1(\kappa,c).
\]
\end{lemma}

\begin{proof}
Since $B$ has i.i.d.\ standard Gaussian entries, its distribution
is invariant under right multiplication by an arbitrary deterministic
orthogonal matrix. Consequently, in the spectral decomposition
$
S=O\Lambda O^\top,
$
the eigenvector matrix $O$ may be chosen Haar distributed on the orthogonal group
$\mathsf O(\p)$ and independent of the eigenvalues
$\Lambda=\diag(\lambda_1,\ldots,\lambda_\p)$.
Equivalently, this follows from the orthogonal invariance of the
right singular vectors of a standard Gaussian matrix.

Fix a bounded continuous function $\phi$ and, conditional on
$\Lambda$, set
\[
D:=\phi(\Lambda)
=
\diag\bigl(\phi(\lambda_1),\ldots,\phi(\lambda_\p)\bigr).
\]
Then
\[
\phi(S)=ODO^\top.
\]

We first consider quadratic forms. Let $w\in\mathbb R^\p$ be
deterministic. The assertion is trivial if $w=0$, so suppose
$w\neq0$ and write
$
x:=O^\top w/\|w\|.
$
Conditional on $\Lambda$, the vector $x$ is uniformly distributed
on the unit sphere $\mathbb S^{\p-1}$. Recall the standard moment
identities
\[
\E[x_i^2]=\frac1\p,
\qquad
\E[x_i^4]=\frac{3}{\p(\p+2)},
\qquad
\E[x_i^2x_j^2]=\frac{1}{\p(\p+2)},
\quad i\neq j.
\]
Therefore, conditional on $\Lambda$,
\begin{align*}
\E\!\left[
w^\top\phi(S)w\mid\Lambda
\right]
&=
\|w\|^2\E[x^\top D x\mid\Lambda]\\
&=
\frac{\tr D}{\p}\|w\|^2.
\end{align*}
Moreover,
\begin{align*}
\mathrm{Var}(x^\top D x\mid\Lambda)
&=
\frac{2}{\p(\p+2)}
\left\{
\tr(D^2)-\frac{(\tr D)^2}{\p}
\right\}\\
&\le
\frac{2}{\p(\p+2)}\tr(D^2)\\
&\le
\frac{2\|D\|^2}{\p+2}.
\end{align*}
Hence
\begin{equation}
\E\!\left[
\left|
w^\top\phi(S)w
-
\frac{\tr D}{\p}\|w\|^2
\right|^2
\Bigm|\Lambda
\right]
\le
\frac{2\|w\|^4\|D\|^2}{\p+2}.
\label{eq:haar_quad_bound}
\end{equation}
Since $\phi$ is bounded,
$
\|D\|\le \|\phi\|_\infty,
$
and therefore the right-hand side of
\eqref{eq:haar_quad_bound} is $O(\p^{-1})$, uniformly in
$\Lambda$, whenever $\|w\|=O(1)$.

We now pass to bilinear forms. Since $\phi(S)$ is symmetric,
the bilinear form can be written as a difference of quadratic forms:
\begin{align*}
u_\p^\top\phi(S)v_\p
=
\frac14
\Big[
(u_\p+v_\p)^\top\phi(S)(u_\p+v_\p)
-
(u_\p-v_\p)^\top\phi(S)(u_\p-v_\p)
\Big].
\end{align*}
Applying \eqref{eq:haar_quad_bound} to $u_\p+v_\p$ and
$u_\p-v_\p$, whose norms are uniformly bounded, yields
\begin{equation}
u_\p^\top\phi(S)v_\p
-
\frac{\tr D}{\p}\,u_\p^\top v_\p
\xrightarrow{P}0.
\end{equation}
Since
$
\tr D=\tr\{\phi(S)\},
$
this proves \eqref{eq:isotropic_spectral_function}.

Next, the empirical spectral distribution of $S$ converges almost
surely to $\mu_{MP,c}$. Since $\phi$ is bounded and continuous,
\[
\frac1\p\tr\{\phi(S)\}
=
\int\phi(x)\,d\mu_S(x)
\rightarrow
\int\phi(x)\,d\mu_{MP,c}(x)
\qquad\text{a.s.}
\]
Combining this convergence with
\eqref{eq:isotropic_spectral_function} proves
\eqref{eq:isotropic_spectral_MP}.

It remains to specialize the general result. Taking
$
\phi_\kappa(x):=(x+\kappa)^{-1}
$
gives
\[
u_\p^\top G_\kappa v_\p
=
f_1(\kappa,c)\,u_\p^\top v_\p+o_P(1).
\]

Next take
$
\psi_\kappa(x):=(x+\kappa)^{-2}.
$
Since $\kappa>0$, the function $\psi_\kappa$ is bounded on
$[0,\infty)$, and hence
\[
u_\p^\top G_\kappa^2v_\p
=
\left(
\int\frac1{(x+\kappa)^2}\,
d\mu_{MP,c}(x)
\right)
u_\p^\top v_\p+o_P(1).
\]
On the other hand,
\[
f_1(\kappa,c)
=
\int\frac1{x+\kappa}\,d\mu_{MP,c}(x),
\]
and differentiation under the integral is justified by dominated
convergence, since
$
\frac1{(x+\kappa)^2}\le \frac1{\kappa^2}.
$
Thus
\[
\partial_\kappa f_1(\kappa,c)
=
-
\int\frac1{(x+\kappa)^2}\,d\mu_{MP,c}(x),
\]
which proves
\[
u_\p^\top G_\kappa^2v_\p
=
-\partial_\kappa f_1(\kappa,c)\,
u_\p^\top v_\p+o_P(1).
\]

Finally, noting that
$
G_\kappa S G_\kappa
=
G_\kappa-\kappa G_\kappa^2,
$
we obtain
\begin{align*}
u_\p^\top G_\kappa S G_\kappa v_\p
&=
\left[
f_1(\kappa,c)
+
\kappa\partial_\kappa f_1(\kappa,c)
\right]
u_\p^\top v_\p
+o_P(1)\\
&=
f_2(\kappa,c)\,
u_\p^\top v_\p+o_P(1).
\end{align*}
This proves the result.
\end{proof}
\subsection{Conditional convergence result}

\begin{lemma}\label{lem:cond_gauss_shift_sum}
Let $Z_n$ be a random element taking values in the measurable space $(\mathcal{Z},
\mathcal{A})$.
Let $f:\mathcal Z\to[0,\infty)$ and $h:\mathcal Z\to\mathbb R$ be measurable functions and define
\[
Y_n := h(Z_n).
\]
Assume that, conditional on $Z_n$, the random variable $X_n$ is Gaussian with mean $0$ and variance
$f(Z_n)$, i.e.
\[
X_n \mid Z_n \sim \mathcal N\!\big(0,f(Z_n)\big).
\]
Suppose further that
\[
f(Z_n)\xrightarrow{P}\sigma_1^2 \qquad\text{for some }\sigma_1^2\ge 0,
\]
and
\[
Y_n \xrightarrow{d} \mathcal N(0,\sigma_2^2)\qquad\text{for some }\sigma_2^2\ge 0.
\]
Then
\[
X_n+Y_n \xrightarrow{d} \mathcal N\!\big(0,\sigma_1^2+\sigma_2^2\big).\]
\end{lemma}

\begin{proof}
Let $\xi\sim \mathcal N(0,1)$ be a standard normal random variable that is independent of $Z_n$ for
each $n$. Define
\[
\widetilde X_n := \sqrt{f(Z_n)}\,\xi .
\]
Conditionally on $Z_n$, $\widetilde X_n$ is centered Gaussian with variance $f(Z_n)$, i.e.
\[
\widetilde X_n \mid Z_n \sim \mathcal N\!\big(0,f(Z_n)\big),
\]
and hence $(\widetilde X_n,Z_n)$ has the same joint law as $(X_n,Z_n)$. In particular,
$\widetilde X_n+Y_n$ has the same law as $X_n+Y_n$. Therefore, it suffices to prove that
\[
\widetilde X_n+Y_n \xrightarrow{d} \mathcal N\!\big(0,\sigma_1^2+\sigma_2^2\big).
\]

By assumption, $f(Z_n)\xrightarrow{P}\sigma_1^2$, and since $x\mapsto \sqrt{x}$ is continuous on
$[0,\infty)$, we also have
\[
\sqrt{f(Z_n)}\xrightarrow{P}\sigma_1.
\]
Because $\xi$ is independent of $Z_n$, it is also independent of $Y_n=h(Z_n)$. Hence the pair
$(\xi,Y_n)$ is independent in the sense that conditional on $Y_n$ the variable $\xi$ remains
standard normal. We claim that
\[
(\widetilde X_n,Y_n) = (\sqrt{f(Z_n)}\,\xi,\; Y_n)
 \xrightarrow{d}\ (\sigma_1\xi,\; Y),
\]
where $Y\sim \mathcal N(0,\sigma_2^2)$ and $Y$ is independent of $\xi$.
To see this, fix $(t_1,t_2)\in\mathbb R^2$ and write the joint characteristic function:
\begin{align*}
\mathbb E\exp\!\big(i t_1\sqrt{f(Z_n)}\,\xi + i t_2 Y_n\big)
&= \mathbb E\Big[\exp(i t_2 Y_n)\,
\mathbb E\big[\exp\!\big(i t_1\sqrt{f(Z_n)}\,\xi\big)\mid Z_n\big]\Big]\\
&= \mathbb E\Big[\exp(i t_2 Y_n)\,
\exp\!\Big(-\frac{t_1^2}{2}f(Z_n)\Big)\Big].
\end{align*}
Since $0\le \exp(-t_1^2 f(Z_n)/2)\le 1$ and $f(Z_n)\xrightarrow{P}\sigma_1^2$, we have
\[
\exp\!\Big(-\frac{t_1^2}{2}f(Z_n)\Big)\xrightarrow{P}
\exp\!\Big(-\frac{t_1^2}{2}\sigma_1^2\Big),
\]
and by dominated convergence,
\[
\mathbb E\Big[\exp(i t_2 Y_n)\,
\exp\!\Big(-\frac{t_1^2}{2}f(Z_n)\Big)\Big]
-
\exp\!\Big(-\frac{t_1^2}{2}\sigma_1^2\Big)\,\mathbb E\exp(i t_2 Y_n)
\to 0.
\]
Using $Y_n\xrightarrow{d} Y\sim \mathcal N(0,\sigma_2^2)$, we get
\[
\mathbb E\exp(i t_2 Y_n)\to \exp\!\Big(-\frac{t_2^2}{2}\sigma_2^2\Big).
\]
Therefore,
\[
\mathbb E\exp\!\big(i t_1\sqrt{f(Z_n)}\,\xi + i t_2 Y_n\big)
\to
\exp\!\Big(-\frac{t_1^2}{2}\sigma_1^2\Big)\,
\exp\!\Big(-\frac{t_2^2}{2}\sigma_2^2\Big),
\]
which is the joint characteristic function of $(\sigma_1\xi,Y)$ with $\xi\perp Y$. This proves the
claimed joint convergence.

Finally, by the continuous mapping theorem applied to the map $(x,y)\mapsto x+y$,
\[
\widetilde X_n + Y_n \xrightarrow{d} \sigma_1\xi + Y.
\]
Since $\sigma_1\xi\sim \mathcal N(0,\sigma_1^2)$, $Y\sim \mathcal N(0,\sigma_2^2)$, and they are
independent, their sum is Gaussian with variance $\sigma_1^2+\sigma_2^2$, i.e.
\[
\sigma_1\xi + Y \sim \mathcal N\!\big(0,\sigma_1^2+\sigma_2^2\big).
\]
Because $\widetilde X_n+Y_n$ has the same law as $X_n+Y_n$, the conclusion follows.
\end{proof}

\subsection{Small-aspect-ratio limits of the asymptotic variance and efficiency results}
\label{subsec:seq_limits_var}

In this subsection we take the variance expression in \eqref{lim_var_defn_explicit} and
compute the sequential limit
\[
\lim_{\kappa_1\to 0}\ \lim_{c_1\to 0}\ \lim_{c_2\to 0}\ \sigma^2,
\]

\begin{lemma}{\label{low_dim_var_lemma}}
    Let $\sigma^2$ be the limiting variance  defined as in \eqref{lim_var_defn_explicit}. Then the following limit holds

\begin{equation}\label{eq:seq_limit_final}
\lim_{\kappa_1\to 0}\ \lim_{c_1\to 0}\ \lim_{c_2\to 0}\ \sigma^2
=
\mathfrak b_A+F_5(0)\,\mathfrak t_M.
\end{equation}
where 
\[
F_5(0)=\lim_{n\to \infty}e_1^\top\Sigma_B^{-1}e_1 = \tau\in(0,\infty).
\]
\end{lemma}
\begin{proof}
For notational convenience, write
\[
F_5(c_2):=f^{\Sigma_B}_5(c_2),
\qquad
f_1:=f_1(\kappa_1,c_1),
\qquad
f_1':=\partial_{\kappa} f_1(\kappa_1,c_1),
\qquad
f_2:=f_2(\kappa_1,c_1)=f_1+\kappa_1 f_1'.
\]
Then \eqref{lim_var_defn_explicit} can be written as
\begin{equation*}
\sigma^2
=
\Big(1-\kappa_1 f_1\Big)^{-2}\Big[T_1(c_1,c_2,\kappa_1)+T_2(c_1,c_2,\kappa_1)+T_3(c_1,c_2,\kappa_1)\Big],
\end{equation*}
where
\begin{align}
T_1
&:= f_2\Big(c_1 F_5(c_2)+(1-c_2)^{-1}\mathfrak b_A\Big),
\label{eq:T1_def}\\
T_2
&:= \Big(1-2\kappa_1 f_1-\kappa_1^2 f_1'\Big)\,F_5(c_2)\,\mathfrak t_M,
\label{eq:T2_def}\\
T_3
&:= \kappa_1^2\,\frac{(\kappa_1 f_1)^{4}}{c_1(1-c_2)}
h_1(-\kappa_1,c_1)\big(\mathfrak g^2+\mathfrak b_A\mathfrak t_M
\Big).
\label{eq:T3_def}
\end{align}

\textbf{The $c_2\rightarrow 0$ limit}

\medskip
Letting $c_2\to 0$ first, we use $(1-c_2)^{-1}\to 1$ and $F_5(c_2)\to F_5(0)$ to obtain
\begin{align}
T_1 &\to f_2\big(c_1 F_5(0)+\mathfrak b_A\big),
\label{eq:T1_c2}\\
T_2 &\to \Big(1-2\kappa_1 f_1-\kappa_1^2 f_1'\Big)\,F_5(0)\,\mathfrak t_M,
\label{eq:T2_c2}\\
T_3 &\to \kappa_1^2\,\frac{(\kappa_1 f_1)^{4}}{c_1}
h_1(-\kappa_1,c_1)\big(\mathfrak g^2+\mathfrak b_A\mathfrak t_M\big).
\label{eq:T3_c2}
\end{align}

\textbf{The $c_1\rightarrow 0$ limit}

\medskip

We now let $c_1\to 0$ with $\kappa_1>0$ fixed. In this regime the MP law degenerates to a point mass at $1$,
so the function $f_1(\kappa,c)$ satisfies
\begin{equation}\label{eq:f1_c0}
f_1(\kappa,0)=\frac{1}{1+\kappa},
\qquad
\partial_{\kappa} f_1(\kappa,0)=-\frac{1}{(1+\kappa)^2},
\qquad
f_2(\kappa,0)=\frac{1}{(1+\kappa)^2}.
\end{equation}
Moreover,
\begin{equation}\label{eq:prefactor_c0}
\Big(1-\kappa f_1(\kappa,0)\Big)^{-2}
=
(1+\kappa)^2.
\end{equation}

\paragraph{Limit of the $T_1$ contribution.}
Using \eqref{eq:T1_c2}, $c_1F_5(0)\to 0$, and \eqref{eq:f1_c0}-\eqref{eq:prefactor_c0},
\[
\Big(1-\kappa_1 f_1(\kappa_1,c_1)\Big)^{-2}\,T_1
\ \rightarrow\
(1+\kappa_1)^2\cdot \frac{\mathfrak b_A}{(1+\kappa_1)^2}
=
\mathfrak b_A.
\]

\paragraph{Limit of the $T_2$ contribution.}
Using \eqref{eq:T2_c2} and \eqref{eq:f1_c0}, we compute
\[
1-2\kappa f_1(\kappa,0)-\kappa^2\partial_{\kappa}f_1(\kappa,0)
=
1-\frac{2\kappa}{1+\kappa}+\frac{\kappa^2}{(1+\kappa)^2}
=
\frac{1}{(1+\kappa)^2}.
\]
Combining with \eqref{eq:prefactor_c0} yields
\[
\Big(1-\kappa_1 f_1(\kappa_1,c_1)\Big)^{-2}\,T_2
\ \rightarrow\
(1+\kappa_1)^2\cdot \frac{1}{(1+\kappa_1)^2}\,F_5(0)\,\mathfrak t_M
=
F_5(0)\,\mathfrak t_M.
\]

\paragraph{Limit of the $T_3$ contribution.}
By Lemma~\ref{h1_h2_lem}, for every fixed $\kappa_1>0$,
\[
\frac{h_1(-\kappa_1,c_1)}{c_1}
\rightarrow
\frac{1}{\kappa_1^4},
\qquad
c_1\to0.
\]
Therefore using
$f_1(\kappa_1,c_1)
\rightarrow
\frac{1}{1+\kappa_1}$, we obtain
\begin{align*}
\lim_{c_1\to0}T_3
&=
\kappa_1^2
\left(\frac{\kappa_1}{1+\kappa_1}\right)^4
\frac{1}{\kappa_1^4}
\left(
\mathfrak g^2+\mathfrak b_A\mathfrak t_M
\right)\\
&=
\frac{\kappa_1^2}{(1+\kappa_1)^4}
\left(
\mathfrak g^2+\mathfrak b_A\mathfrak t_M
\right).
\end{align*}
Since
\[
\left(
1-\kappa_1f_1(\kappa_1,c_1)
\right)^{-2}
\rightarrow
(1+\kappa_1)^2,
\]
the complete contribution of $T_3$ to the variance satisfies
\[
\lim_{c_1\to0}
\left(
1-\kappa_1f_1(\kappa_1,c_1)
\right)^{-2}T_3
=
\frac{\kappa_1^2}{(1+\kappa_1)^2}
\left(
\mathfrak g^2+\mathfrak b_A\mathfrak t_M
\right).
\]
Combining the three contributions gives
\[
\lim_{c_1\to0}\lim_{c_2\to0}\sigma^2
=
\mathfrak b_A
+
F_5(0)\mathfrak t_M
+
\frac{\kappa_1^2}{(1+\kappa_1)^2}
\left(
\mathfrak g^2+\mathfrak b_A\mathfrak t_M
\right).
\]
Letting $\kappa_1\to0$ proves
\[
\lim_{\kappa_1\to0}
\lim_{c_1\to0}
\lim_{c_2\to0}\sigma^2
=
\mathfrak b_A+F_5(0)\mathfrak t_M.
\]
\end{proof}

\begin{lemma}[Small-$c$ behavior of $h_1$]
\label{h1_h2_lem}
Let $K$ be a compact subset of $\mathbb C\setminus[0,\infty).$
For $c>0$ sufficiently small, interpret the square root in
\begin{equation}
f_1(-z,c)
=
\frac{(1-c-z)-\sqrt{(1-c-z)^2-4cz}}{2cz}
\label{eq:f1_minus_z_branch}
\end{equation}
as the analytic branch satisfying
\begin{equation}
\sqrt{(1-c-z)^2-4cz}
\rightarrow
1-z,
\qquad c\downarrow0,
\label{eq:sqrt_branch_convention}
\end{equation}
uniformly for $z\in K$.

Then, uniformly for $z\in K$,
\begin{align}
f_1(-z,c)
&=
\frac{1}{1-z}
+
\frac{c}{(1-z)^3}
+
O_K(c^2),
\label{eq:F_small_c}\\
\partial_z f_1(-z,c)
&=
\frac{1}{(1-z)^2}
+
\frac{3c}{(1-z)^4}
+
O_K(c^2).
\label{eq:Fprime_small_c}
\end{align}
Moreover,
\begin{equation}
h_1(z,c)
=
\frac{c}{z^4}+O_K(c^2).
\label{eq:h1_small_c}
\end{equation}
In particular, for every fixed $\kappa>0$,
\[
\lim_{c\downarrow0}
\frac{h_1(-\kappa,c)}{c}
=
\frac{1}{\kappa^4}.
\]
\end{lemma}

\begin{proof}
We first make precise the square-root convention in
\eqref{eq:f1_minus_z_branch}. Write
\[
\Delta_c(z)
:=
(1-c-z)^2-4cz.
\]

Define
\[
u_c(z)
:=
\frac{-2c(1+z)+c^2}{(1-z)^2}.
\]
so that
\[
\Delta_c(z)
=
(1-z)^2
\left(
1+
u_c(z)
\right),
\]

Because
$K\subset\mathbb C\setminus[0,\infty)$,
the quantity $1-z$ is uniformly bounded away from zero on $K$.
Hence
\[
\sup_{z\in K}|u_c(z)|\rightarrow0,
\qquad c\downarrow0.
\]
For all sufficiently small $c$, define
\begin{equation}
s_c(z)
:=
(1-z)\bigl(1+u_c(z)\bigr)^{1/2},
\label{eq:sc_def}
\end{equation}
where the square root on the right-hand side denotes the analytic
branch in a neighborhood of $1$ satisfying $1^{1/2}=1$.
Then
\[
s_c(z)^2=\Delta_c(z),
\qquad
s_c(z)\rightarrow1-z
\]
uniformly on $K$. Thus $s_c(z)$ is precisely the square-root branch
used in \eqref{eq:f1_minus_z_branch}. In particular, when
$z=-\kappa<0$,
\[
s_c(-\kappa)
=
\sqrt{(1-c+\kappa)^2+4c\kappa}>0,
\]

Using \eqref{eq:sc_def} and the Taylor expansion of the square root
around $1$, uniformly for $z\in K$,
\begin{align*}
s_c(z)
&=
(1-z)
\left[
1+\frac12u_c(z)+O_K(c^2)
\right]\\
&=
(1-z)
-c\frac{1+z}{1-z}
+O_K(c^2).
\end{align*}
Rationalizing \eqref{eq:f1_minus_z_branch} gives
\begin{equation}
f_1(-z,c)
=
\frac{2}
{(1-c-z)+s_c(z)}.
\label{eq:f1_rationalized}
\end{equation}
Moreover,
\begin{align*}
(1-c-z)+s_c(z)
&=
2(1-z)
-c\left(
1+\frac{1+z}{1-z}
\right)
+O_K(c^2)\\
&=
2(1-z)
-\frac{2c}{1-z}
+O_K(c^2).
\end{align*}
Therefore
\begin{align*}
f_1(-z,c)
&=
\frac{1}
{(1-z)-c(1-z)^{-1}+O_K(c^2)}\\
&=
\frac{1}{1-z}
+
\frac{c}{(1-z)^3}
+
O_K(c^2),
\end{align*}
which proves \eqref{eq:F_small_c}.

For all sufficiently small $c$, the functions involved are analytic
in $z$ on a fixed open neighborhood of $K$. Hence the preceding
expansion may be differentiated uniformly on $K$, yielding
\[
\partial_z f_1(-z,c)
=
\frac{1}{(1-z)^2}
+
\frac{3c}{(1-z)^4}
+
O_K(c^2).
\]

Indeed choose a compact set $K_+$ such that
\[
K\subset \operatorname{int}(K_+)
\subset \mathbb C\setminus[0,\infty).
\]
For all sufficiently small $c$, the chosen square-root branch and hence
$f_1(-z,c)$ are analytic in $z$ on a fixed open neighborhood of
$K_+$. The expansion above holds uniformly on $K_+$ and can be
written as
\[
f_1(-z,c)
=
\frac{1}{1-z}
+
\frac{c}{(1-z)^3}
+
c^2R_c(z),
\]
where $R_c$ is analytic and
\[
\sup_{z\in K_+}|R_c(z)|=O(1).
\]
By Cauchy's estimate,
\[
\sup_{z\in K}|\partial_zR_c(z)|=O(1).
\]
Consequently the expansion may be differentiated termwise on $K$,
giving
\[
\partial_z f_1(-z,c)
=
\frac{1}{(1-z)^2}
+
\frac{3c}{(1-z)^4}
+
O_K(c^2).
\]
This proves \eqref{eq:Fprime_small_c}.

Recall now that
\[
g_1(z,c)=c f_1(-z,c).
\]
It follows that
\begin{align}
g_1(z,c)
&=
\frac{c}{1-z}
+
O_K(c^2),
\label{eq:g1_small_c}\\
\partial_z g_1(z,c)
&=
\frac{c}{(1-z)^2}
+
O_K(c^2).
\label{eq:g1prime_small_c}
\end{align}
Consequently,
\begin{equation}
\frac{g_1(z,c)}{\partial_zg_1(z,c)}+z
=
1+O_K(c).
\label{eq:g1_ratio_small_c}
\end{equation}

Next, by definition,
\[
g_2(z,c)
=
-z^{-1}\bigl(f_1(-z,c)\bigr)^{-1}-1.
\]
From \eqref{eq:F_small_c},
\[
\bigl(f_1(-z,c)\bigr)^{-1}
=
(1-z)-\frac{c}{1-z}+O_K(c^2),
\]
and hence
\begin{equation}
g_2(z,c)
=
-\frac1z
+
\frac{c}{z(1-z)}
+
O_K(c^2).
\label{eq:g2_small_c}
\end{equation}
Differentiating gives
\begin{equation}
\partial_z g_2(z,c)
=
\frac1{z^2}
-
c\frac{1-2z}{z^2(1-z)^2}
+
O_K(c^2).
\label{eq:g2prime_small_c}
\end{equation}
Therefore
\begin{equation}
g_2(z,c)+z\,\partial_zg_2(z,c)
=
\frac{c}{(1-z)^2}
+
O_K(c^2).
\label{eq:g2_combo_small_c}
\end{equation}
Since
\[
\partial_zg_2(z,c)
=
\frac1{z^2}+O_K(c),
\]
and $K$ is bounded away from $0$, we obtain
\begin{equation}
\frac{g_2(z,c)}{\partial_zg_2(z,c)}+z
=
\frac{cz^2}{(1-z)^2}
+
O_K(c^2).
\label{eq:g2_ratio_small_c}
\end{equation}

Recall that
\[
d(z,c)
=
\frac{1}{z^2}
\left(
\frac{g_1(z,c)}{\partial_zg_1(z,c)}+z
\right)
\left(
\frac{g_2(z,c)}{\partial_zg_2(z,c)}+z
\right).
\]
Using \eqref{eq:g1_ratio_small_c} and
\eqref{eq:g2_ratio_small_c}, we obtain
\begin{equation}
d(z,c)
=
\frac{c}{(1-z)^2}
+
O_K(c^2).
\label{eq:d_small_c}
\end{equation}
In particular,
\[
1-d(z,c)=1+O_K(c),
\]
and hence $1-d(z,c)$ is uniformly bounded away from zero on $K$
for all sufficiently small $c$.

Finally, from the definition of $h_1$,
\begin{align*}
h_1(z,c)
&=
\frac{c}{z^4(1-d(z,c))}
\left[
\frac{g_2(z,c)}
     {\partial_zg_1(z,c)}
+
\frac{z\,\partial_zg_2(z,c)}
     {\partial_zg_1(z,c)}
\right]\\
&=
\frac{c}{z^4(1-d(z,c))}
\frac{g_2(z,c)+z\,\partial_zg_2(z,c)}
     {\partial_zg_1(z,c)}.
\end{align*}
By \eqref{eq:g2_combo_small_c} and
\eqref{eq:g1prime_small_c},
\[
\frac{g_2(z,c)+z\,\partial_zg_2(z,c)}
     {\partial_zg_1(z,c)}
=
\frac{
c(1-z)^{-2}+O_K(c^2)
}{
c(1-z)^{-2}+O_K(c^2)
}
=
1+O_K(c).
\]
Combining this with \eqref{eq:d_small_c} yields
\[
h_1(z,c)
=
\frac{c}{z^4}\bigl(1+O_K(c)\bigr)
=
\frac{c}{z^4}+O_K(c^2).
\]
This proves \eqref{eq:h1_small_c}.
\end{proof}

\begin{lemma}[Fisher information for the mediation functional]
\label{fisher_info}
Fix dimensions $p$ and $d$ and consider the Gaussian linear
structural equation model
\begin{align}
\stM &= \beta_A\stA+\beta_C\stC+\epsilon,\\
\stY &= \theta_M^\top\stM+\theta_A\stA+\theta_C^\top\stC+\delta,
\end{align}
where $\epsilon\sim N(0,I_p)$ and $\delta\sim N(0,1)$
are mutually independent and independent of $(\stA,\stC)$.
Assume $\E[(\stA,\stC^\top)^\top]=0$ and the covariance matrix
\[
\Sigma_B:=\E\!\left[
\begin{pmatrix}\stA\\\stC\end{pmatrix}
\begin{pmatrix}\stA&\stC^\top\end{pmatrix}
\right]
\]
is positive definite. Let the parameter be
\[
\eta=(\theta_M,\theta_A,\theta_C,\beta_A,\beta_C)
\]
and define the mediation functional $\gamma=\beta_A^\top\theta_M$.
Then its efficient asymptotic variance for estimating $\gamma$ is:
\[
V_\gamma
=\|\beta_A\|^2
+\frac{\|\theta_M\|^2}{\sigma_{A\mid C}^2},
\]
where
\[
\sigma_{A\mid C}^2
:=\E[\stA^2]
-\E[\stA\stC^\top]
\E[\stC\stC^\top]^{-1}
\E[\stC\stA]>0.
\]
Here $\sigma_{A\mid C}^2$ is the population residual
variance from the best linear projection of $\stA$ onto
$\stC$. When $\|\beta_A\|^2
+
\|\theta_M\|^2>0$, i.e. when $V_{\gamma}\neq 0$, the Fisher information for $\gamma$ is
$
I_\gamma = V_\gamma^{-1}.
$
\end{lemma}

\begin{proof}

The joint log-likelihood of $(\stY,\stM)$ conditional on $(\stA,\stC)$ factorizes as
\[
\ell
=
\ell_{\stM}(\beta_A,\beta_C;\stM\mid \stA,\stC)
+
\ell_{\stY}(\theta_M,\theta_A,\theta_C;\stY\mid \stM,\stA,\stC).
\]
It is easy to see from the above factorization that the scores with respect to $\beta_A$ and $\theta_M$ are orthogonal,
so the Fisher information matrix is block diagonal:
\[
I(\eta)=
\begin{pmatrix}
I_\theta & 0\\
0 & I_\beta
\end{pmatrix},
\]
where $I_\theta$ corresponds to $(\theta_M,\theta_A,\theta_C)$ and $I_\beta$ to $(\beta_A,\beta_C)$.

\medskip
\noindent
 First we compute $I_{\beta}^{-1}$. Let $\stZ=(\stA,\stC^\top)^\top$. Then a routine calculation shows
\[
I_\beta = -\E[\nabla^2_{\beta}\ell]=\E[\stZ\stZ^\top] \otimes I,
\qquad
I_\beta^{-1} = \E[\stZ\stZ^\top]^{-1} \otimes I.
\]
Writing
\[
\E[\stZ\stZ^\top]
=
\begin{pmatrix}
a & b^\top\\
b & \Sigma_{CC}
\end{pmatrix},
\quad
a=\E[\stA^2],\; b=\E[\stC\stA],\; \Sigma_{CC}=\E[\stC\stC^\top].
\]
By the block inversion formula, the block corresponding to $\beta_A$ is given by
$
(a-b^\top\Sigma_{CC}^{-1}b)^{-1}I_p.
$
Hence
\[
\big(I_\beta^{-1}\big)_{\beta_A\beta_A}
=
\frac{1}{\sigma_{A\mid C}^2}\,I,
\quad
\sigma_{A\mid C}^2 = a - b^\top \Sigma_{CC}^{-1} b.
\]

\medskip
\noindent

Next we compute $I_{\theta}^{-1}$. Let $\stW=(\stM^\top,\stA,\stC^\top)^\top$. Then once again routine computation gives
$I_\theta^{-1} = \E[\stW\stW^\top]^{-1}$.
Partition $\E[\stW\stW^\top]$ as follows:
\[
\E[\stW\stW^\top]
=
\begin{pmatrix}
\E[\stM\stM^\top] & \E[\stM\stV^\top]\\
\E[\stV\stM^\top] & \E[\stV\stV^\top]
\end{pmatrix},
\quad
\stV=(\stA,\stC^\top)^\top.
\]
By block inversion, the $\theta_M$-block of $\E[\stW\stW^\top]^{-1}$ is
\[
\Big(\E[\stM\stM^\top] - \E[\stM\stV^\top]\E[\stV\stV^\top]^{-1}\E[\stV\stM^\top]\Big)^{-1}.
\]
Since
\[
\stM = \beta_A \stA + \beta_C \stC + \epsilon,
\]
the residual of $\stM$ after projection on $(\stA,\stC)$ is $\epsilon$, so
\[
\E[\stM\stM^\top] - \E[\stM\stV^\top]\E[\stV\stV^\top]^{-1}\E[\stV\stM^\top]
=
I.
\]
Therefore
\[
\big(I_\theta^{-1}\big)_{\theta_M\theta_M}
=
I.
\]

The gradient of $\gamma$ with respect to $\eta$ is
\[
\nabla_\eta \gamma
=
(\beta_A,0,0,\theta_M,0)^T
\]
Thus
\[
V_\gamma
=
\nabla_\eta \gamma^\top I(\eta)^{-1}\nabla_\eta \gamma
=
\beta_A^\top \big(I_\theta^{-1}\big)_{\theta_M\theta_M} \beta_A
+
\theta_M^\top \big(I_\beta^{-1}\big)_{\beta_A\beta_A} \theta_M,
\]
which yields
\[
V_\gamma
=
\|\ba\|^2
+
\frac{1}{\sigma_{A\mid C}^2}\,\|\tm\|^2.
\]
\end{proof}

\section{Verification of the design assumptions}{\label{sec:lip_conc}}

 We say that a sequence of random matrices $B_n\in\mathbb R^{n\times(d+1)}$
is $\propto\mathcal E_2(1)$-concentrated with respect to the Frobenius norm if
there exist constants $C,c>0$, independent of $n$ and $d$, such that for every
$1$-Lipschitz function
$f:\mathbb R^{n\times(d+1)}\to\mathbb R$,
\[
\mathbb P\left(|f(B_n)-\mathbb E f(B_n)|>t\right)
\le
C e^{-c t^2},
\qquad t\ge0.
\]

In this section, we verify that a nonlinear row model constructed through a
logistic link satisfies this concentration assumption. In particular, for the
design matrix $B_n$ constructed below, we show that $B_n$ is
$\propto\mathcal E_2(1)$-concentrated with respect to the Frobenius norm.

\subsection{A logistic-link column model under general dependence}

Let $Y\sim N(0,I_d)$, and let $\Sigma=\Sigma_d\in\mathbb R^{d\times d}$ be
deterministic with
$
\|\Sigma\|\le C_\Sigma.
$
Let $h:\mathbb R^d\to\mathbb R$ be globally Lipschitz with constant $L_h$, and
define
\[
\sigma(t)=\frac{1}{1+e^{-t}}.
\]
Define
\[
X=
\begin{pmatrix}
\sigma\!(h(\Sigma^{1/2}Y))-\E[\sigma\!\big(h(\Sigma^{1/2}Y)\big)]\\
\Sigma^{1/2}Y
\end{pmatrix}
\in\mathbb R^{d+1}.
\]
Let $X_1,\ldots,X_n$ be i.i.d. copies of $X$, and define
\[
\mathbf X=(X_1,\ldots, X_n)^T
\in\mathbb R^{n\times(d+1)}.
\]

\subsubsection{Concentration properties}

\begin{definition}
A random vector $Z$ in a Euclidean space satisfies Gaussian concentration with
constant $K$ if, for every $1$-Lipschitz function $\varphi$,
\[
\mathbb P\left(|\varphi(Z)-\mathbb E\varphi(Z)|\ge t\right)
\le
2\exp\left(-\frac{t^2}{2K^2}\right),
\qquad t\ge0.
\]
\end{definition}

Define
\[
F(y)=
\begin{pmatrix}
\sigma\!\big(h(\Sigma^{1/2}y)\big)-\E[\sigma\!\big(h(\Sigma^{1/2}Y)\big)]\\
\Sigma^{1/2}y
\end{pmatrix}.
\]

\begin{lemma}[Lipschitz property]
The map $F:\mathbb R^d\to\mathbb R^{d+1}$ is globally Lipschitz with constant
\[
K_x
=
\|\Sigma^{1/2}\|
\sqrt{1+\frac{L_h^2}{16}}.
\]
\end{lemma}

\begin{proof}
For $y,y'\in\mathbb R^d$,
\[
\|F(y)-F(y')\|^2
=
\left|
\sigma(h(\Sigma^{1/2}y))
-
\sigma(h(\Sigma^{1/2}y'))
\right|^2
+
\|\Sigma^{1/2}(y-y')\|^2.
\]
Since $\sup_t |\sigma'(t)|\le 1/4$,
\[
|\sigma(u)-\sigma(v)|\le \frac14 |u-v|.
\]
Using the Lipschitz property of $h$,
\[
|h(\Sigma^{1/2}y)-h(\Sigma^{1/2}y')|
\le
L_h\|\Sigma^{1/2}(y-y')\|.
\]
Therefore
\[
\left|
\sigma(h(\Sigma^{1/2}y))
-
\sigma(h(\Sigma^{1/2}y'))
\right|^2
\le
\frac{L_h^2}{16}
\|\Sigma^{1/2}(y-y')\|^2.
\]
Hence
\[
\|F(y)-F(y')\|^2
\le
\left(1+\frac{L_h^2}{16}\right)
\|\Sigma^{1/2}(y-y')\|^2
\le
K_x^2 \|y-y'\|^2.
\]
This proves the claim.
\end{proof}

\begin{lemma}[Matrix concentration]\label{lem:matrix}
Let $\mathbf X$ be defined as above. Then $\mathbf X$ satisfies Gaussian
concentration with respect to the Frobenius norm with constant
\[
K_x
=
\|\Sigma^{1/2}\|
\sqrt{1+\frac{L_h^2}{16}}.
\]
That is, for every $1$-Lipschitz function
$\Phi:\mathbb R^{n\times(d+1)}\to\mathbb R$ with respect to the Frobenius norm,
\[
\mathbb P\left(|\Phi(\mathbf X)-\mathbb E\Phi(\mathbf X)|\ge t\right)
\le
2\exp\left(-\frac{t^2}{2K_x^2}\right).
\]
\end{lemma}

\begin{proof}
Let
\[
\mathbf Y=
(Y_1,\ldots,Y_n)^T
\in\mathbb R^{n\times d},
\]
where $Y_1,\ldots,Y_n$ are i.i.d. $N(0,I_d)$. Define
\[
\mathcal F(\mathbf Y)
=
(F(Y_1),\ldots,F(Y_n))^T.
\]
Then we observe that
$
\mathbf X=\mathcal F(\mathbf Y)
$
and for any $\mathbf y,\mathbf y'\in\mathbb R^{n\times d}$,
\[
\|\mathcal F(\mathbf y)-\mathcal F(\mathbf y')\|_F^2
=
\sum_{i=1}^n \|F(y_i)-F(y_i')\|^2
\le
K_x^2\sum_{i=1}^n \|y_i-y_i'\|^2
=
K_x^2\|\mathbf y-\mathbf y'\|_F^2.
\]
Thus $\mathcal F$ is $K_x$-Lipschitz.

Let $\Phi:\mathbb R^{n\times(d+1)}\to\mathbb R$ be $1$-Lipschitz with respect
to the Frobenius norm, and define
\[
\Psi:=\Phi\circ\mathcal F.
\]
Then $\Psi$ is $K_x$-Lipschitz on $\mathbb R^{n\times d}$. Since $\mathbf Y$ is
a standard Gaussian vector in $\mathbb R^{nd}$, Gaussian concentration gives
\[
\mathbb P\left(|\Psi(\mathbf Y)-\mathbb E\Psi(\mathbf Y)|\ge t\right)
\le
2\exp\left(-\frac{t^2}{2K_x^2}\right).
\]
Since $\Psi(\mathbf Y)=\Phi(\mathbf X)$, the result follows.
\end{proof}

\begin{corollary}\label{cor_conc}
If $\|\Sigma\|\le C_\Sigma$ and $L_h=O(1)$, then
$K_x=O(1)$. Consequently, $\mathbf X$ is
$\propto\mathcal E_2(1)$-concentrated with respect to the Frobenius norm.
\end{corollary}

\subsubsection{Uniform conditioning of the population covariance}
\label{sec:logistic_covariance_bounds}

We next verify that the covariance matrix of the centered logistic-link
design is uniformly well-conditioned under the assumptions imposed in the discussion after
Remark~\ref{rem:ass_1}. This result will also be used in the subsequent
verification of the Weak Tail Projection property.

Let
$
C_d=\Sigma_d^{1/2}Y_d, Y_d\sim N(0,I_d),$
and define
\[
G_d:=\sigma(h_d(C_d)),
\qquad
A_d:=G_d-\E[G_d],
\]
where $\sigma(t)=(1+e^{-t})^{-1}$. Let
\[
B_d:=(A_d, C^T_d)^T
\qquad
\Sigma_{B,d}:=\E[B_dB_d^\top].
\]
Assume that the covariance matrices of the Gaussian covariates are
uniformly well conditioned:
\begin{equation}
0<\underline{\sigma}_C I_d
\preceq
\Sigma_d
\preceq
\overline{\sigma}_C I_d
<\infty.
\label{eq:Sigma_C_uniform_condition}
\end{equation}

The following proposition shows that, under
\eqref{eq:Sigma_C_uniform_condition}, uniform conditioning of
$\Sigma_{B,d}$ is equivalent to requiring that the logistic exposure
does not become asymptotically degenerate.

\begin{proposition}[Characterization of population covariance conditioning]
\label{prop:logistic_cov_conditioning}
Under \eqref{eq:Sigma_C_uniform_condition}, the following conditions
are equivalent:
\begin{enumerate}[(i)]
\item
there exists $\underline{\sigma}_B>0$ such that
\[
\inf_d\lambda_{\min}(\Sigma_{B,d})
\ge
\underline{\sigma}_B;
\]

\item Let $
v_d:=\mathrm{Var}(G_d),
q_d:=\E[C_dA_d]$, then 
the Schur complement
$
\nu_d
:=
v_d-q_d^\top\Sigma_d^{-1}q_d
$
satisfies
$
\inf_d\nu_d>0;
$

\item
\[
\inf_d
\mathrm{Var}\!\left[
\sigma\!\left(h_d(\Sigma_d^{1/2}Y_d)\right)
\right]
>0.
\]
\end{enumerate}
Moreover, the upper spectral bound is automatic:
\[
\sup_d\lambda_{\max}(\Sigma_{B,d})
\le
\frac14+\overline{\sigma}_C<\infty.
\]
Consequently, under \eqref{eq:Sigma_C_uniform_condition},
\[
0<\underline{\sigma}_B
\le
\lambda_{\min}(\Sigma_{B,d})
\le
\lambda_{\max}(\Sigma_{B,d})
\le
\overline{\sigma}_B<\infty
\]
uniformly in $d$ if and only if
\begin{equation}
\inf_d
\mathrm{Var}\left[
\sigma\!\left(h_d(\Sigma_d^{1/2}Y_d)\right)
\right]
>0.
\label{eq:logistic_non_degeneracy}
\end{equation}
\end{proposition}

\begin{proof}
Since $\E[A_d]=\E[C_d]=0$, we have
\[
\Sigma_{B,d}
=
\begin{pmatrix}
v_d&q_d^\top\\
q_d&\Sigma_d
\end{pmatrix}.
\]

\medskip

\noindent
\textbf{Upper spectral bound.}
Since $0<G_d<1$,
$
v_d=\mathrm{Var}(G_d)\le\frac14.
$
For $x=(a,u^\top)^\top\in\mathbb R^{d+1}$,
\begin{align*}
x^\top\Sigma_{B,d}x
&=
\E[(aA_d+u^\top C_d)^2]\\
&\le
\left(
|a|\sqrt{v_d}
+
\sqrt{u^\top\Sigma_du}
\right)^2\\
&\le
(v_d+\overline{\sigma}_C)
(a^2+\|u\|^2).
\end{align*}
Therefore
\[
\lambda_{\max}(\Sigma_{B,d})
\le
\frac14+\overline{\sigma}_C.
\]

\medskip

\noindent
\textbf{Equivalence of (i) and (ii).}
The quantity
\[
\nu_d
=
v_d-q_d^\top\Sigma_d^{-1}q_d
\]
is the Schur complement of $\Sigma_d$ in $\Sigma_{B,d}$.
Equivalently,
\begin{equation}
\nu_d
=
\inf_{b\in\mathbb R^d}
\E[(A_d-b^\top C_d)^2].
\label{eq:nu_projection}
\end{equation}

If (i) holds, then
$
\|\Sigma_{B,d}^{-1}\|
\le
\underline{\sigma}_B^{-1}.
$
The block inversion formula gives
$
(\Sigma_{B,d}^{-1})_{11}=\nu_d^{-1},
$
and hence
$
\nu_d^{-1}
\le
\underline{\sigma}_B^{-1}.
$
Thus
$
\nu_d\ge\underline{\sigma}_B,
$
which proves (ii).

Conversely, suppose that
\[
\nu_d\ge\underline{\nu}>0.
\]
Define
$
b_d:=\Sigma_d^{-1}q_d.
$
Then
\[
v_d=\nu_d+b_d^\top\Sigma_db_d,
\]
and for $x=(a,u^\top)^\top$,
\begin{equation}
x^\top\Sigma_{B,d}x
=
\nu_da^2
+
(u+ab_d)^\top\Sigma_d(u+ab_d).
\label{eq:block_quad_logistic}
\end{equation}

Since
\[
b_d^\top\Sigma_db_d
=
q_d^\top\Sigma_d^{-1}q_d
\le
v_d
\le
\frac14,
\]
the lower bound in \eqref{eq:Sigma_C_uniform_condition} gives
\[
\|b_d\|
\le
\frac{1}{2\sqrt{\underline{\sigma}_C}}
=:K_b.
\]
Therefore
\[
a^2+\|u+ab_d\|^2
\ge
\frac{a^2+\|u\|^2}{(1+K_b)^2}.
\]
Combining this with \eqref{eq:block_quad_logistic},
\[
x^\top\Sigma_{B,d}x
\ge
\frac{
\min\{\underline{\nu},\underline{\sigma}_C\}
}{
(1+K_b)^2
}
\|x\|^2.
\]
Hence
\[
\inf_d\lambda_{\min}(\Sigma_{B,d})>0,
\]
which proves (i).

\medskip

\noindent
\textbf{Equivalence of (ii) and (iii).}
Since $\nu_d$ is a residual variance,
\[
0\le\nu_d\le v_d.
\]
Thus (ii) immediately implies (iii).

For the converse, suppose that
$
v_d\ge v_0>0
$
uniformly in $d$. Let
\[
L_d:=b_d^\top C_d,
\qquad
R_d:=A_d-L_d.
\]
Since $L_d$ is the $L^2$-projection of $A_d$ onto the linear span of
$C_d$,
\[
\E[L_dR_d]=0, \qquad
v_d=\E[L_d^2]+\E[R_d^2].
\]
Note that,
$
\E[R_d^2]=\nu_d.
$
We write
\[
s_d^2:=\E[L_d^2]=v_d-\nu_d.
\]
Because $C_d$ is Gaussian, $L_d$ is Gaussian:
\[
L_d\sim N(0,s_d^2).
\]

If $s_d^2\le v_0/2$, then
\[
\nu_d
=
v_d-s_d^2
\ge
v_0/2.
\]

Suppose instead that $s_d^2>v_0/2$. Since $0<G_d<1$, the centered
variable $A_d=G_d-\E G_d$ satisfies
$
|A_d|\le1.
$
Therefore
\[
|A_d-L_d|
\ge
(|L_d|-1)_+,
\]
and hence, for $Z\sim N(0,1)$,
\begin{align*}
\nu_d
&=
\E[(A_d-L_d)^2]\\
&\ge
\E[(|L_d|-1)_+^2]\\
&=
\E[(s_d|Z|-1)_+^2]\\
&\ge
\E\left[
\left(
\sqrt{\frac{v_0}{2}}\,|Z|-1
\right)_+^2
\right]
=:c(v_0)>0.
\end{align*}
Consequently,
\[
\nu_d
\ge
\min\left\{
\frac{v_0}{2},c(v_0)
\right\}>0
\]
uniformly in $d$. This proves (ii), and hence establishes the
equivalence.
\end{proof}

Next we give a counterexample where condition (iii) fails, so the lower spectral bound for $\Sigma_{B,d}$ does not hold. Take $\Sigma_d=I_d$ and
\[
h_d(c)=d+c_1.
\]
Then $\operatorname{Lip}(h_d)=1$ uniformly in $d$, but
\[
\sigma(h_d(C_d))=\sigma(d+Z)\xrightarrow{P}1,
\qquad Z\sim N(0,1),
\]
and consequently
\[
\mathrm{Var}\{\sigma(h_d(C_d))\}\to0.
\]
Thus the smallest eigenvalue of $\Sigma_{B,d}$ converges to zero.

\subsubsection{Weak Tail Projection property and the lower spectral edge}
\label{sec:logistic_wtp}

We next verify a condition that will be used to control the smallest
eigenvalue of the sample covariance matrix associated with the
logistic-link design. The relevant notion is the Weak Tail Projection
(WTP) property introduced in \cite{chafai_tikhomirov}.

Recall that, for each $d$, we consider
$
Y\sim N(0,I_d)$
and we denote 
\[
X^\circ
:=
F(Y)
=
\begin{pmatrix}
\sigma\!\bigl(h(\Sigma^{1/2}Y)\bigr)
-
\E\!\left[\sigma\!\bigl(h(\Sigma^{1/2}Y)\bigr)\right]\\
\Sigma^{1/2}Y
\end{pmatrix}.
\]
Let $
\Sigma_B
:=
\E[X^\circ(X^\circ)^\top].
$
Throughout this subsection we assume that, uniformly in $d$,
\begin{equation}
0<\underline{\sigma}
\le
\lambda_{\min}(\Sigma_B)
\le
\lambda_{\max}(\Sigma_B)
\le
\overline{\sigma}<\infty.
\label{eq:logistic_cov_bounds_wtp}
\end{equation}
We also retain from Corollary \ref{cor_conc} the assumptions
$
\|\Sigma\|\le C_\Sigma,
L_h\le L
$
for constants $C_\Sigma,L<\infty$ independent of $d$.

Define the whitened row
\begin{equation}
Z
:=
\Sigma_B^{-1/2}X^\circ
\in\mathbb R^{d+1}.
\label{eq:whitened_logistic_row}
\end{equation}
Then
\[
\E[Z]=0,
\qquad
\E[ZZ^\top]=I_{d+1},
\]
so $Z$ is isotropic.

We first record the concentration property inherited by $Z$. Recall
from Lemma~\ref{lem:matrix} and the preceding Lipschitz calculation
that
$
F:\mathbb R^d\to\mathbb R^{d+1}
$
is Lipschitz with constant
\[
K_x
=
\|\Sigma^{1/2}\|
\sqrt{1+\frac{L_h^2}{16}}.
\]
Centering by the deterministic vector $\E[F(Y)]$ does not change the
Lipschitz constant. Hence the map
\[
y\mapsto
\Sigma_B^{-1/2}
\bigl(F(y)-\E[F(Y)]\bigr)
\]
has Lipschitz constant at most
\[
K_z
:=
\|\Sigma_B^{-1/2}\|K_x
\le
\underline{\sigma}^{-1/2}
C_\Sigma^{1/2}
\sqrt{1+\frac{L^2}{16}}
=:K_0,
\]
where $K_0<\infty$ is independent of $d$.

Consequently, for every $K_0$-Lipschitz functional of $Z$, Gaussian
concentration can be applied through the underlying standard Gaussian
vector $Y$.

\medskip

\noindent
\textbf{Weak Tail Projection property.}
We recall Definition~1.1 of \cite{chafai_tikhomirov}. A sequence of
centered isotropic random vectors $Z_d\in\mathbb R^{d+1}$ is said to
satisfy the Weak Tail Projection property if the following two
conditions hold.

\begin{enumerate}
\item[(WTP-a)] The squared one-dimensional marginals are uniformly
integrable:
\[
\lim_{M\to\infty}
\sup_d\sup_{v\in S^d}
\E\left[
\langle Z_d,v\rangle^2
\mathbf 1_{\{\langle Z_d,v\rangle^2\ge M\}}
\right]
=0.
\]

\item[(WTP-b)] There exist functions
$
f:\mathbb N\to[0,1],
\qquad
g:\mathbb N\to\mathbb R_+
$
such that
$
f(r)\to0,
g(r)\to0
$
as $r\to\infty$, and such that for every $d$ and every nonzero
orthogonal projection $P$ on $\mathbb R^{d+1}$, writing
$r=\operatorname{rank}(P)$,
\[
\mathbb P\left(
\|PZ_d\|^2-r\ge f(r)r
\right)
\le g(r).
\]
\end{enumerate}

We now verify these two conditions for the centered and whitened
logistic-link design.

\begin{lemma}[WTP for the logistic-link design]
\label{lem:logistic_wtp}
Under \eqref{eq:logistic_cov_bounds_wtp}, together with
$\|\Sigma\|\le C_\Sigma$ and $L_h\le L$, the sequence of whitened
random vectors $Z$ defined in \eqref{eq:whitened_logistic_row}
satisfies the Weak Tail Projection property.
\end{lemma}

\begin{proof}
We verify the two parts separately.

\medskip

\noindent
\textbf{Verification of (WTP-a).}
Let $v\in S^d$. Since the map
\[
y\mapsto
\left\langle
v,
\Sigma_B^{-1/2}
\bigl(F(y)-\E[F(Y)]\bigr)
\right\rangle
\]
is $K_0$-Lipschitz and
$
\E[\langle v,Z\rangle]=0
$. Gaussian concentration gives
\begin{equation}
\mathbb P\left(
|\langle v,Z\rangle|\ge t
\right)
\le
2\exp\left(
-\frac{t^2}{2K_0^2}
\right),
\qquad t\ge0,
\label{eq:linear_marginal_tail}
\end{equation}
uniformly in $d$ and $v\in S^d$.

In particular, the fourth moments of all one-dimensional marginals
are uniformly bounded:
\begin{equation}
\sup_d\sup_{v\in S^d}
\E|\langle v,Z\rangle|^4
\le C K_0^4
<\infty
\label{eq:uniform_fourth_marginal}
\end{equation}
for a universal constant $C$. Therefore
\begin{equation*}
\E\left[
\langle v,Z\rangle^2
\mathbf 1_{\{\langle v,Z\rangle^2\ge M\}}
\right]\le
\frac{1}{M}
\E|\langle v,Z\rangle|^4
\le
\frac{CK_0^4}{M}.
\end{equation*}
Taking the supremum over $d$ and $v\in S^d$, and then letting
$M\to\infty$, yields
\[
\lim_{M\to\infty}
\sup_d\sup_{v\in S^d}
\E\left[
\langle v,Z\rangle^2
\mathbf 1_{\{\langle v,Z\rangle^2\ge M\}}
\right]
=0.
\]
Thus (WTP-a) holds.

\medskip

\noindent
\textbf{Verification of (WTP-b).}
Let $P$ be any nonzero orthogonal projection on $\mathbb R^{d+1}$
and set
$
r:=\operatorname{rank}(P).
$
Consider the function
\[
\Phi_P(y)
:=
\left\|
P\Sigma_B^{-1/2}
\bigl(F(y)-\E[F(Y)]\bigr)
\right\|.
\]
Since an orthogonal projection is a contraction and the Euclidean
norm is $1$-Lipschitz,
\[
\operatorname{Lip}(\Phi_P)\le K_0.
\]
Therefore Gaussian concentration gives
\begin{equation}
\mathbb P\left(
\Phi_P(Y)-\E[\Phi_P(Y)]\ge t
\right)
\le
\exp\left(
-\frac{t^2}{2K_0^2}
\right),
\qquad t\ge0.
\label{eq:projection_norm_concentration}
\end{equation}

Since $Z$ is isotropic,
\[
\E\|PZ\|^2
=
\tr\!\left(P\E[ZZ^\top]\right)
=
\tr(P)
=
r.
\]
Hence, by Jensen's inequality,
\begin{equation}
\E\|PZ\|
\le
\sqrt{\E\|PZ\|^2}
=
\sqrt r.
\label{eq:projection_mean_bound}
\end{equation}

For $r>1$, choose
$
f(r):=r^{-1/4}$.
Then $f(r)\in[0,1]$ and $f(r)\to0$. On the event
\[
\|PZ\|^2-r\ge f(r)r=r^{3/4},
\]
we have
\[
\|PZ\|
\ge
\sqrt{r+r^{3/4}}
=
\sqrt r\,\sqrt{1+r^{-1/4}}.
\]
Combining this with \eqref{eq:projection_mean_bound},
\begin{align*}
\|PZ\|-\E\|PZ\|
&\ge
\sqrt r
\left(
\sqrt{1+r^{-1/4}}-1
\right)\\
&=
\frac{r^{1/4}}
{\sqrt{1+r^{-1/4}}+1}.
\end{align*}
Using the bound $
\sqrt{1+r^{-1/4}}+1
\le
1+\sqrt2,
$ for $r\ge1$,
we have 
\begin{equation}
\|PZ\|-\E\|PZ\|
\ge
\frac{r^{1/4}}{1+\sqrt2}.
\label{eq:projection_deviation_lower}
\end{equation}

Applying \eqref{eq:projection_norm_concentration} gives
\begin{align*}
\mathbb P\left(
\|PZ\|^2-r\ge r^{3/4}
\right)
&\le
\mathbb P\left(
\|PZ\|-\E\|PZ\|
\ge
\frac{r^{1/4}}{1+\sqrt2}
\right)\\
&\le
\exp\left(
-\frac{\sqrt r}
{2K_0^2(1+\sqrt2)^2}
\right).
\end{align*}
Thus (WTP-b) holds with, for example,
\[
f(r)=r^{-1/4},
\qquad
g(r)
=
\exp\left(
-\frac{\sqrt r}
{2K_0^2(1+\sqrt2)^2}
\right).
\]
Clearly $f(r)\to0$ and $g(r)\to0$, so the sequence of whitened logistic-link rows satisfies the Weak Tail
Projection property.
\end{proof}

We next use the WTP property to establish the lower spectral edge
needed in the inverse-Gram analysis.

\begin{lemma}[Lower spectral edge for the logistic-link design]
\label{lem:logistic_hard_edge}
Let $Z_1,\ldots,Z_n$ be i.i.d.\ copies of the whitened vector $Z$
defined in \eqref{eq:whitened_logistic_row}, and let
\[
\mathbf Z
=
(
Z_1,
\ldots,
Z_n
)^{\top}
\in\mathbb R^{n\times(d+1)}.
\]
Assume
$
\frac{d+1}{n}\rightarrow c\in(0,1).$
Then, for every
\[
0<\varepsilon<1-\sqrt c,
\]
\begin{equation}
\mathbb P\left(
\lambda_{\min}\left(
\frac1n\mathbf Z^\top\mathbf Z
\right)
\ge
(1-\sqrt c-\varepsilon)^2
\right)
\rightarrow1.
\label{eq:whitened_hard_edge}
\end{equation}

Consequently, if $\mathbf B$ denotes the matrix whose rows are the
centered logistic-link vectors $X_i^\circ$, then
\begin{equation}
\mathbb P\left(
\lambda_{\min}\left(
\frac1n\mathbf B^\top\mathbf B
\right)
\ge
\underline\sigma
(1-\sqrt c-\varepsilon)^2
\right)
\rightarrow1.
\label{eq:logistic_hard_edge}
\end{equation}
In particular, there exists a deterministic constant $a>0$ such that
\[
\mathbb P\left(
\lambda_{\min}\left(
\frac1n\mathbf B^\top\mathbf B
\right)\ge a
\right)
\rightarrow1.
\]
\end{lemma}

\begin{proof}
By Lemma~\ref{lem:logistic_wtp}, the centered isotropic vectors $Z$
satisfy the WTP property. Let
$
q:=d+1
$
and define
$
A_q:=\mathbf Z^\top\mathbf Z
$.
Theorem~1.5 of \cite{chafai_tikhomirov} applies because
$q/n\to c<1$ and yields
\begin{equation}
\liminf_{d\to\infty}
\frac{
\E[\lambda_{\min}(A_q)]
}{
(\sqrt n-\sqrt q)^2
}
\ge1.
\label{eq:CT_expected_edge}
\end{equation}
Since
$
\lambda_{\min}(A_q)
=
s_{\min}(\mathbf Z)^2,
$
this is equivalently
\begin{equation}
\liminf_{d\to\infty}
\frac{
\E[s_{\min}(\mathbf Z)^2]
}{
(\sqrt n-\sqrt q)^2
}
\ge1.
\label{eq:CT_expected_smin}
\end{equation}

We now use the Gaussian-Lipschitz structure of the logistic design
to convert this expected lower-edge bound to a high-probability
statement.

Let
\[
\mathbf Y
=
(Y_1,\ldots,Y_n)^\top
\in\mathbb R^{n\times d},
\]
where the $Y_i$ are i.i.d.\ $N(0,I_d)$. The map
\[
\mathbf Y\mapsto\mathbf Z
\]
is $K_0$-Lipschitz with respect to the Frobenius norm. Indeed, by the
row-wise Lipschitz estimate,
\begin{align*}
\|\mathbf Z(\mathbf Y)-\mathbf Z(\mathbf Y')\|_F^2
&=
\sum_{i=1}^n
\|Z(Y_i)-Z(Y_i')\|^2\\
&\le
K_0^2
\sum_{i=1}^n
\|Y_i-Y_i'\|^2\\
&=
K_0^2
\|\mathbf Y-\mathbf Y'\|_F^2.
\end{align*}

Moreover, for any two matrices $M,N$ of the same dimensions,
\[
|s_{\min}(M)-s_{\min}(N)|
\le
\|M-N\|
\le
\|M-N\|_F.
\]
Therefore
$
\mathbf Y\mapsto s_{\min}(\mathbf Z(\mathbf Y))
$
is $K_0$-Lipschitz. Gaussian concentration thus yields
\begin{equation}
\mathbb P\left(
\left|
s_{\min}(\mathbf Z)
-
\E s_{\min}(\mathbf Z)
\right|
\ge t
\right)
\le
2\exp\left(
-\frac{t^2}{2K_0^2}
\right).
\label{eq:smin_concentration}
\end{equation}

Integrating the tail bound gives
\begin{align}
\mathrm{Var}(s_{\min}(\mathbf Z))
&\le
\E\left[
\left(
s_{\min}(\mathbf Z)
-
\E s_{\min}(\mathbf Z)
\right)^2
\right]\nonumber\\
&\le
4K_0^2.
\label{eq:smin_variance_bound}
\end{align}

Writing
$
\mu_d:=\E s_{\min}(\mathbf Z),
$
we therefore obtain
\[
\mu_d^2
=
\E[s_{\min}(\mathbf Z)^2]
-
\mathrm{Var}(s_{\min}(\mathbf Z))
\ge
\E[s_{\min}(\mathbf Z)^2]-4K_0^2.
\]
Combining this with \eqref{eq:CT_expected_smin} gives
\[
\liminf_{d\to\infty}
\frac{\mu_d}{\sqrt n}
\ge
1-\sqrt c.
\]

For all sufficiently large $d$ and $
0<\varepsilon<1-\sqrt c,
$
\[
\mu_d
\ge
\sqrt n
\left(
1-\sqrt c-\frac{\varepsilon}{2}
\right).
\]
Hence
\begin{equation*}
\mathbb P\left(
\frac{s_{\min}(\mathbf Z)}{\sqrt n}
<
1-\sqrt c-\varepsilon
\right)\le
\mathbb P\left(
s_{\min}(\mathbf Z)-\mu_d
<
-\frac{\varepsilon}{2}\sqrt n
\right)\le
2\exp\left(
-\frac{n\varepsilon^2}{8K_0^2}
\right)
\rightarrow0,
\end{equation*}
by \eqref{eq:smin_concentration}. Since
\[
\lambda_{\min}\left(
\frac1n\mathbf Z^\top\mathbf Z
\right)
=
\frac{s_{\min}(\mathbf Z)^2}{n},
\]
this proves \eqref{eq:whitened_hard_edge}.

Finally, since
$
X_i^\circ=\Sigma_B^{1/2}Z_i,
$
we have
\[
\frac1n\mathbf B^\top\mathbf B
=
\Sigma_B^{1/2}
\left(
\frac1n\mathbf Z^\top\mathbf Z
\right)
\Sigma_B^{1/2}.
\]
Consequently,
\begin{align*}
\lambda_{\min}\left(
\frac1n\mathbf B^\top\mathbf B
\right)
&\ge
\lambda_{\min}(\Sigma_B)
\lambda_{\min}\left(
\frac1n\mathbf Z^\top\mathbf Z
\right)\\
&\ge
\underline\sigma
\lambda_{\min}\left(
\frac1n\mathbf Z^\top\mathbf Z
\right).
\end{align*}
Combining this inequality with
\eqref{eq:whitened_hard_edge} proves
\eqref{eq:logistic_hard_edge}.

For example, taking
$
\varepsilon
=
\frac{1-\sqrt c}{2}
$
gives the deterministic choice
\[
a
=
\frac{\underline\sigma}{4}
(1-\sqrt c)^2>0
\]
and hence
\[
\mathbb P\left(
\lambda_{\min}\left(
\frac1n\mathbf B^\top\mathbf B
\right)\ge a
\right)
\rightarrow1.
\]
\end{proof}

\section{The degenerate case for the pathwise mediation effect}
\label{sec:pathwise_nonregular}
The point
$
(\beta_{A,j},\theta_{M,j})=(0,0)
$
requires separate treatment. Indeed, the pathwise mediation functional is
$
g(b,t)=bt,
$
whose gradient is
\[
\nabla g(b,t)=(t,b).
\]
Thus, whenever
\[
\beta_{A,j}^2+\theta_{M,j}^2>0,
\]
the functional is locally regular and the first-order fluctuations of
$\hat{\beta}_{A,j}$ and $\hat{\theta}_{M,j}$ yield the
$\sqrt m$-Gaussian limit in Theorem~\ref{thm_plugin_path}. This includes
the two regular null cases in which exactly one of
$\beta_{A,j}$ and $\theta_{M,j}$ is zero.

In contrast, when
\[
\beta_{A,j}=\theta_{M,j}=0,
\]
the first derivative of $g$ vanishes. Consequently, the
$\sqrt m$-scaled pathwise estimator is asymptotically degenerate,
and the first nontrivial fluctuation occurs at the $m$ scale
(equivalently, at the $n=m/2$ scale). Using the fold normalization $n=m/2$ and the corresponding coordinatewise
limits in our model, we have 
\[
\sqrt n\,\hat{\beta}_{A,j}\xrightarrow{d}U,
\qquad
\sqrt n\,\hat{\theta}_{M,j}\xrightarrow{d}V,
\]
where $U$ and $V$ are centered nondegenerate Gaussian random variables.
For the sample-split estimator the two variables are independent,
because the two component estimators are constructed on independent
folds. Hence
\[
n\hat{\eta}_j
=
\bigl(\sqrt n\,\hat{\beta}_{A,j}\bigr)
\bigl(\sqrt n\,\hat{\theta}_{M,j}\bigr)
\xrightarrow{d}UV.
\]
Equivalently, $m\hat{\eta}_j\xrightarrow{d}2UV$. The product $UV$ is not Gaussian. Thus the failure of asymptotic
normality at $(0,0)$ is intrinsic to the product structure of the
estimand rather than to the particular debiasing procedure.

Accordingly, Gaussian inference for the pathwise mediation effect developed in this paper is stated under the
regularity condition
\[
\beta_{A,j}^2+\theta_{M,j}^2>0.
\]
Inference at $(\beta_{A,j},\theta_{M,j})=(0,0)$ is a distinct
second-order problem and is not considered here.

\section{Finite-dimensional and limiting bias corrections in the Gaussian experiments}
\label{sec:centering_correction}

The estimators in the main theorems use the finite-dimensional
bias-correction factor
\[
\mathfrak c_{1,n}
=
\{1-\kappa_1 f_{1,n}(\kappa_1)\}^{-1},
\]
whereas the Gaussian experiments use
\[
\mathfrak c_1
=
\{1-\kappa_1 f_1(\kappa_1,c_1)\}^{-1}.
\]
This appendix justifies that substitution for the proportional
dimension choices used in the experiments. The two factors need not
be equal at finite dimensions. We show instead that their difference
is $O(n^{-1})$, which makes the resulting changes to the global
and regular pathwise estimators negligible at the $\sqrt m$ scale,
where $m=2n$. 

\subsection{Notation and the Gaussian finite-dimensional equation}

Throughout this appendix, $\kappa_1>0$ is fixed, and
\[
r_n:=n-d-1,\qquad
c_{1,n}:=\frac{p}{r_n},\qquad
c_{2,n}:=\frac{d+1}{n}.
\]
We assume $c_{1,n}\to c_1\in(0,\infty)$ and
$c_{2,n}\to c_2\in(0,1)$. In particular,
$r_n\asymp p\asymp n$. To distinguish the two sources of
finite-dimensional discrepancy, set
\begin{equation}
f_n:=f_{1,n}(\kappa_1),\qquad
f_n^\star:=f_1(\kappa_1,c_{1,n}),\qquad
f_\infty:=f_1(\kappa_1,c_1).
\label{eq:cc:three_centers}
\end{equation}
Thus $f_n$ is the finite-system center used in the paper,
$f_n^\star$ is the Mar\v{c}enko--Pastur center evaluated at the
actual residual aspect ratio, and $f_\infty$ uses the limiting
aspect ratio. 

Recall the Gaussian specialization of the finite system in the proof
of Lemma~\ref{MP_fluc}, based on \cite{yin2023limiting}.
Here that lemma is used with $\mathfrak n=r_n=n-d-1$,
$\mathfrak p=p$, and $\kappa=\kappa_1$, as specified after
\eqref{eq:MP_fluc_finite_center}. Thus $f_n$ is the same
finite-system center as in Lemma~\ref{MP_fluc}, evaluated at
the residual sample size.
A standard Gaussian vector in $\mathbb R^p$ has squared radius
$\rho_p^2\sim\chi_p^2$, with
\[
\mathbb E[\rho_p^2]=p,\qquad
\mathbb E[\rho_p^4]=p^2+2p.
\]
Accordingly, the radial distribution entering that system is the law of
\[
T_p:=\frac{\rho_p^2}{\sqrt{p^2+2p}},
\qquad
a_p:=\mathbb E[T_p]
=\left(1+\frac2p\right)^{-1/2}.
\]
This is a population distribution, not the empirical distribution of
observed radii. All expectations involving $T_p$ below are with
respect to this distribution; in particular, $f_n$ is deterministic.

Since the population covariance of the mediator errors is $I_p$,
the finite equations at $z=-\kappa_1$ reduce to
\begin{equation}
\begin{aligned}
\kappa_1 g^0_{1,r_n}(-\kappa_1)
&=\frac{c_{1,n}}{1+g^0_{2,r_n}(-\kappa_1)},\\
\kappa_1 g^0_{2,r_n}(-\kappa_1)
&=\mathbb E\left[
\frac{T_p}{1+g^0_{1,r_n}(-\kappa_1)T_p}
\right].
\end{aligned}
\label{eq:cc:finite_system}
\end{equation}
By the definition of $f_{1,n}(\kappa_1)$,
\[
f_n=\frac{1}{\kappa_1\{1+g^0_{2,r_n}(-\kappa_1)\}},
\qquad
g^0_{1,r_n}(-\kappa_1)=c_{1,n}f_n.
\]
Eliminating the two auxiliary functions gives
\begin{equation}
\frac1{f_n}
=
\kappa_1+
\mathbb E\left[\frac{T_p}{1+c_{1,n}f_nT_p}\right].
\label{eq:cc:scalar_finite}
\end{equation}
Equivalently, $f_n$ solves $H_n(f_n)=0$, where
\begin{equation}
H_n(x)
:=
\kappa_1x+
\mathbb E\left[\frac{xT_p}{1+c_{1,n}xT_p}\right]-1,
\qquad x\ge0.
\label{eq:cc:Hn}
\end{equation}
Differentiating under the expectation, which is justified by
the integrable bound $T_p$, yields
\[
H_n'(x)
=
\kappa_1+
\mathbb E\left[\frac{T_p}{(1+c_{1,n}xT_p)^2}\right]
\ge\kappa_1.
\]
Also $H_n(0)=-1$ and $H_n(1/\kappa_1)>0$.
Consequently, the positive solution is unique and satisfies
\begin{equation}
0<f_n<\kappa_1^{-1}.
\label{eq:cc:root_interval}
\end{equation}
In particular, the finite-dimensional correction factor is well defined.

For comparison, the function $f_1$ in \eqref{fixed_pt_1}
is the unique positive solution of
\begin{equation}
\kappa_1f+\frac{f}{1+cf}=1,
\end{equation}
or, equivalently,
\begin{equation}
\kappa_1cf^2+(\kappa_1+1-c)f-1=0.
\label{eq:cc:MP_equation}
\end{equation}
Thus $f_n^\star$ solves the same scalar equation as $f_n$
with $T_p$ replaced by the constant $1$, while retaining
the actual aspect ratio $c_{1,n}$.

\subsection{Quantitative comparison of the centering constants}

\begin{proposition}[Gaussian finite-centering comparison]
\label{prop:cc:comparison}
For fixed $\kappa_1>0$, under the aspect-ratio assumptions above,
\begin{equation}
f_n-f_n^\star=O(p^{-1}).
\label{eq:cc:finite_to_MP}
\end{equation}
Moreover,
\begin{equation}
|f_n-f_\infty|
=
O\left(p^{-1}+|c_{1,n}-c_1|\right),
\label{eq:cc:finite_to_limit}
\end{equation}
and
\begin{equation}
|\mathfrak c_{1,n}-\mathfrak c_1|
=
O\left(p^{-1}+|c_{1,n}-c_1|\right).
\label{eq:cc:factor_comparison}
\end{equation}
\end{proposition}

\begin{proof}
We first control the radial perturbation, then the aspect-ratio
perturbation, and finally the bias-correction factors.

\paragraph{The radial perturbation.}
The normalization of $T_p$ gives
\[
\mathbb E[T_p^2]=1,\qquad
\operatorname{Var}(T_p)=1-a_p^2=\frac{2}{p+2}.
\]
Also,
\[
0<1-a_p\le\frac1p,\qquad
\mathbb E[(T_p-1)^2]=2(1-a_p)\le\frac2p.
\]
Choose $C_c<\infty$ such that $c_{1,n}\le C_c$ for all
sufficiently large $n$. For these $n$, define
\[
h_n(t):=\frac{f_n^\star t}{1+c_{1,n}f_n^\star t},
\qquad t\ge0,
\]
whose derivatives satisfy
\[
h_n'(t)
=
\frac{f_n^\star}{(1+c_{1,n}f_n^\star t)^2},
\qquad
h_n''(t)
=
-\frac{2c_{1,n}(f_n^\star)^2}
{(1+c_{1,n}f_n^\star t)^3}.
\]
Since $0<f_n^\star<\kappa_1^{-1}$,  uniformly over $t\ge0$, we have
\[
|h_n'(t)|\le\kappa_1^{-1},
\qquad
|h_n''(t)|\le 2C_c\kappa_1^{-2}.
\]
Taylor's theorem around $t=1$ therefore implies
\[
|h_n(t)-h_n(1)-h_n'(1)(t-1)|
\le C_c\kappa_1^{-2}(t-1)^2.
\]
Taking expectations and using the preceding moment bounds gives
\begin{align}
|\mathbb E[h_n(T_p)]-h_n(1)|
&\le
\kappa_1^{-1}|\mathbb E[T_p]-1|
+
C_c\kappa_1^{-2}\mathbb E[(T_p-1)^2]\nonumber\\
&\le
\frac{\kappa_1^{-1}+2C_c\kappa_1^{-2}}{p}.
\label{eq:cc:Taylor_bound}
\end{align}
Because $f_n^\star$ satisfies \eqref{eq:cc:MP_equation},
\[
H_n(f_n^\star)
=
\mathbb E[h_n(T_p)]-h_n(1).
\]
The derivative bound $H_n'\ge\kappa_1$ and
$H_n(f_n)=0$ now yield
\begin{equation}
|f_n-f_n^\star|
\le
\frac{|H_n(f_n^\star)|}{\kappa_1}
\le
\frac{\kappa_1^{-2}+2C_c\kappa_1^{-3}}{p}.
\label{eq:cc:explicit_root_bound}
\end{equation}
This proves \eqref{eq:cc:finite_to_MP}.

\paragraph{The aspect-ratio perturbation.}
Implicit differentiation of
$\kappa_1f_1+f_1/(1+cf_1)=1$ gives
\begin{equation}
\partial_c f_1(\kappa_1,c)
=
\frac{f_1(\kappa_1,c)^2}
{1+\kappa_1\{1+cf_1(\kappa_1,c)\}^2}.
\label{eq:cc:c_derivative}
\end{equation}
This derivative is positive and at most $\kappa_1^{-2}$.
The mean value theorem therefore gives
\[
|f_n^\star-f_\infty|
\le\kappa_1^{-2}|c_{1,n}-c_1|.
\]
Together with \eqref{eq:cc:explicit_root_bound},
this proves \eqref{eq:cc:finite_to_limit}.

\paragraph{The bias-correction factors.}
The limiting equation implies
\[
1-\kappa_1f_\infty
=
\frac{f_\infty}{1+c_1f_\infty}>0.
\]
Since $f_n\to f_\infty$, the quantities
$1-\kappa_1f_n$ are bounded away from zero for all
sufficiently large $n$. The exact identity
\[
\mathfrak c_{1,n}-\mathfrak c_1
=
\frac{\kappa_1(f_n-f_\infty)}
{(1-\kappa_1f_n)(1-\kappa_1f_\infty)}
\]
then proves \eqref{eq:cc:factor_comparison}.
\end{proof}

\begin{remark}[Asymptotic equivalence is not finite-dimensional equality]
\label{rem:cc:not_exact}
Even if $c_{1,n}=c_1$ exactly, the finite-system center
need not equal $f_1(\kappa_1,c_1)$.
Indeed, $h_n$ in the preceding proof is increasing and
concave, and $\mathbb E[T_p]=a_p<1$. Thus
\[
\mathbb E[h_n(T_p)]
\le h_n(a_p)<h_n(1).
\]
Consequently $H_n(f_n^\star)<0$, and strict monotonicity
of $H_n$ gives $f_n>f_n^\star$.
What matters for the central limit theorems is the
quantitative bound \eqref{eq:cc:finite_to_MP}, not equality
of the two finite-dimensional formulas.
\end{remark}

\paragraph{Covariance-normalization convention.}
The finite system recalled above is associated with the
normalized Gaussian covariance matrix $a_p S$, whereas
the residualized regression uses $S$ itself.
This distinction is already accounted for in the proof of
Lemma~\ref{MP_fluc}; it also admits a direct bound here.
For any positive semidefinite matrix $S$, we have
\begin{align*}
\left\|(a_pS+\kappa_1I)^{-1}
-(S+\kappa_1I)^{-1}\right\|
&\le
|a_p-1|
\sup_{x\ge0}
\frac{x}{(a_px+\kappa_1)(x+\kappa_1)}\\
&\le
\frac{|a_p-1|}{a_p\kappa_1}
=
O(p^{-1}).
\end{align*}
Multiplying this bound by $\sqrt p$ gives $o(1)$.
Thus the normalization convention introduces no additional
first-order fluctuation.

\subsection{Exact proportional dimensions}

Suppose first that the residual aspect ratios are exactly
proportional:
\begin{equation}
d+1=nc_2,\qquad p=nc_1(1-c_2),
\label{eq:cc:exact_dimensions}
\end{equation}
along sample sizes for which these quantities are integers.
Then
\[
r_n=n(1-c_2),\qquad c_{1,n}=c_1.
\]
Proposition~\ref{prop:cc:comparison} and $p\asymp n$ give
\begin{equation}
f_{1,n}(\kappa_1)-f_1(\kappa_1,c_1)=O(n^{-1}),
\qquad
\mathfrak c_{1,n}-\mathfrak c_1=O(n^{-1}).
\label{eq:cc:rate_exact}
\end{equation}

The same conclusion holds for the rounded dimensions used in
Section~\ref{sim_frame}:
\begin{equation}
p=\lfloor nc_1(1-c_2)\rfloor,\qquad
d=\lfloor nc_2\rfloor-1.
\label{eq:cc:rounded_dimensions}
\end{equation}
Since $p=nc_1(1-c_2)+O(1)$ and $r_n=n(1-c_2)+O(1)$,
we have $c_{1,n}-c_1=O(n^{-1})$. The bounds in
Proposition~\ref{prop:cc:comparison} therefore imply
\eqref{eq:cc:rate_exact} for these rounded dimensions as well.

For deterministic vectors $u_n,v_n$ of bounded Euclidean norms,
\eqref{eq:cc:rate_exact} also implies
\[
\sqrt p\,
\left|
\{f_{1,n}(\kappa_1)-f_1(\kappa_1,c_1)\}
u_n^\top v_n
\right|
=O(n^{-1/2})\longrightarrow0.
\]
Thus the finite and limiting centers are interchangeable
inside the bilinear resolvent CLT for these dimension sequences.

\subsection{Effect on the mediation estimators}

We now show that the scalar comparison gives the required
comparison of the actual estimators.
For $s\in\{a,b\}$, let $\breve\theta_{M,s}$ be the
uncorrected partial-ridge estimator on fold $s$, and let
$\widehat\beta_{A,s}$ be the OLS estimator on that fold.
Define
\[
\widehat\theta_{M,s}^{\mathrm{fin}}
:=
\mathfrak c_{1,n}\breve\theta_{M,s},
\qquad
\widehat\theta_{M,s}^{\mathrm{lim}}
:=
\mathfrak c_1\breve\theta_{M,s}.
\]
The superscripts distinguish only the correction factors;
the data, fold assignments, and ridge penalties are identical.

For $u\in\{\mathrm{fin},\mathrm{lim}\}$ and $s\ne t$, set
\begin{align*}
\widehat\gamma_{(s,t)}^u
&:=
\widehat\beta_{A,t}^{\top}\widehat\theta_{M,s}^u,\\
\widehat\eta_{j,(s,t)}^u
&:=
\widehat\beta_{A,t,j}\widehat\theta_{M,s,j}^u.
\end{align*}
The cross-fitted versions are
\[
\widehat\gamma_{cf}^u
=
\frac{\widehat\gamma_{(a,b)}^u+
\widehat\gamma_{(b,a)}^u}{2},
\qquad
\widehat\eta_{j,cf}^u
=
\frac{\widehat\eta_{j,(a,b)}^u+
\widehat\eta_{j,(b,a)}^u}{2}.
\]
As in the main analysis, inverse-based identities are understood
on the full-column-rank events, with the same extension convention
for both versions. Those events have probability tending to one.

\begin{corollary}[Equivalence at the central-limit scale]
\label{cor:cc:estimators}
Under the Gaussian model and the assumptions of the main global
mediation theorems, suppose that
$c_{1,n}-c_1=O(n^{-1})$. Then, for either fold assignment,
\begin{align}
\sqrt m\,
(\widehat\gamma_{(s,t)}^{\mathrm{lim}}
-\widehat\gamma_{(s,t)}^{\mathrm{fin}})
&\xrightarrow{P}0,\nonumber\\
\sqrt m\,
(\widehat\gamma_{cf}^{\mathrm{lim}}
-\widehat\gamma_{cf}^{\mathrm{fin}})
&\xrightarrow{P}0.
\label{eq:cc:global_equivalence}
\end{align}
For each fixed $j$, the analogous comparisons hold:
\begin{align}
\sqrt m\,
(\widehat\eta_{j,(s,t)}^{\mathrm{lim}}
-\widehat\eta_{j,(s,t)}^{\mathrm{fin}})
&\xrightarrow{P}0,\nonumber\\
\sqrt m\,
(\widehat\eta_{j,cf}^{\mathrm{lim}}
-\widehat\eta_{j,cf}^{\mathrm{fin}})
&\xrightarrow{P}0.
\label{eq:cc:path_equivalence}
\end{align}
Under Assumption~\ref{ass:path}, the limiting-correction pathwise
estimators therefore have the same nondegenerate limits as the
finite-correction estimators.
\end{corollary}

\begin{proof}
Lemmas~\ref{ba2_lemma} and~\ref{tm2_lemma}, applied to each fold,
give
\begin{equation}
\|\widehat\beta_{A,s}\|=O_P(1),
\qquad
\|\breve\theta_{M,s}\|=O_P(1).
\label{eq:cc:norm_bounds}
\end{equation}

By Cauchy--Schwarz,
\[
|\widehat\beta_{A,t}^{\top}\breve\theta_{M,s}|
\le
\|\widehat\beta_{A,t}\|\,
\|\breve\theta_{M,s}\|
=O_P(1).
\]
Likewise,
\[
|\widehat\beta_{A,t,j}\breve\theta_{M,s,j}|
\le
\|\widehat\beta_{A,t}\|\,
\|\breve\theta_{M,s}\|
=O_P(1).
\]
The exact differences are
\begin{align*}
\widehat\gamma_{(s,t)}^{\mathrm{lim}}
-\widehat\gamma_{(s,t)}^{\mathrm{fin}}
&=
(\mathfrak c_1-\mathfrak c_{1,n})
\widehat\beta_{A,t}^{\top}\breve\theta_{M,s},\\
\widehat\eta_{j,(s,t)}^{\mathrm{lim}}
-\widehat\eta_{j,(s,t)}^{\mathrm{fin}}
&=
(\mathfrak c_1-\mathfrak c_{1,n})
\widehat\beta_{A,t,j}\breve\theta_{M,s,j}.
\end{align*}
Since $\mathfrak c_1-\mathfrak c_{1,n}=O(n^{-1})$,
both differences are $O_P(n^{-1})$.
Multiplication by $\sqrt m=\sqrt{2n}$ gives $o_P(1)$.
Averaging the two fold assignments proves the cross-fitted
comparisons. This last step requires no independence between
the two cross-fitted components.
\end{proof}

Write $\widehat\gamma^u=\widehat\gamma_{(a,b)}^u$ and
$\widehat\eta_j^u=\widehat\eta_{j,(a,b)}^u$.
The main theorems and Slutsky's theorem now give
\begin{align*}
\sqrt m(\widehat\gamma^{\mathrm{lim}}-\gamma)
&\xrightarrow{d}N(0,2\sigma^2),&
\sqrt m(\widehat\gamma_{cf}^{\mathrm{lim}}-\gamma)
&\xrightarrow{d}N(0,\sigma^2),\\
\sqrt m(\widehat\eta_j^{\mathrm{lim}}-\eta_j)
&\xrightarrow{d}N(0,2\sigma_j^2),&
\sqrt m(\widehat\eta_{j,cf}^{\mathrm{lim}}-\eta_j)
&\xrightarrow{d}N(0,\sigma_j^2).
\end{align*}
The pathwise statements use Assumption~\ref{ass:path}.
In every case, centering remains at the actual
sample-size-dependent target $\gamma=\beta_A^\top\theta_M$
or $\eta_j=\beta_{A,j}\theta_{M,j}$, not at its limiting value.

\subsection{Variance estimation and studentization}

Let $\widehat\sigma_{\mathrm{fin}}^2$ and
$\widehat\sigma_{\mathrm{lim}}^2$ denote the variance estimators
in \eqref{lim_var_est} evaluated with the corresponding
versions of the mediation estimator. Define
$\widehat\sigma_{j,\mathrm{fin}}^2$ and
$\widehat\sigma_{j,\mathrm{lim}}^2$ similarly using
\eqref{est_var_sigma_j}. All other computational choices
are held fixed.

The estimators $\widehat{\mathfrak t}_M$,
$\widehat{\mathfrak b}_A$, and
$\widehat{\mathfrak r}$ are unchanged: they depend on the
uncorrected ridge estimator, the OLS estimator, and the
design, rather than on the scalar correction applied to
$\breve\theta_M$.
In the global variance estimator, only the squared mediation
estimate changes. Consequently,
\[
\widehat\sigma_{\mathrm{lim}}^2-
\widehat\sigma_{\mathrm{fin}}^2
=
\varphi_4
\left\{
(\widehat\gamma^{\mathrm{lim}})^2
-
(\widehat\gamma^{\mathrm{fin}})^2
\right\}
=
O_P(n^{-1}),
\]
where the arguments $(\kappa_1,c_1,c_2)$ of $\varphi_4$
are suppressed. Here both mediation estimates are $O_P(1)$
and their difference is $O_P(n^{-1})$.

For the pathwise estimator,
\[
\widehat\theta_{M,a,j}^{\mathrm{lim}}
-
\widehat\theta_{M,a,j}^{\mathrm{fin}}
=
(\mathfrak c_1-\mathfrak c_{1,n})
\breve\theta_{M,a,j}
=
O_P(n^{-1}).
\]
The variance formula therefore gives
\begin{align*}
\widehat\sigma_{j,\mathrm{lim}}^2-
\widehat\sigma_{j,\mathrm{fin}}^2
={}&
\left(
\widehat{\mathfrak r}
+
\psi_2\widehat\beta_{A,b,j}^{\,2}
\right)\\
&\quad\times
\left\{
(\widehat\theta_{M,a,j}^{\mathrm{lim}})^2
-
(\widehat\theta_{M,a,j}^{\mathrm{fin}})^2
\right\}
=
O_P(n^{-1}).
\end{align*}

The coefficient multiplying the squared-coordinate difference is
$O_P(1)$. The same comparisons hold when $c_{1,n}$ and $c_{2,n}$
are substituted for the limiting aspect ratios in both variance
estimators, as in the implementation convention following
Theorem~\ref{fluc_est_var}: the corresponding coefficient functions are bounded near the stipulated
interior
limits. These substitutions preserve consistency by continuity;
no root-$n$ rate is needed for variance consistency itself.
It follows from Theorems~\ref{fluc_est_var}
and~\ref{fluc_est_var_path} that
\[
\widehat\sigma_{\mathrm{lim}}^2\xrightarrow{P}\sigma^2,
\qquad
\widehat\sigma_{j,\mathrm{lim}}^2\xrightarrow{P}\sigma_j^2.
\]
With the positive-variance convention used in the paper,
the standardized statistics satisfy
\begin{align*}
\frac{\sqrt m(\widehat\gamma^{\mathrm{lim}}-\gamma)}
{\sqrt2\,\widehat\sigma_{\mathrm{lim}}}
&\xrightarrow{d}N(0,1),&
\frac{\sqrt m(\widehat\gamma_{cf}^{\mathrm{lim}}-\gamma)}
{\widehat\sigma_{\mathrm{lim}}}
&\xrightarrow{d}N(0,1),\\
\frac{\sqrt m(\widehat\eta_j^{\mathrm{lim}}-\eta_j)}
{\sqrt2\,\widehat\sigma_{j,\mathrm{lim}}}
&\xrightarrow{d}N(0,1),&
\frac{\sqrt m(\widehat\eta_{j,cf}^{\mathrm{lim}}-\eta_j)}
{\widehat\sigma_{j,\mathrm{lim}}}
&\xrightarrow{d}N(0,1).
\end{align*}
Thus the limiting correction used in the Gaussian QQ-plots
preserves both the central limit theorems and their feasible
standardizations. This is an asymptotic statement; it does
not assert equality of the finite-sample statistics.

\subsection{Scope of the replacement}

Exact proportionality is sufficient, but not necessary.
Because $p\asymp n$, Proposition~\ref{prop:cc:comparison}
shows that the same central-limit-scale replacement is valid whenever
\[
|c_{1,n}-c_1|=o(n^{-1/2}).
\]
In contrast, the ordinary convergence $c_{1,n}\to c_1$
alone need not give an $o(n^{-1/2})$ difference between
the correction factors. For example, the positive derivative
in \eqref{eq:cc:c_derivative} implies that aspect-ratio
deviations of order $n^{-1/4}$ generally produce
correction-factor deviations of that same order.
The constant limiting correction is therefore justified here
by the prescribed dimension sequence, not merely by its
first-order limit.

There is also a simple finite-aspect-ratio alternative:
\[
\mathfrak c_{1,n}^{\mathrm{MP}}
:=
\{1-\kappa_1f_1(\kappa_1,c_{1,n})\}^{-1}.
\]
By \eqref{eq:cc:finite_to_MP} and the same denominator bound,
\[
\mathfrak c_{1,n}^{\mathrm{MP}}-\mathfrak c_{1,n}
=
O(p^{-1}).
\]

 Thus this explicit correction is equivalent to the
finite-system correction at the central-limit scale in the
Gaussian setting even without a rate for $c_{1,n}\to c_1$.
The constant $\mathfrak c_1$ used in the present experiments
is additionally equivalent because
$c_{1,n}-c_1=O(n^{-1})$.
This argument justifies the Gaussian experiments only; the
non-Gaussian experiments remain empirical robustness checks and
are not covered by this centering argument.

\section{Auxiliary numerical experiments}\label{sec:num_exp_extra}

\subsection{Robustness to non-Gaussian errors: pathwise mediation}

Figure~\ref{fig:pathwise_nongaussian_qq} presents the
corresponding QQ-plots for the standardized pathwise
mediation-effect estimators. Their broadly linear
behavior suggests that the Gaussian approximation may
remain useful under the moment-matched non-Gaussian
error configurations considered here.

\input{plots_univ/figure_pathwise_nongaussian.tex}

\subsection{Choice of the tuning parameter}
Both limiting variances $\sigma^2$ and $\sigma_j^2$, defined in \eqref{lim_var_defn} and \eqref{lim_var_sigma_j}, depend nontrivially on the scaled regularization parameter $\kappa_1$. As a result, optimizing the limiting variance with respect to $\kappa_1$ does not appear to yield a useful closed-form expression.

A natural practical strategy is to select $\kappa_1$ by minimizing an estimate of the limiting variance over a grid of candidate values. For example, one may estimate $\sigma^2(\kappa_1)$ or $\sigma_j^2(\kappa_1)$ using an independent auxiliary sample, choose the minimizer over a grid, and then construct the final estimator on the primary sample. A systematic analysis of this data-driven tuning rule is left for future work.

\begin{figure}[tb]
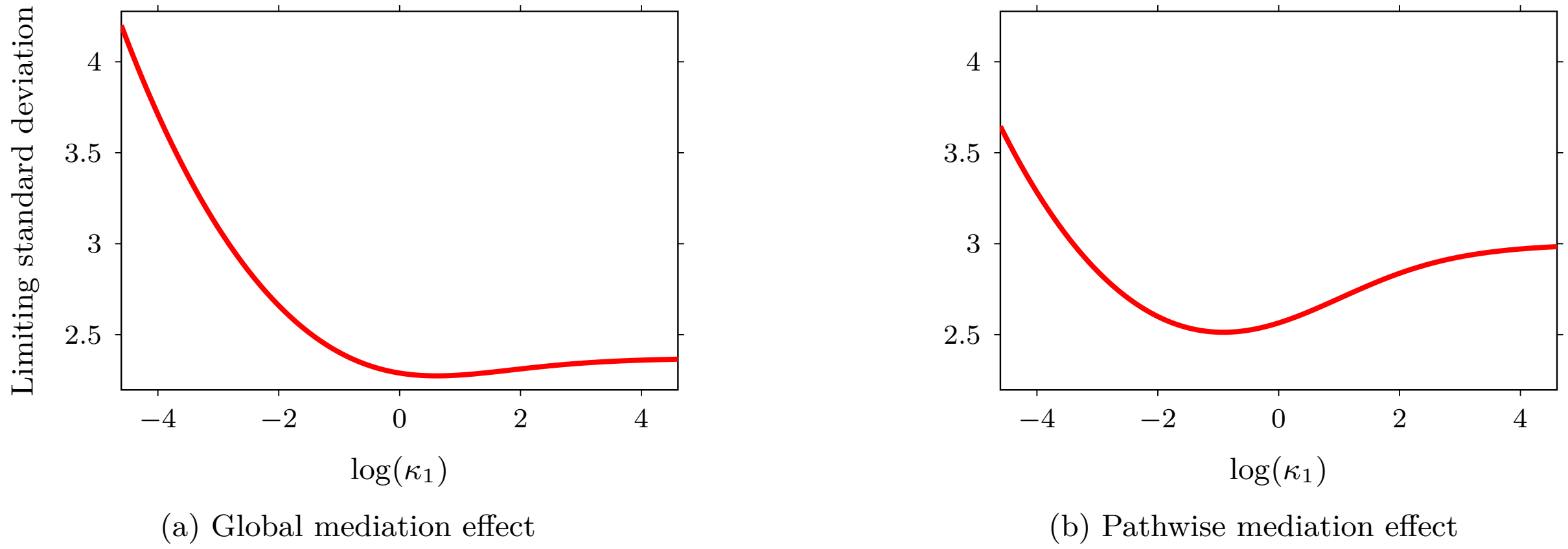

\centering
\captionsetup[subfigure]{font=footnotesize,skip=3pt}
\begin{subfigure}[t]{0.44\textwidth}
\centering
\input{plots/sd_vs_kap.tex}
\caption{Global mediation effect}
\label{sd_vs_kappa}
\end{subfigure}
\hfill
\begin{subfigure}[t]{0.44\textwidth}
\centering
\input{plots/sd_vs_kap_path.tex}
\caption{Pathwise mediation effect}
\label{sd_vs_kappa_path}
\end{subfigure}
\caption{Limiting standard deviations as functions of $\log(\kappa_1)$ with $c_1=1$ and $c_2=0.5$.}
\label{fig:sd_vs_kappa}
\end{figure}

Here we examine the shapes of the limiting standard deviations as
functions of $\kappa_1$. We fix $c_1=1$ and $c_2=0.5$. For the global
mediation effect, we use the signal configuration from
Section~\ref{glob_exp}. For the pathwise mediation effect, we use the
signal configuration from the pathwise experiment, in which
$\beta_{A,1}$ and $\theta_{M,1}$ remain of order one while the remaining
coordinates scale as $p^{-1/2}$. Figure~\ref{fig:sd_vs_kappa} displays
the resulting global and pathwise limiting standard deviations as
functions of $\log(\kappa_1)$. The functions $\sigma(\kappa_1)$ and $\sigma_j(\kappa_1)$ are
the limiting standard deviations of the $\sqrt{m}$-normalized
cross-fitted estimators. The corresponding plug-in limiting
standard deviations are $\sqrt{2}\sigma(\kappa_1)$ and
$\sqrt{2}\sigma_j(\kappa_1)$. This constant multiplier does not
change the minimizing value of $\kappa_1$.
% Generated by render_tuning_plots.R

\end{document}